\documentclass[10pt]{smfart} 
\usepackage{a4wide}

\usepackage{baskervald}
\usepackage[utf8]{inputenc} 
\usepackage[francais,english]{babel}
\usepackage{textcomp} 
\usepackage{stmaryrd}
\usepackage{graphicx} 
\usepackage{verbatim} 
\usepackage[retainorgcmds]{IEEEtrantools}
\usepackage[nice]{nicefrac}  
\usepackage{url}
\usepackage{url}
\usepackage{tikz}
\usetikzlibrary{matrix,arrows}
\usepackage{mathabx,epsfig}
\definecolor{cqcqcq}{rgb}{0.752941176471,0.752941176471,0.752941176471}
\definecolor{ffqqqq}{rgb}{1.,0.,0.}

\usepackage[all]{xy}
 \def\dar[#1]{\ar@<2pt>[#1]\ar@<-2pt>[#1]}
 \def\tar[#1]{\ar@<4pt>[#1]\ar@<0pt>[#1]\ar@<-4pt>[#1]}

\usepackage{amsthm} 
\usepackage{amsmath} 
\usepackage{amsfonts}
\usepackage{amssymb}
\usepackage{todonotes}

\renewcommand\epsilon\varepsilon 
\renewcommand\phi\varphi 
\def\quotient#1#2{
  \raise0ex\hbox{$#1$}\big/\!\lower1ex\hbox{$#2$}
}

\newcommand\NN{\mathbb{N}} 
\newcommand\ZZ{\mathbb{Z}} 
\newcommand\QQ{\mathbb{Q}} 
\newcommand\RR{\mathbb{R}} 
\newcommand\CC{\mathbb{C}} 
 
\newcommand\FF{\mathbb{F}} 
\newcommand\FP{\mathbb{F}_p}

\newcommand\fleche{\longrightarrow} 
\newcommand{\dans}{ \!\in\! }					
\newcommand{\priv}{ \!\smallsetminus\! }		
\newcommand{\defi}{ \begin{defin} }
\newcommand{\edefi}{ \end{defin} }
\newcommand{\exo}{ \begin{exercice} }		
\newcommand{\eexo}{ \end{exercice} }
\newcommand{\thr}{ \begin{theor} }
\newcommand{\ethr}{ \end{theor} }
\newcommand{\hyp}{ \begin{hypothese} }
\newcommand{\ehyp}{ \end{hypothese} }
\newcommand{\ques}{ \begin{question} }
\newcommand{\eques}{ \end{question} }
\newcommand{\nott}{ \begin{note} }
\newcommand{\enott}{ \end{note} }
\newcommand{\pro}{ \begin{prop} }
\newcommand{\epro}{ \end{prop} }
\newcommand{\prodefi}{ \begin{propdefi} }
\newcommand{\eprodefi}{ \end{propdefi} }
\newcommand{\rem}{ \begin{rema} }
\newcommand{\erem}{ \end{rema} }
\newcommand{\eppte}{ \end{propriete} }
\newcommand{\ppte}{ \begin{propriete} }
\newcommand{\cor}{ \begin{corr} }
\newcommand{\ecor}{ \end{corr} }
\newcommand{\lem}{ \begin{lemm} }
\newcommand{\elem}{ \end{lemm} }
\newcommand{\exen}{ \begin{exemple} }
\newcommand{\eexen}{ \end{exemple} }
\newcommand{\exe}{ \begin{exemple} }
\newcommand{\eexe}{ \end{exemple} }
\newcommand{\rap}{ \begin{rappel} }
\newcommand{\erap}{ \end{rappel} }
\newcommand{\que}{ \begin{quest} }
\newcommand{\eque}{ \end{quest} }
\newcommand{\fact}{ \begin{factt} }
\newcommand{\efact}{ \end{factt} }
\newcommand{\conj}{ \begin{conjecture} }
\newcommand{\econj}{ \end{conjecture} }
\newcommand{\dem}{ \begin{proof}[\textsc{Démonstration}] }  
\newcommand{\edem}{ \end{proof} } 

\DeclareMathOperator{\Ker}{Ker}
\DeclareMathOperator{\Card}{Card}

\DeclareMathOperator{\Coker}{Coker}

\DeclareMathOperator{\Isom}{Isom}

\DeclareMathOperator{\Aut}{Aut}

\DeclareMathOperator{\tr}{tr}

\DeclareMathOperator{\rg}{rg}
\DeclareMathOperator{\Ht}{ht}
\DeclareMathOperator{\Hom}{Hom}

\DeclareMathOperator{\Id}{Id}

\DeclareMathOperator{\Spec}{Spec}
\DeclareMathOperator{\Spf}{Spf}

\DeclareMathOperator{\id}{id}

\DeclareMathOperator{\Gal}{Gal}

\DeclareMathOperator{\Fil}{Fil}
\DeclareMathOperator{\im}{Im}

\DeclareMathOperator{\Gr}{Gr}
\DeclareMathOperator{\Ha}{Ha}
\DeclareMathOperator{\Frob}{Frob}
\DeclareMathOperator{\Lie}{Lie}
\DeclareMathOperator{\Hdg}{Hdg}
\DeclareMathOperator{\Fitt}{Fitt}
\DeclareMathOperator{\Newt}{\mathcal N\!ewt}

\DeclareMathOperator{\Def}{Def}
\DeclareMathOperator{\Cris}{Cris}
\DeclareMathOperator{\Diag}{Diag}
\DeclareMathOperator{\End}{End}

\DeclareMathOperator{\GL}{GL}

\theoremstyle{definition} 
\newtheorem{defin}{Définition}[section]
\newtheorem{propriete}[defin]{Propriété}
\newtheorem{hypothese}[defin]{Hypothèse}
\newtheorem{question}[defin]{Question}
\newtheorem{exercice}[defin]{Exercice}

\theoremstyle{plain} 
\newtheorem{theor}[defin]{Théoreme}
\newtheorem{lemm}[defin]{Lemme}  
\newtheorem{prop}[defin]{Proposition}
\newtheorem{propdefi}[defin]{Proposition-Définition}
\newtheorem{corr}[defin]{Corollaire}

\newtheorem{rappel}[defin]{Rappel}
\newtheorem{conjecture}[defin]{Conjecture}

\theoremstyle{remark} 
\newtheorem{rema}[defin]{Remarque}
\newtheorem{quest}[defin]{Question}

\newtheorem*{exemple*}{Exemple}
\newtheorem{exemple}[defin]{Exemple}
\newtheorem{note}[defin]{Note}
\newtheorem{factt}[defin]{Fait}

\title{ \textsc{Invariants de Hasse $\mu$-ordinaires}  }
\date{}
\author{Valentin Hernandez}
\address{Bureau 509, Tour 15-16\\4 Place Jussieu\\75005 Paris}
\email{valentin.hernandez@imj-prg.fr}
\urladdr{}

\begin{document}

\frontmatter
\altkeywords{}
\thanks{}
\maketitle

\selectlanguage{english}

\begin{abstract}
In this article, we construct in a purely local way partial (Hasse) invariants for $p$-divisible groups with given endomorphisms, using crystalline cohomology. Theses invariants generalises the classical Hasse invariant, and allow us to study families of such groups. We also study a few geometric properties of theses invariants.
Used in the context of Shimura varieties, for example, theses invariants are detecting some Newton strata, including the $\mu$-ordinary locus.
\end{abstract}

\selectlanguage{french}
\begin{abstract} Dans cet article on se propose de construire d'une manière purement locale des invariants partiels pour des groupes $p$-divisibles munis d'endomorphismes, 
en utilisant des résultats de cohomologie cristalline. Ces invariants généralisent l'invariant de Hasse, et permettent d'étudier des familles de tels groupes. On étudie aussi différentes 
propriétés géométriques de ces invariants. Appliqués – par exemple – à certaines variétés de Shimura, ces invariants détectent certaines strates de Newton, notamment la strate 
$\mu$-ordinaire. 
\end{abstract}

\tableofcontents
\section{Introduction}

\subsection{Chronologie}

L'introduction de l'invariant de Hasse remonte à 1936, et à la construction par Hasse et Witt, \cite{HaWi}, pour une courbe algébrique $C$ (non-singulière) 
de genre $g$ sur un corps fini de caractéristique $p$,
d'une matrice $H$ de taille $g\times g$, dite de Hasse-Witt.
En particulier pour les courbes elliptiques sur une base de caractéristique $p$, traité par Hasse un peu avant le cas général de Hasse-Witt, Katz, dans \cite{Ka}, 
donne alors la définition suivante, de ce qu'on appelle invariant de Hasse,

\defi[Katz]
Soit $S$ un schéma de caractéristique $p$, et $f : E\fleche S$ un schéma elliptique. Le Frobenius relatif $F : E \fleche E^{(p)}$ de $E/S$ induit une application,
\[ F^* : F^*(R^1f_*\mathcal O_E)) = (R^1f_*\mathcal O_E)^{\otimes p} \fleche R^1f_*\mathcal O_E,\]
qui en passant au dual induit une section $\widetilde{\Ha} \dans H^0(S,\omega_{E/S}^{\otimes(p-1)})$, appelée \textit{invariant de Hasse} de $E/S$.
\edefi

Cette définition s'étend simplement aux schémas abéliens et aux groupes $p$-divisibles sur un schéma de caractéristique $p$, en regardant le déterminant du Frobenius
sur l'algèbre de Lie du dual de Cartier. 
Une première application de l'invariant de Hasse apparait dans les travaux de Deligne et Serre \cite{DS} sur la construction des représentations galoisiennes associées aux 
formes modulaires de poids 1. Cette application utilise de manière cruciale la série d'Eisenstein $E_{p-1}$ (si $p \geq 5$) pour construire des congruences entre formes 
modulaires,en utilisant le résultat de Deligne selon lequel le développement de la série d'Eisenstein $E_{p-1}$ est congrue à 1 modulo $p$. Le fait est alors que $E_{p-1}$ 
relève l'invariant de Hasse.
En particulier, il y a un lien très fort entre l'invariant de Hasse et les congruences entre formes modulaires. Citons par exemple Hida, \cite{Hida}, qui sur une courbe 
modulaire $X$ sur $\ZZ_p$, utilise l'invariant de Hasse et le lieu ordinaire de $X \otimes \FP$ pour construire des espaces de formes modulaires $p$-adiques, qui réalisent les 
congruences (modulo des puissances de $p$) entre des formes modulaires.

Pour une variété abélienne sur un corps de caractéristique $p$, l'invariant de Hasse est non nul, 
exactement lorsque la variété abélienne est ordinaire, c'est à dire quand son $p$-rang est maximal. Si une variété est définie sur $\mathcal O_C$, $C = \widehat{\overline \QQ}_p$, et que son $p$-rang est maximal, il y a alors un sous-groupe dans la $p$-torsion (de type multplicatif) qui relève canoniquement la partie multiplicative (i.e. le noyau du Frobenius) modulo $p$.
On peut alors caractériser l'existence d'un tel sous-groupe, sous la condition que la valuation (tronquée) de l'invariant de Hasse (vu comme élément de $\mathcal O_C/p$) est plus 
petite que $\frac{1}{2}$. C'est un théorème de Lubin dans le cas d'une courbe elliptique, et de nombreux autres auteurs en dimension supérieure, en particulier Fargues 
(cf. \cite{Ka},\cite{Far}). Notons le lien avec les congruences entre formes modulaires, puisque le sous-groupe précédent permet de construire des opérateurs compacts sur 
des espaces de formes modulaires $p$-adiques (\cite{Hida},\cite{Ka},\cite{Pi}).

Le cas le plus intéressant d'un schéma $S$ sur lequel on dispose d'une variété Abélienne (ou un groupe $p$-divisible), est probablement celui où $S$ est une variété 
de Shimura (c'est le cas de \cite{Ka},\cite{Hida}). Il arrive néanmoins que pour de nombreuses variétés de Shimura, le lieu des points ordinaires (le lieu ordinaire) soit vide, ou 
de manière équivalente, l'invariant de Hasse soit nul en tout point. C'est par exemple le cas des variétés de Picard, pour $U(2,1)$, associé à un corps CM $E$, lorsque $p$ 
est inerte dans $E$.
Bien que le lieu ordinaire soit vide, il existe un ouvert, appelé lieu $\mu$-ordinaire, qui peut se définir à l'aide de la théorie des groupes (cf. \cite{RR}), ou dans le cas 
PEL, à l'ensemble des points de la variété de Shimura dont le groupe $p$-divisible a un polygone de Newton "minimal" (voir section \ref{sect3}). 
Lorsque la donnée de Shimura est PEL et non ramifiée, ce lieu est étudié par Wedhorn dans \cite{Wed}, qui prouve que c'est un ouvert dense, et qu'il coïncide avec le lieu 
ordinaire lorsque ce dernier est non vide.
Dans le cas contraire, que l'on appellera \textit{cas $\mu$-ordinaire}, pour une variété de Shimura de type unitaire, Goldring et Nicole ont construit dans \cite{GN} un invariant, dit invariant de Hasse $\mu$-ordinaire, qui détecte exactement le lieu $\mu$-ordinaire.

La stratégie de \cite{GN} est basée sur l'action du Frobenius sur la cohomologie de De Rahm modulo $p$ de la variété abélienne universelle, et (comme le lieu ordinaire est 
vide) de montrer qu'il existe une "division" de l'action de Frobenius par une bonne puissance de $p$, division qui n'est alors pas identiquement nulle sur $S$.  L'argument de 
\cite{GN} est alors le suivant, on regarde la variété universelle au dessus du lieu $\mu$-ordinaire, dans ce cas la cohomologie crystalline est relativement explicite, et on peut 
diviser une puissance du Frobenius, et étendre cette construction par densité à tout la variété grâce à un résultat de De Jong. Le même argument de densité est utilisé pour montrer que cette division passe au quotient à la cohomologie de De Rahm.
Le problème de l'argument de densité est qu'il rend la construction peu explicite hors du lieu ordinaire, et donc peu adaptée aux calculs que l'on espère faire. De plus, la construction est globale, alors que l'invariant de Hasse est un objet local (il ne dépend que du groupe $p$-divisible). Il y a aussi des travaux récents sur le sujet concernant d'autres invariants de Hasse (sur des strates plus générales que notre strate $\mu$-ordinaire), \cite{Bo,KW,GK}, dont l'un des invariants, celui associé à la strate $\mu$-ordinaire, est un produit des invariants partiels construit dans cet article. Mais ces constructions (bien que locales pour \cite{KW,GK}) ne sont pas aussi explicite que l'on espère, mais parfaitement adaptées à généraliser le travail de \cite{DS}. On propose alors une autre construction, sur laquelle on peut explicitement calculer les invariants de Hasse partiels, et leur produit l'invariant $\mu$-ordinaire, et en déduire des propriétés interressantes : le fait qu'ils soient réduits dans certains cas, et la compatibilité à la dualité. Dans un futur travail, \cite{Her2}, on utilisera cette construction explicite pour relier les invariants à la construction d'une filtration dite canonique dans les points de torsion d'un groupe $p$-divisible, dans un cadre plus général que celui \cite{Far}.

\subsection{Résultats}

Notre construction est dans l'esprit très proche de la construction originale de \cite{GN}, puisqu'elle est basée sur la cohomologie cristalline.
Plutôt que de se donner une variété de Shimura, on se donne $F/\QQ_p$ une extension non ramifiée, $\mathcal O$ son anneau d'entiers, $S$ un $\mathcal O$-schéma de caractéristique $p$, $G \fleche S$ un groupe de Barsotti-Tate ($p$-divisible) tronqué 
d'échelon $r$, et on suppose que $G$ est un $\mathcal O$-module, $\iota : \mathcal O \fleche \End_S(G)$. L'action de $\mathcal O$ sur $G$ donne une action de $\mathcal O$ sur son faisceau conormal, ainsi que celui de $G^D$, le dual de Cartier de $G$, et comme cet anneau est non ramifié, on peut
décomposer selon $\mathcal I = \Hom(F,\overline{\QQ_p})$,
\[ \omega_G = \bigoplus_{\tau \in \mathcal I} \omega_{G,\tau} \quad \text{et} \quad \omega_{G^D} = \bigoplus_{\tau \in \mathcal I} \omega_{G^D,\tau}.\]
Ces faisceaux sont localement libre, et quitte à prendre une union de composantes connexes de $S$, on peut supposer leur rang constant. 
On note alors $p_\tau = \rg_{\mathcal O_S} \omega_{G,\tau}$ et $q_\tau = \rg_{\mathcal O_S} \omega_{G^D,\tau}$, la signature de $(G,\iota)$.


On a alors le théorème suivant, qui rassemble plusieurs énoncés du corps de l'article : théorème \ref{thrdes}, théorème \ref{thrreduit}, théorème \ref{thrdual}.

\thr
Notons $f = [\mathcal O:\ZZ_p]$. Supposons $r$ assez grand relativement à la signature de $G$ (voir Théorème \ref{thrdes} et la remarque qui le précède).
Il existe une section \[\widetilde{^\mu\Ha}(G) \dans H^0(S,\det(\omega_{G^D})^{\otimes (p^{f}-1)}),\] telle que pour $s \dans S$, $\widetilde{^\mu\Ha}(G)_s$ 
est inversible si et seulement si
$G_s$ est $\mu$-ordinaire. La section $\widetilde{^\mu\Ha}(G)$ coïncide après multiplication par une puissance de $p$ explicite, dépendant de la signature de $G$, avec le 
déterminant du Veschiebung agissant sur le cristal de Berthelot-Breen-Messing de $G$. De plus on a les compatibilités suivantes  
\begin{enumerate}
\item Pour tout $\pi : S' \fleche S$, on a l'identification $\pi^*\widetilde{^\mu\Ha}(G) = \widetilde{^\mu\Ha}(\pi^*G)$,
\item $\widetilde{^\mu\Ha}(G) = \widetilde{^\mu\Ha}(G^D)$, où $G^D$ est le dual de Cartier de $G$.
\end{enumerate}
\ethr

La compatibilité à la dualité de l'invariant de Hasse "classique" avait été établie par Fargues dans \cite{Far}, proposition 2, voir aussi \cite{BijHa} dans un cas plus général, 
et on arrive, malgré quelques difficultés supplémentaires d'ordre principalement combinatoire, à établir cette compatibilité dans notre cas, cf. théorème \ref{thrdual}.

Le théorème précédent vaut en particulier lorsque $S$ est une variété de Shimura de type PEL.
Lorsque $S$ est une variété de Shimura PEL simple, comme dans la section \ref{sect9}, (essentiellement $p$ inerte dans le corps totalement réel), on a alors le théorème 
suivant,

\thr
L'invariant de Hasse $^\mu\Ha \dans H^0(S,(\det\omega)^{\otimes(p^f-1)})$ définit un diviseur de Cartier de $S$. De plus, lorsque tous les $p_\tau$ 
(sauf éventuellement $0,h$) sont distincts, ce diviseur de Cartier est réduit.
\ethr

Dans le cas de la courbe modulaire, l'invariant de Hasse est un diviseur de Cartier réduit c'est un résultat dû à Igusa, qui a été généralisé dans le cas 
des variétés de Siegel par Pilloni (\cite{Pi1}); Théorème A.4. Dans le cas d'autres variétés de Shimura dont le lieu ordinaire est non vide, l'invariant de Hasse peut être non 
réduit (par exemple le cas des variétés de Hilbert). Notre condition est nécessaire et suffisante, cf. remarque \ref{remred}.

\subsection{Construction} 

Essayons d'expliquer un peu plus en détails la construction du $\mu$-invariant de Hasse, ainsi que la démonstration de ses propriétés.
Etant donnée $(G,\iota)$ un groupe $p$-divisible tronqué d'échelon $r$ sur $S$ muni d'une action de $\mathcal O$, et signature $(p_\tau,q_\tau)$, on a les décompositions,
\[ \omega_G = \bigoplus_{\tau \in \mathcal I} \omega_{G,\tau} \quad \text{et} \quad \omega_{G^D} = \bigoplus_{\tau \in \mathcal I} \omega_{G^D,\tau}\]
dont les rangs des modules sont donnés par la signature. De plus, le Frobenius qui agit sur $G$, induit des morphismes,
\[ \omega_{G^D,\tau} \overset{V}{\fleche} \omega_{G^D\sigma^{-1}\tau}^{(p)}.\]
Ce morphisme ne peut être injectif si $q_{\sigma^{-1}\tau} < q_\tau$. L'obstruction principale à la non-vacuité du lieu ordinaire réside dans le fait que la signature n'est pas 
parallèle ($q_\tau = q_{\tau'}$ pour tout $\tau,\tau'$). Si la signature n'est pas parallèle, alors $V^f : \omega_{G^D,\tau} \fleche \omega_{G^D,\tau}^{(p^f)}$ n'est pas inversible, 
car un de ses composant n'est pas injectif, bien que ces deux faisceaux aient même rang.
L'idée, déjà présente dans \cite{GN}, est que si $S$ se relève sur $\ZZ_p^{nr}$ en un schéma sans torsion $\widetilde S$, sur lequel on a une variété abélienne $\widetilde A$ tel que 
$G = \widetilde A[p^r]\otimes_{\widetilde S} S$, on a peut alors regarder $\mathcal H^1_{dR}(\widetilde A/\widetilde S)$, qui se décompose selon les plongements $\tau$, et sur lequel on a un Verschiebung $\sigma^{-1}$-linéaire, qui est injectif. On peut alors considérer $V^f : \mathcal H^1_{dR}(\widetilde A/\widetilde S)_\tau \fleche \mathcal H^1_{dR}(\widetilde A/\widetilde S)_\tau^{(p^f)}$, et regarder la filtration de Hodge,
\[0 \fleche \omega_{\widetilde A^D/\widetilde S,\tau} \fleche \mathcal H^1_{dR}(\widetilde A/\widetilde S) \fleche \mathcal \omega_{\widetilde A/\widetilde S,\tau}^\vee\fleche 0\]
et $\omega_{\widetilde A^D/\widetilde S,\tau}$ se réduit modulo $p$ sur $\omega_{G^D,\tau}$. On peut alors espérer diviser par une bonne puissance de $p$ la puissance 
extérieure $q_\tau$-ième de $V^f$ agissant sur $\mathcal H^1_{dR}(\widetilde A/\widetilde S)_\tau$ qui est localement libre sur $\widetilde S$, lui même sans $p$-torsion, que cette division 
préserve (la puissance extérieur de) la filtration de Hodge (au moins modulo $p$), et soit inversible dessus si $G$ est $\mu$-ordinaire, puis réduire modulo $p$. On en déduira une "division de $V^f$", qui nous donnera un invariant partiel, $\widetilde{\Ha_\tau}$. Le produit de ces invariant sera l'invariant de Hasse $\mu$-ordinaire.

En général il ne semble pas possible de trouver un tel $\widetilde A$, mais on a un bon substitut de sa cohomologie de De Rham, le cristal de Bethelot-Breen-Messing 
$\mathbb D(G)$ de $G$. Le problème du site cristallin est qu'il semble difficile en général de maîtriser les différents ouverts mais lorsque $S$ est affine et lisse, on arrive à trouver "de bons ouverts" sur lesquels on arrive à diviser le Verschiebung. C'est pourquoi, dans un premier temps, on suppose que $S$ est lisse (c'est le cas d'une variété de Shimura PEL non ramifiée par exemple, cf. théorème \ref{thruni}).
Pour obtenir le cas général, remarquons que le champ des $\mathcal{BT}_r^\mathcal O$ (groupes de Barsotti-Tate tronqués d'échelon $r$ avec action de $\mathcal O$, dont on fixe la 
signature) est lisse (cf. \cite{Wed2}), il existe donc un schéma $X$, lisse, et une flèche lisse $X \fleche \mathcal{BT}_r^\mathcal O$, équivariante sous l'action de 
$\GL_{p^{rhf}}$. On peut alors construire les invariants de Hasse sur $X$ et les descendre à $\mathcal{BT}_r^\mathcal O$. On en déduit le cas d'une base générale.

Dans le cas ou $S = \Spec(\mathcal O_C/p)$, qui n'est pas lisse, on peut aussi construire les invariants de Hasse, de manière plus directe, en regardant le cristal sur $A_{cris}$. C'est le cas qui nous interressera dans un futur travail, cf. \cite{Her2}. 
On peut alors montrer que ce sont les mêmes que précédemment (essentiellement car $A_{cris}$ est sans torsion). On peut faire de même sur un espace de déformation d'un groupe $p$-divisible avec action de $\mathcal O$ (voir partie \ref{ssect83}, et section \ref{sect10}).

Le fait que le $\mu$-invariant $\widetilde{^\mu\Ha}$ est un diviseur de Cartier est essentiellement un théorème de Wedhorn \cite{Wed}. Pour montrer qu'il est réduit sur le champ $\mathcal{BT}_r$, on donne une donnée de Shimura PEL "simple" adaptée à $\mathcal O,(p_\tau,q_\tau)$ dans la section \ref{sect9}, on note $S$ la variété (sans niveau en $p$) associée.
On peut alors montrer que l'invariant $\mu$-ordinaire sur $S$ est réduit, en étudiant la géométrie des 
strates de Newton de $S$, dont on montre que sur celles qui nous intéressent, de codimension 1, l'invariant de Hasse $\mu$-ordinaire y définit une forme linéaire non nulle dans l'espace tangent. L'interêt est que sur ces strates, on peut utiliser la décomposition de Hodge-Newton rappelée en 
section \ref{sect2} ainsi qu'une forme explicite sur chaque anneau local de la déformation universelle en utilisant les displays, et donc de calculer effectivement l'invariant de Hasse. Au passage on montre aussi que ces strates sont des strates "minimales", c'est à dire que ce sont aussi des Strates d'Ekedahl-Oort.

La démonstration de la dualité repose sur les mêmes méthodes.

\medbreak
Explicitons les différentes sections. La section \ref{sect3} traite le cas d'un groupe $p$-divisible $G$ avec action de $\mathcal O$ sur un corps parfait, qui est extrêmement 
simple puisqu'à $G$ est associé son cristal de Dieudonné, qui est un $W(k)$-module libre. On rappelle alors les objets standards associés à $G$, les théorèmes connus sur 
les polygones (théorème de Mazur) et la décomposition Hodge-Newton. Dans une seconde partie, on explicite aussi dans ce cas la construction de nos invariants, à la fois
pour expliquer au lecteur la philosophie de la construction, mais aussi puisque c'est ce cas qui a servit de fil conducteur à la suite de l'article. La section \ref{sect4} traite le cas 
d'une base lisse, qui contient le corps residuel de $F$, auquel on peut se ramener par extension des scalaires. La section suivante redescend la construction à une base lisse 
quelconque. La section \ref{sect6} montre la compatibilité au changement de base (c'est essentiellement due à la même compatibilité pour le cristal de \cite{BBM}). Dans la 
section \ref{sect7}, on donne un exemple où les invariants ne sont pas compatibles au produit (c'est du à un problème combinatoire sur les signatures), et on donne une 
condition pour que cette compatibilité ait lieu. La section \ref{sect8} est dédiée à la descente des invariants sur le champ $\mathcal{BT}_r^\mathcal O$, qui permet donc de 
faire la construction sur une base quelconque. On y donne aussi des exemples de groupes $\mu$-ordinaires, ainsi que des calculs sur des schémas de base utiles 
pour la suite. Dans la section \ref{sect10}, on définit des données de Shimura simples, qui seront utilisées ensuite pour montrer (sous une hypothèse 
précise) que nos invariants sont réduits. Dans la section \ref{sect11}, on montre la compatibilité à la dualité. Enfin, dans la section \ref{sect12}, on donne une forme explicite 
pour les $\mathcal O$-modules $p$-divisibles sur $\mathcal O_C$. Les appendices sont des rappels bien connus de notations et définitions à propos des fibrés inversibles sur 
des champs, ainsi que des diviseurs de Cartier.

Remerciements : Je voudrai remercier Laurent Fargues et Vincent Pilloni pour m'avoir introduit à ce sujet et aussi pour leurs explications, leurs encouragements et leur aide au cours de l'écriture de ce travail. Je voudrai aussi remercier Stephane Bijakowski pour de nombreux conseils et de nombreuses discussions.

\section{Notations et décompositions}
\label{sect2}
Soit $p$ un nombre premier, $S$ un schéma de caractéristique $p$ (i.e. tel que $p\mathcal O_{S} = 0$), $F/\QQ_p$ 
une extension non ramifiée, de degré $f$, et on note $p^f$ le cardinal du corps résiduel 
$\kappa_F$ et $\mathcal O = \mathcal O_F$ son anneau des entiers. On note aussi $\widebar{\FP}$ une cloture algébrique de $\FP$. 
On suppose que $S$ est défini au-dessus de $\kappa_F$.
Soit $G$ un groupe de Barsotti-Tate tronqué d'échelon $r$ sur $S$, muni d'une action de 
$\mathcal O$, 
\[ \iota : \mathcal O \fleche \End_{S}(G).\]
On note pour la suite $\Sigma_n = \Spec(W_n(\kappa_F))$,
$\Sigma = \Spec(W(\kappa_F))$, et $\underline{G}$ le faisceau sur $\Cris(S/\Sigma)$ induit par $G$.
On note 
\[\mathcal E := \mathcal{E}xt^1_{S/\Sigma}(\underline{G^D},\mathcal O_{S/\Sigma}),\]
le cristal de Dieudonné (covariant) de $G$, au sens de Berthelot-Breen-Messing (\cite{BBM}), qui est munie de deux applications, 
\[ F : \mathcal E^{(p)} \fleche \mathcal E \quad \text{et} \quad V : \mathcal E \fleche \mathcal E^{(p)},\]
qui sont induites respectivement par le Vershiebung et le Frobenius de $G$. C'est un cristal en 
$\mathcal O_{S/\Sigma}$-modules localement libres.

Comme $G$ est muni d'une action de $\mathcal O$, $\mathcal E$ l'est aussi. Notons,
\[\mathcal I = \Hom_{\ZZ_p}(\mathcal O,W_n(\kappa_F)) \simeq \Hom_\ZZ(\kappa_F,\overline{\FP}).\]

Comme $S$ est un $\kappa_F$-schéma, on a un isomorphisme,
\[ \mathcal O \otimes_{\ZZ_p} \mathcal O_{S} \simeq \prod_{\tau \in \mathcal I} \mathcal O_{S}.\]

\lem
Il existe une décomposition canonique de $\mathcal E$ en sous-cristaux,
\[\mathcal E = \bigoplus_{\tau \in \mathcal I} \mathcal E_{\tau},\]
telle que,
\[ V : \mathcal E_\tau \fleche \mathcal E_{\sigma^{-1}\tau}^{(p)}\quad \text{et} \quad
F : \mathcal E_\tau^{(p)} \fleche \mathcal E_{\sigma\tau},\]
où on note $\mathcal E_\tau^{(p)} = (\mathcal E_\tau)^{(p)}$.
De même, on peut décomposer canoniquement les 
$\mathcal O_S$-modules,
\[ \omega_G = \bigoplus_{\tau \in \mathcal I} \omega_{G,\tau} \quad \text{et} \quad 
\omega_{G^D} = \bigoplus_{\tau \in \mathcal I} \omega_{G^D,\tau}.\]
\elem

\defi
Les modules $\omega_{G^D,\tau}$ et $\omega_{G,\tau}$, pour tout $\tau \in \mathcal I$, sont localement libres. Supposons que les dimensions de chacun des $\omega_{G,\tau},\omega_{G^D,\tau}, \tau \dans \mathcal I$, soient constantes sur la base $S$. On définit alors la signature de 
$G$ par la donnée $(p_\tau,q_\tau)$, où,
\[ p_\tau = \rg_{\mathcal O_{S}} \omega_{G,\tau} \quad \text{et} \quad 
q_\tau = \rg_{\mathcal O_{S}} \omega_{G^D,\tau}.\]
On a de plus, \[p_\tau + q_\tau = \frac{\Ht G}{f} =: \Ht_\mathcal O(G) = h,\]
la $\mathcal O$-hauteur de $G$.
\edefi

On va utiliser les décompositions précédentes pour construire des invariants de Hasse partiels, associés aux plongements $\tau \dans \mathcal I$ (qui tiendront donc compte de l'action
de $\mathcal O$ sur $G$).

\section{Le cas d'un corps parfait}
\label{sect3}

On rappelle ici les constructions classiques de polygones associées à des cristaux, en restant le plus élémentaire possible. On rappelle aussi le travail de Mantovan-Viehmann
sur la décomposition de Hodge-Newton (dans le cas simple d'un corps parfait) et on fait le lien avec la signature. Enfin, on donne une construction explicite de 
l'invariant $\mu$-ordinaire (plus exactement des invariants partiels associés aux plongements $\tau$) qui servira de fil conducteur dans la section suivante. Mis à part peut-être
cette dernière partie, les résultats de cette section sont bien connus, mais plutôt que de renvoyer à de nombreux articles, et pour fixer les notations, on a préféré les réexposer 
ici.

\subsection{Cristaux et polygones.}

Supposons ici que $S = \Spec k$, où $k$ est un corps parfait de caractéristique $p$ tel que $\kappa_F \subset k$ (mais on ne fixe pas de tel plongement). 
Supposons de plus que $G/S$ est un groupe de Barsotti-Tate (non tronqué). Plutôt que de regarder son cristal de Dieudonné $\mathcal E$ au sens de Berthelot-Breen-Messing, \cite{BBM}, regardons son module de Dieudonné covariant au sens de Fontaine, \cite{Fon},
\[\mathbb D(G) = \Hom(G^D,CW).\]

C'est un $V$-cristal au sens suivant,
\defi
Soit $m \dans \NN^*$. Un $V^m$-cristal sur $k$ est un couple $(M,V)$, où $M$ est un $W(k)$-module libre de type fini, et 
\[V : M \fleche M,\]
est un morphisme injectif, $\Frob^{-m}$-linéaire. On note $Cris^m/k$ la catégorie des $V^{m}$-cristaux sur $k$. 
Un $V^m$-isocristal est un couple $(M,V)$ où $M$ est un $W(k)[\tfrac{1}{p}]$-espace vectoriel et $V$ est un isomorphisme $\Frob^{-m}$-linéaire. En particulier, si
$(M,V)$ est un $V^m$-cristal, $(M[\tfrac{1}{p}],V)$ est un $V^m$-isocristal.
\edefi

Notons $p-Div/k$ la catégorie des groupes p-divisibles sur $k$. Notons $Cris^{1}_{[0,1]}/k$ la catégorie des $V$-cristaux tels que $pM \subset VM \subset M$.
D'après \cite{Fon}, Chapitre III, on a le théorème suivant :
\thr
Le foncteur,
\[ \mathbb D : 
\begin{array}{ccc}
 p-Div/k & \fleche  & Cris^{1}_{[0,1]}/k  \\
 G & \rightsquigarrow  &\mathbb D(G)   
\end{array}
\]
est une équivalence de catégories telle que $\rg \mathbb D(G) = \Ht G$ et $\quad [\mathbb D(G) : V(\mathbb D(G))] = \dim G$.
On a aussi,
\[ \Lie(G) = \omega_{G}^\vee \simeq\quotient{\mathbb D(G)}{V(\mathbb D(G))} \quad \text{et} \quad \omega_{G^D} \simeq \quotient{V(\mathbb D(G))}{p\mathbb D(G)}.\]
\ethr

Le cristal de Dieudonné au sens de Fontaine est alors relié au cristal au sens de Berthelot-Breen-Messing par,

\pro[\cite{BBM} chapitre 4, Théorème 4.2.14 ]
\label{proBBMFON}
Soit $G$ un groupe $p$-divisible tronqué d'échelon $r$ sur $k$, un corps parfait de caractéristique $p$. Notons $
\mathbb D(G)$ son module de Dieudonné covariant et $\mathcal E = \mathcal Ext^1(G^D,\mathcal O_{\Spec(k)/\Spec(W(k))})$ son cristal de Dieudonné, et $\sigma$ le Frobenius de $k$. On a alors,
\[ \mathcal E_{(W(k) \twoheadrightarrow k)} \simeq \mathbb D(G)^{(\sigma)}.\]
\epro

On a supposé que $G$ était muni d'une action de $\mathcal O$, on a donc une action de $\mathcal O$ sur $\mathbb D(G)$. Décomposons alors,
\[ \mathbb D(G) = M  = \bigoplus_{\tau \in \mathcal I} M_\tau, \quad \omega_G = \bigoplus_{\tau \in \mathcal I}\omega_{G,\tau} \quad \text{et} \quad \omega_{G^D} = \bigoplus_{\tau \in \mathcal I}\omega_{G^D,\tau}. \]

Afin d'alléger les notations, on fixe un isomorphisme $\mathcal I \simeq \ZZ/f\ZZ = \{1,\dots,f\}$ (on notera donc désormais $i$ pour un élément de $\mathcal I$), tel que, \[\forall i \dans \mathcal I, \quad \Frob \circ i = i +1.\]
On a alors,
\[V : M_i \fleche M_{i-1} \quad \text{et}\quad F : M_i \fleche M_{i+1}.\]
Notons,
\[ \phi_i = M_i \overset{V}{\fleche} M_{i-1} \overset{V}{\fleche} M_{i-2} \overset{V}{\fleche} 
\cdots\overset{V}{\fleche} M_{i+1} \overset{V}{\fleche} M_{i}.\]
C'est à dire $\phi_i = V^f_{|M_i}.$ 
On a alors que $(M_i, \phi_i)$ est un $V^f$-cristal, pour tout $i \dans \mathcal I$. De plus,
$F$ est une isogénie de $\Frob^*(M_i,\phi_i)$ vers $(M_{i+1},\phi_{i+1})$, on en déduit donc 
que leurs isocristaux sont égaux, et donc que le polygone de Newton de $(M_i,\phi_i)$ ne dépend
pas de $i \dans \mathcal I$. 

\defi
On appelle polygone de Newton de $G$ (ou de $\mathbb D(G)$) avec $\mathcal O$-action (ou encore $\mathcal O$-polygone de Newton), que l'on note 
$\Newt_{\mathcal O}(G)$ (ou $\Newt_{\mathcal O}(\mathbb D(G))$), la fonction convexe affine par morceaux,
\[\Newt_{\mathcal O}(G) = \frac{1}{f} \Newt(M_i,\phi_i),\]
où $\Newt(M_i,\phi_i)$ est la fonction affine par morceaux défini grâce aux pentes de la classification de Dieudonné-Manin des isocristaux.
\edefi

\defi
Soit $M$ un $W(k)$-module de longueur finie. On appelle polygone de Hodge renversé de $M$, noté 
$\Hdg^\diamond(M)$, le polygone concave à abscisses de ruptures entières (que l'on identifie à une fonction affine par morceaux sur $[0,+\infty[$), tel que,
\[\Hdg^\diamond(M)(i) = \deg M - v(\Fitt_i M), \forall i \dans \NN.\]
\edefi

Par exemple, si 
\[M \simeq \bigoplus_{i=1}^n \quotient{W(k)}{p^{a_i}W(k)},\quad a_1 \geq \dots \geq a_n,\]
Alors $\Hdg^\diamond(M)$ commence en $(0,0)$ et est de pente $a_i$ sur $[i-1,i]$, et de pente 0 sur $[n,+\infty[$.

\pro
Soit $M_1,M_2,M_3$ trois $W(k)$-modules libres de rang fini, et $u,v$ deux morphismes 
(semi-linéaires vis-à-vis de $\sigma \dans\Aut(k))$,
\[ M_1 \overset{u}{\fleche} M_2  \overset{v}{\fleche} M_3,\]
qui deviennent des isomorphismes après inversion de $p$ (c'est à dire $u,v$ injectifs et de conoyau de longueur finie). Alors,
\[\Hdg^\diamond(\Coker(v \circ u)) \leq \Hdg^\diamond(\Coker u) + \Hdg^\diamond(\Coker v).\]
\epro

\dem
Choisissons des bases de $M_1,M_2$ et $M_3$. Soit $i \dans \NN$. Les coefficients des matrices de $\bigwedge^i u$, $\bigwedge^i v$ et $\bigwedge^i(v\circ u)$ dans ces bases engendrent les idéaux $\Fitt_i(\Coker u), \Fitt_i(\Coker v)$ et $\Fitt_i(\Coker v\circ u)$. Mais comme $\bigwedge^i(v\circ u) = \bigwedge^i v \circ \bigwedge^i u$, on a,
\[ \Hdg^\diamond(\Coker(v \circ u))(i) \leq \Hdg^\diamond(\Coker u)(i) + \Hdg^\diamond(\Coker v)(i).\qedhere\]
\edem

\defi
Soit $(M,V)$ un $V^m$-cristal sur $k$ de rang $h$, alors on définit son polygone de Hodge renversé par,
\[\Hdg^\diamond(M,V) = \Hdg^\diamond(M/V(M)).\] 
On définit ensuite son polygone de Hodge $\Hdg(M,V)$ comme le renversé de $\Hdg^{\diamond}(M,V)_{|[0,h]}$, voir remarque qui suit. 
\edefi

\rem
Notons que le polygone de Hodge de $(M,V)$ ne dépend pas uniquement de $M/V(M)$, contrairement à $\Hdg^\diamond(M,V)$ : on veut contrôler la pente 0 (qui correspond à la partie étale de $G$), ce qui explique l'introduction de $h$ dans la définition. C'est à dire, si les pentes (avec multiplicités, et éventuellement nulles) de $\Hdg^\diamond(M,V) = \Hdg^\diamond(M/V(M))$ sont 
$b_1\geq b_2 \geq \dots \geq b_r \geq\underbrace{0 \geq\dots\geq 0}_{h-r \text{ fois}}$, 
celles de $\Hdg(M)$ sont $0 \leq\dots \leq 0 \leq b_r \leq b_{r-1} \leq \dots \leq b_1$.

Cette normalisation est faite pour que l'on ne perde pas les pentes nulles lors du renversement de $\Hdg(M,V)$.
\erem

Revenons à $\mathbb D(G)$. On a la décomposition,
$\mathbb D(G) = \bigoplus_{i \in \mathcal I} \mathbb D(G)_i$,
et, par la proposition précédente,
\[\Hdg^\diamond(\mathbb D(G)_i,\phi_i) := \Hdg^\diamond(\mathbb D(G)_i/V^f(\mathbb D(G)_i)) \leq \sum_{i \in \mathcal I} \Hdg^\diamond(\mathbb D(G)_i/V(\mathbb D(G)_{i+1})).\]
De plus, on sait d'après l'inégalité de Mazur, cf \cite{Kat} Theorem 1.4.1, que,
\begin{equation}
\label{eqnewthdg}
 \Newt(\mathbb D(G)_i[\tfrac{1}{p}],\phi_i[\tfrac{1}{p}]) \geq \Hdg(\mathbb D(G)_i,\phi_i) \geq 
\sum_{i \in \mathcal I} \Hdg_i(\mathbb D(G),V),\end{equation}
où l'on note $\Hdg_i(\mathbb D(G),V)$ le renversé de $\Hdg^\diamond(\mathbb D(G)_i/V(\mathbb D(G)_{i+1}))_{|[0,h]}$.

\rem
Bien sûr, en general $\Hdg_i(\mathbb D(G),V) \neq \Hdg(\mathbb D(G)_i,V^f)$.
\erem

On en déduit donc que,
\[ \Newt_{\mathcal O}(\mathbb D(G)) \geq \frac{1}{f} 
\sum_{i \in \mathcal I} \Hdg_i(\mathbb D(G),V).\]
 
 \defi
Soit  $M$ un cristal sur $k\supset \kappa_F$ avec action de $\mathcal O = \mathcal O_F$. 
On appelle $\mathcal O$-polygone de Hodge de $M$, le polygone,
\[ \Hdg_{\mathcal O}(M) = \frac{1}{f} 
\sum_{i \in \mathcal I} \Hdg_i(M,V).\]
Si $G$ est un groupe $p$-divisible sur $k\supset \kappa_F$ avec action de $\mathcal O = \mathcal O_F$, on appelle $\mathcal O$-polygone de Hodge de $G$, noté $\Hdg_{\mathcal O}(G)$, celui de $\mathbb D(G)$.
 \edefi

On retrouve donc la proposition classique suivante, voir aussi \cite{RR},

\pro
\label{proNewHdg}
Le polygone de Newton $\Newt_\mathcal O(G)$ est au-dessus du polygone de Hodge $\Hdg_\mathcal O(G)$, et ils ont mêmes points terminaux.
\epro

 \begin{figure}[h]
 \caption{Exemple de polygones de Hodge et Newton dans le cas de $\mathcal O = \ZZ_{p^2}$}
 \label{figNewtHdg}
 \begin{center}
 \begin{tikzpicture}[line cap=round,line join=round,>=triangle 45,x=1cm,y=1cm]
\draw[->,color=black] (-0.5,0.) -- (5.,0.);
\foreach \x in {,1.,2.,3.,4.}
\draw[shift={(\x,0)},color=black] (0pt,2pt) -- (0pt,-2pt);
\draw[->,color=black] (0.,-0.5) -- (0.,4.);
\foreach \y in {-0.5,0.5,1.,1.5,2.,2.5,3.,3.5}
\draw[shift={(0,\y)},color=black] (2pt,0pt) -- (-2pt,0pt);
\clip(-0.5,-0.5) rectangle (5.,4.);
\draw (0.,0.)-- (1.,0.);
\draw (1.,0.)-- (3.,1.);
\draw (3.,1.)-- (4.00163934426,2.20659971306);
\draw [color=ffqqqq] (0.,0.)-- (2.44169908507,0.720849542536);
\draw [color=ffqqqq] (2.44169908507,0.720849542536)-- (4.00163934426,2.20659971306);
\draw (0.209836065574,2.55093256815)-- (0.583606557377,2.55093256815);
\draw [color=ffqqqq] (0.209836065574,2.02725968436)-- (0.588524590164,2.02008608321);
\draw [dash pattern=on 1pt off 1pt] (3.,1.)-- (3.,0.);
\draw (2.99344262295,-0.0961262553802) node[anchor=north west] {$q_{\tau_2}$};
\draw (1.03606557377,-0.117647058824) node[anchor=north west] {$q_{\tau_1}$};
\draw (0.64262295082,2.15638450502) node[anchor=north west] {$\Newt_{\mathcal O}$};
\draw (0.637704918033,2.7230989957) node[anchor=north west] {$\Hdg_{\mathcal O}$};
\draw [dash pattern=on 1pt off 1pt] (4.00163934426,2.20659971306)-- (4.,0.);
\draw (3.96229508197,-0.0602582496413) node[anchor=north west] {$h$};
\end{tikzpicture}
\end{center}
\end{figure}
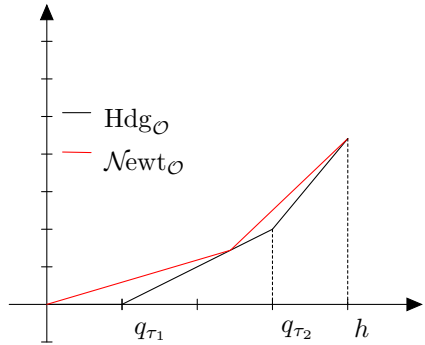
  
  On rappelle la définition suivante, cf.\cite{Wed} par exemple

\defi
Soit $G/\Spec(k)$ un $\mathcal O$-module $p$-divisible, où $k$ est un corps algébriquement clos de caractéristique $p >0$. Il est dit $\mu$-ordinaire si,
\[ \Newt_{\mathcal O}(G) = \Hdg_{\mathcal O}(G).\]
Soit $G/S$ une famille de $\mathcal O$-module $p$-divisible sur une base $S$ telle que $p\mathcal O_S = 0$, alors $G$ est dit $\mu$-ordinaire si pour tout
point géométrique $x = \Spec(k) \fleche S$, $G_x$ est $\mu$-ordinaire.
\edefi

\subsection{Rappels sur la décomposition de Hodge-Newton}

On rappelle dans cette section le résultat principal de \cite{MV}. Supposons maintenant que $\Newt_{\mathcal O}(M)$ et $\Hdg_{\mathcal O}(M)$ se touchent en un point d'abscisse $x_0$, et que cette abscisse soit une abscisse de rupture 
pour $\Newt_{\mathcal O}(M)$ (ce qui est par exemple le cas sur la figure \ref{figNewtHdg}). D'après les inégalités (\ref{eqnewthdg}) 
il en est de même pour $\Newt(M_i,\phi_i)$ et 
$\Hdg(M_i,\phi_i), \forall i \in I$. On en déduit donc par le théorème de Katz sur la rupture Hodge-Newton (cf. \cite{Katz} Theorem 1.6.1), que,
\[ M_i = M_i' \oplus M_i'', \quad \forall i\dans \mathcal I,\]
où $M_i',M_i''$ sont stables par $V^f$. Notons, 
\[N_i = M_i[\tfrac{1}{p}], \quad N_i' = M_i'[\tfrac{1}{p}], \quad N_i'' = M_i''[\tfrac{1}{p}].\]
Ce sont des isocristaux, et les pentes de $N_i'$ (respectivement $N_i''$) correspondent à celles de $N_i$ avant (respectivement après) le point de rupture d'abscisse $x_0$, et donc
les pentes de Newton de $N_i'$ sont toutes strictement inférieures à celles de $N_i''$ .
Mais l'isogénie,
\[ V : M_i[\tfrac{1}{p}] \fleche M_{i-1}[\tfrac{1}{p}],\]
envoie composantes isoclines sur composantes isoclines de mêmes pentes, donc,
\[ V(N_i' \cap M_i) \subset N_{i-1}' \cap M_{i-1} \qquad V(N_i'' \cap M_i) \subset N_{i-1}'' \cap M_{i-1}\]
et donc, $V(M_i') \subset M_{i-1}'$ et $V(M_i'') \subset M_{i-1}''$.
Notons alors, 
\[M' = \bigoplus_{i \in \mathcal I} M_i' \quad \text{et} \quad M'' = \bigoplus_{i \in \mathcal I} M_i''.\]
Munis de $V$, ce sont des $F$-cristaux sur $k$.

\lem
Soit $M', M''$ deux $V$-cristaux sur $k$ de rang $r,s$ munis d'une action de $\mathcal O$. Alors pour tout $i \dans \mathcal I$,
\[ \Hdg_i(M'\oplus M'',V) = \Hdg_i(M',V) \star \Hdg_i(M'',V),\]
où $\star$ désigne la \textit{concaténation} des polygones, i.e. réordonne les pentes par ordre croissant.
\elem

\dem
Il suffit de remarquer que $V(M'_{i-1} \oplus M_{i-1}'') = VM_{i-1}' \oplus VM_{i-1}''$, et donc on peut écrire,
\[M'_i/VM'_{i-1}= \bigoplus_{i=1}^r W(k)/p^{a_i}W(k) \quad \text{et} \quad M''_i/VM''_{i-1}= \bigoplus_{i=1}^s W(k)/p^{b_i}W(k),\]
et \[(M'_i \oplus M_i'')/V(M'_{i-1} \oplus M_{i-1}'') = \bigoplus_{i=1}^r W(k)/p^{a_i}W(k) \oplus \bigoplus_{j=1}^s W(k)/p^{b_i}W(k).\]
de telle sorte que $r= \rg M'$ et $s= \rg M''$ (i.e. on fait apparaitre les pentes 0).
Le polygone de Hodge associé à $i$ de $M'\oplus M''$ est alors donné (par ses pentes) en réordonnant de manière croissante l'ensemble $\{a_i,b_j : i,j\}$.
\edem

En particulier, $\Hdg_i(M' \oplus M'')(r) \leq \Hdg_i(M')(r)$, avec égalité si et seulement si les polygones coïncident sur $[0,r]$.

On en déduit en particulier le théorème suivant, voir aussi \cite{MV},

\thr[Mantovan-Viehmann]
\label{thrHN}
Soit $(M,V)$ un cristal sur $k$ muni d'une action de $\mathcal O$. 
Alors si $\Newt_{\mathcal O}(M)$ et $\Hdg_{\mathcal O}(M)$ se touchent en un point d'abscisse $x$ qui est un point de rupture pour $\Newt_{\mathcal O}(M)$, alors il existe des sous-cristaux $M'$ et $M''$ munis de la $\mathcal O$-structure induite de $M$ tels que,
\[ M = M' \oplus M''.\]
De plus les pentes de Newton (respectivement de Hodge) de $M'$ (respectivement de $M''$) correspondent à celles de $M$ avant (respectivement après) le point d'abscisse $x$, et 
$\rg M' = x$.
\ethr

\dem
L'assertion sur la décomposition et les polygones de Newton découle de précédemment. 
Pour les polygones de Hodge, le lemme précédent nous assure que pour tout $i$, $\Hdg_i(M,V)$ a ses pentes constituées de celles de $\Hdg_i(M',V)$ et $\Hdg_i(M'',V)$, malheureusement ce n'est plus vrai pour $\Hdg_\mathcal O(M,V)$ (voir un contre exemple en section \ref{sect7}).
Néanmoins comme $\Newt_\mathcal O(M)$ et $\Hdg_\mathcal O(M)$ ont même points terminaux pour tout $M$, en appliquant cela à $M'$ on en déduit que
$\sum_i \Hdg_i(M)(x) = \sum_i \Hdg_i(M')(x)$, or comme pour tout $i$ le lemme précédent nous assure que $\Hdg_i(M)(x) \leq \Hdg_i(M')(x)$ on en déduit l'égalité pour tout $i$, et donc l'égalité des polygones pour tout $i$, sur $[0,x]$. En particulier, $\Hdg_\mathcal O(M)_{[0,x]} = \Hdg_\mathcal O(M')$, ainsi que la proposition associée pour $M''$.
%
\edem

On en déduit directement que si $G$ est un groupe de Barsotti-Tate sur $k$, muni d'une action de
 $\mathcal O$, alors si son $\mathcal O$-polygone de Newton et son 
 $\mathcal O$-polygone de 
 Hodge se touchent en un point d'abscisse $x \dans \ZZ$, qui est un point de rupture 
 du polygone de Newton, alors on a une décomposition,
\[ G = G_1 \times G_2,\]
où $G_1,G_2$ sont deux groupes de Barsotti-Tate sur $k$, munis d'une action de 
$\mathcal O$, tels que \[ \Ht G_1 = x.\]

\rem
La réciproque est bien sûr fausse : voir l'exemple \ref{exeprod}. C'est à la base de quelques complications quant à la compatibilité des invariants de Hasse au produit.
\erem

\begin{figure}[h]
\caption{$\mathcal O$-polygone de Hodge associée à la signature $(q_\tau)_{\tau\in \mathcal I}$.}
\label{figHdg}
\begin{center}
\begin{tikzpicture}[line cap=round,line join=round,>=triangle 45,x=1.0cm,y=1.0cm]
\draw[->,color=black] (-0.5,0) -- (10,0);
\foreach \x in {,1,2,3,4,5,6,7,8,9}
\draw[shift={(\x,0)},color=black] (0pt,-2pt);
\draw[->,color=black] (0,-1) -- (0,4.5);
\foreach \y in {-1,1,2,3,4}
\draw[shift={(0,\y)},color=black] (2pt,0pt) -- (-2pt,0pt);
\clip(-0.5,-1) rectangle (10,4.5);
\draw (0,0)-- (1.86,0);
\draw (1.86,0)-- (3.72,0.22);
\draw (3.72,0.22)-- (5,0.62);
\draw [dash pattern=on 1pt off 1pt] (5.82,1.08)-- (6.7,1.52);
\draw (7.64,2)-- (8.42,2.66);
\draw (8.42,2.66)-- (9,3.64);
\draw [dash pattern=on 1pt off 1pt] (3.72,0)-- (3.72,0.22);
\draw [dash pattern=on 1pt off 1pt] (5,0.62)-- (5,0);
\draw [dash pattern=on 1pt off 1pt] (7.64,2)-- (7.62,0);
\draw [dash pattern=on 1pt off 1pt] (8.42,2.66)-- (8.42,0);
\draw (1.82,-0.1) node[anchor=north west] {$q^{(1)}$};
\draw (3.64,-0.1) node[anchor=north west] {$q^{(2)}$};
\draw (4.8,-0.1) node[anchor=north west] {$q^{(3)}$};
\draw (7.2,-0.1) node[anchor=north west] {$q^{(r-1)}$};
\draw (8.2,-0.1) node[anchor=north west] {$q^{(r)}$};
\draw [dash pattern=on 1pt off 1pt] (9,3.64)-- (9,0);
\draw (8.8,-0.2) node[anchor=north west] {$h$};
\draw (0.74,0.5) node[anchor=north west] {0};
\draw (2.82,0.88) node[anchor=north west] {$\frac{n_1}{f}$};
\draw (3.8,1.1) node[anchor=north west] {$\frac{n_1+n_2}{f}$};
\draw (6.60,2.9) node[anchor=north west] {$\frac{n_1 + \dots + n_{r-1}}{f}$};
\draw (8.40,3.68) node[anchor=north west] {$1$};
\end{tikzpicture}
\end{center}
\end{figure}
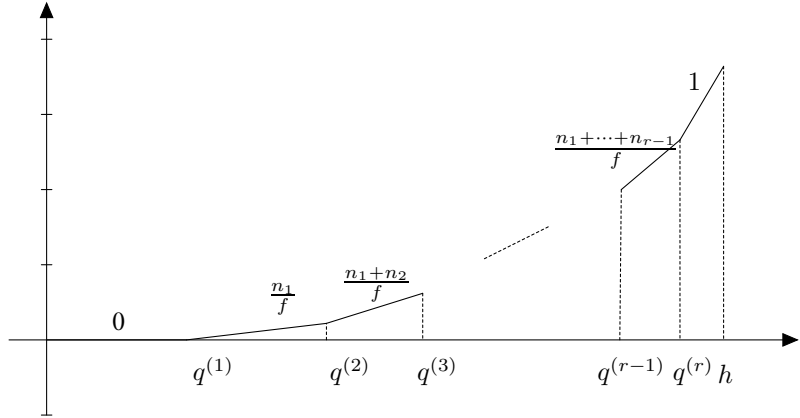

\rem
\label{rem29}
Supposons que nous sommes sous les hypothèses du théorème \ref{thrHN}, pour $M = \mathbb D(G)$ associé à un groupe $p$-divisible (donc $M$ est un $V$-cristal). Essayons de décrire les signatures de $M' = \mathbb D(G')$ et $M''= \mathbb D(G'')$ en fonction de celle de $M$.
Réordonnons les abscisses de ruptures du polygone de Hodge (moyenné) de $M$, $\Hdg_{\mathcal{O}_F}(M)$,
\[ \{ q_\tau : \tau \in \mathcal I\} =: \{ q^{(1)} < q^{(2)} < \dots < q^{(r)}\},\]
et pour tout $i \dans \{ 1,\dots,r\},$ notons $n_i$ la multiplicité de l'abscisse $q^{(i)}$, c'est à dire,
\[ n_i := \Card\{\tau \in \mathcal I : q_\tau = q^{(i)}\}.\]
Le $\mathcal O$-polygone de Hodge de $M$ est alors le polygone sur 
$[0,\Ht_{\mathcal O}(M)]$ défini par les abscisses de ruptures $q^{(i)}, i \dans \{1,\dots,r\}$, et dont 
la pente entre entre $q^{(i)}$ et $q^{(i+1)}$, notée $\nu_i$, est,
\[ \nu_i = \frac{n_1 + \dots + n_i}{f}.\]
Voir figure \ref{figHdg}. Si $q^{(1)} >0$ alors la pente sur $[0,q^{(1)}]$ est nulle, et si $q^{(r)} < \Ht_{\mathcal O}(M)$, la pente sur
$[q^{(r)},\Ht_{\mathcal O}(M)]$ est 1.

On sait par ailleurs que le polygone de Hodge de $M'$ est celui de $M$ entre les abscisses 
$[0,x]$. On sait aussi que,
\[ \omega_{M} = \omega_{M'} \oplus \omega_{M''}.\]
On en déduit les égalités suivantes, où on note $(q_\tau')$ et $(q_\tau'')$ les signatures de $M'$ et $M''$,
\[ q_\tau = q_\tau' + q_\tau'', \quad \{q_\tau' : \tau \in \mathcal I\} =\{q_\tau : \tau \in \mathcal I, q_\tau \leq x\},\]
et, pour tout $i$ tel que  $q^{(i)} < x$, \[\Card\{ \tau : q_\tau' = q^{(i)}\} = \Card\{  \tau : q_\tau = q^{(i)}\} = n_i.\]
On en déduit donc simplement par récurrence sur $i \dans \{1,\dots,r\}$ que,
\[ \forall \tau \in \mathcal I, \quad 
q_\tau' = 
\left\{
\begin{array}{ccc}
q_\tau  &  \text{si } q_\tau \leq x \\
 x &  \text{sinon}
\end{array}
\right.
, \quad 
q_\tau'' = 
\left\{
\begin{array}{ccc}
0 &  \text{si } q_\tau \leq x \\
q_\tau - x &  \text{sinon}
\end{array}
\right.
\]
\erem

\subsection{Construction d'invariants dans le cas d'un point} 
\label{sect22}

Le but de cet article, est de construire des invariants qui vont "mesurer" le défaut d'avoir un point de contact 
(entre les polygones de Hodge et Newton) en une 
abscisse donnée. On commence par présenter dans cette section la construction dans le cas particulier, à la fois plus simple et plus parlant, d'un corps parfait $k$ de caractéristique $p$.
Soit $(M,V) = \mathbb D(G)$ le cristal au sens de Fontaine d'un $k$-groupe $p$-divisible $G$, muni d'une action de $\mathcal O$, et on suppose $\kappa_F \subset k$.
On rappelle que l'on a noté $\mathcal I = \Hom(\kappa_F,k) \simeq \{1,\dots,f\}$ sur lequel on a "un ordre" induit par le Frobenius, et que, 
\[ M \simeq \bigoplus_{i \in \mathcal I} M_i, \quad V : M_i \fleche M_{i-1}.\]
On a aussi $f = \Card(\mathcal I)$, et des entiers $h,p_i,q_i$, pour $i \dans \mathcal I$ et $h = p_i + q_i$.

\rem
Afin de faire la lien avec la section suivante, remarquons que l'on a l'égalité,
\[VM_{i+1} = \Ker \left( M_i \fleche M_i/V(M_{i+1}) \simeq  \omega_{G,i}^\vee\right) \subset M_i.\]
\erem

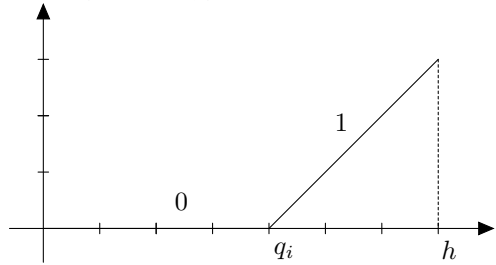
\begin{figure}[h]
\caption{Polygone $\Hdg(M_i/VM_{i-1})$}
\label{figHdgi}
\centering
\begin{tikzpicture}[line cap=round,line join=round,>=triangle 45,x=1.5cm,y=1.5cm]
\draw[->,color=black] (-0.3,0.) -- (4.,0.);
\foreach \x in {,0.5,1.,1.5,2.,2.5,3.,3.5}
\draw[shift={(\x,0)},color=black] (0pt,2pt) -- (0pt,-2pt);
\draw[->,color=black] (0.,-0.3) -- (0.,2.);
\foreach \y in {,0.5,1.,1.5}
\draw[shift={(0,\y)},color=black] (2pt,0pt) -- (-2pt,0pt);
\clip(-0.3,-0.3) rectangle (4.,2.);
\draw (0.,0.)-- (2.,0.);
\draw (2.,0.)-- (3.5,1.5);
\draw [dash pattern=on 1pt off 1pt] (3.5,1.5)-- (3.5,0.);
\draw (1.9606557377,-0.0459110473458) node[anchor=north west] {$q_i$};
\draw (3.44672131148,-0.0243902439024) node[anchor=north west] {$h$};
\draw (2.5,1.1) node[anchor=north west] {$1$};
\draw (1.08278688525,0.4) node[anchor=north west] {$0$};
\end{tikzpicture}
\end{figure}

Considérons $\Hdg_{\mathcal O}(M)$ le $\mathcal O$-polygone de Hodge. Ses abscisses
de ruptures sont, \[\{q_i : i \dans \mathcal I, q_i \not\in \{0,f\}\} \subset ]0,n[ \cap \NN.\]
En effet, chaque $\Hdg_i(M,V)$ a pour seul point de rupture $q_i$, voir figure \ref{figHdgi}.

Choisissons $q_i, i \dans \mathcal I$, une abscisse de rupture, et posons
\[k_{i} = \sum_{j : q_j \leq q_i} (q_i - q_j).\]
L'ordonnée sur $\Hdg_{\mathcal O}(M)$ du point de rupture d'abscisse $q_i$ est $k_i/f$.
En effet, 
\begin{IEEEeqnarray*}{cclcc}
\Hdg_{\mathcal O}(M)(q_i) &= & \frac{1}{f} \sum_{j \in \mathcal I} \Hdg(M_j/VM_{j-1})(q_i) & &  \\
& = & \frac{1}{f} \left( \sum_{q_j \leq q_i} \underbrace{\Hdg(M_j/VM_{j-1})(q_i)}_{q_i-q_j} + \sum_{q_j > q_i}  
\underbrace{\Hdg(M_j/VM_{j-1})(q_i)}_{0}\right) & =&  \frac{k_i}{f}.
\end{IEEEeqnarray*}

On a alors le lemme suivant,

\lem
\label{lem24}
Si $q_i < q_{i_0}$, alors comme réseaux dans $\bigwedge^{q_{i_0}} M_i[\frac{1}{p}]$,
\[ \bigwedge^{q_{i_0}} VM_{i+1} \subset p^{q_{i_0} - q_i} \bigwedge^{q_{i_0}} M_i.\]
\elem

Cela découle du lemme d'algèbre linéaire suivant (voir la démonstration de la proposition \ref{pro34}),
\lem
\label{fact1}
Soit $A$ un anneau, et $M$ un $A$-module, muni d'une filtration à un cran, c'est-à-dire un sous module $N \subset M$. On a une filtration sur $\bigwedge^{n} M$, compatible à la précédente, donnée par,
\[ \forall 0 < i \leq n, \Fil^i\left(\bigwedge^{n}M\right) = \im\left(\bigwedge^{i}N \otimes \bigwedge^{n-i}M\fleche \bigwedge^{n}M\right).\]
Il existe alors une application surjective, 
\[\bigwedge^{i}N \otimes \bigwedge^{n-i}\quotient{M}{N} \fleche \Gr^i\left(\bigwedge^{n}M\right).\]
De plus, c'est un isomorphisme si $M$ est projectif de type fini et $N$ est localement facteur 
direct.
\elem

\dem 
On vérifie directement que l'application surjective,
\[\bigwedge^{i}N \otimes \bigwedge^{n-i}M \fleche \Gr^i\left(\bigwedge^{n}M\right),\]
se factorise par,
\[\bigwedge^{i}N \otimes \bigwedge^{n-i}\quotient{M}{N} \fleche \Gr^i\left(\bigwedge^{n}M\right).\]
Si de plus (localement) $M$ est libre, et $N$ un sous-module facteur direct, $M = N \oplus N'$, alors
\[\Gr^i\left(\bigwedge^{n}M\right) = \im(\bigwedge^i N \otimes \bigwedge^{n-i} N' \fleche \bigwedge^n M) \simeq \bigwedge^i N \otimes \bigwedge^{n-i} N'.\qedhere\]
\edem

\cor
L'application,
\[ V^f : \bigwedge^{q_{i_0}} M_{i_0} \fleche \bigwedge^{q_{i_0}} VM_{i_0+1},\]
est divisible par $p^{k_{i_0}}$, c'est à dire se factorise par, 
\[ p^{k_{i_0}} \bigwedge^{q_{i_0}} VM_{i_0+1}.\]
Comme $VM_{i_0+1}$ est libre sur $W(k)$, qui est sans $p$-torsion,  on en déduit une factorisation,
\begin{center}
\begin{tikzpicture}[description/.style={fill=white,inner sep=2pt}] 
\matrix (m) [matrix of math nodes, row sep=3em, column sep=2.5em, text height=1.5ex, text depth=0.25ex] at (0,0)
{ 
 \bigwedge^{q_{i_0}} M_{i_0} & & \bigwedge^{q_{i_0}} VM_{i_0+1} \\
& \bigwedge^{q_{i_0}} VM_{i_0+1} & \\
 };

\path[->,font=\scriptsize] 
(m-1-1) edge node[auto] {$V^f$} (m-1-3)
(m-1-1) edge node[auto] {$\zeta_0$} (m-2-2)
(m-2-2) edge node[auto] {$p^{k_{i_0}}$} (m-1-3);
\end{tikzpicture}
\end{center}
\ecor

On vérifie de plus que l'application $\zeta_0$ envoie bien,
\[\Ker \left( \bigwedge^{q_{i_0}} VM_{i_0+1} \fleche \bigwedge^{q_{i_0}} \quotient{VM_{i_0+1}}{pM_{i_0}} \simeq \bigwedge^{q_{i_0}} \omega_{G^D,i_0} \right),\]
dans lui même.
On en déduit une application $\Frob^{-f}$-linéaire,
\[\widetilde{\Ha}_{i_0} : \bigwedge^{q_{i_0}} \omega_{G^D,i_0} \fleche \bigwedge^{q_{i_0}} 
\omega_{G^D,i_0}.\]
C'est-à-dire une application linéaire,
\[\widetilde{\Ha}_{i_0} : \bigwedge^{q_{i_0}} \omega_{G^D,i_0} \fleche \bigwedge^{q_{i_0}} 
\omega_{G^D,i_0}^{(p^f)} \simeq (\det\omega_{G^D,i_0})^{\otimes p^f},\]

\defi
On appelle \textit{Invariant de Hasse partiel de $G$ en $i_0$} la section,
\[\widetilde{\Ha}_{i_0}(G) \dans \det\left(\omega_{G^D,i_0}\right)^{\otimes (p^f-1)},\]
associée au morphisme précédent.
\edefi

L'invariant de Hasse ainsi défini mesure bien le défaut de contact des polygones de Hodge et 
Newton en l'abscisse $q_{i_0}$, au sens de la proposition suivante.

\pro
\label{propolygone}
Soit $G$ un groupe $p$-divisible comme précédemment, et $i \dans \mathcal I$. Les conditions
 suivantes sont alors équivalentes,
\begin{enumerate}
\item Les polygones $\Newt_{\mathcal O}(G)$ et $\Hdg_{\mathcal O}(G)$ ont un point 
de contact en l'abscisse $q_i$.
\item L'invariant de Hasse partiel, $\widetilde\Ha_i(G)$, est inversible.
\end{enumerate}
\epro

\dem
Notons $K = W(k)[\frac{1}{p}]$. On a \[\Newt_{\mathcal O}(G) = \frac{1}{f} \Newt(M_i \otimes K,V^f).\]
D'après Katz, cf \cite{Kat} 1.3.4, $\Newt_{\mathcal O}(G)$ a un point de rupture en l'abscisse $q_i$ qui touche le polygone $\Hdg_{\mathcal O}(G)$, si et seulement si, l'isocristal,
\[(\bigwedge^{q_i} M_i \otimes K, p^{-k_i}V^f),\]
a une pente 0 avec multiplicité non nulle.
Notons \[\Lambda = (\bigwedge^{q_i} VM_{i-1}, p^{-k_i}V^f),\]
qui est un réseau (un cristal) dans $(\bigwedge^{q_i} M_i \otimes K, p^{-k_i}V^f)$, et écrivons,
\[ \Lambda = \Lambda^{et} \oplus \Lambda^{nilp}.\] Il suffit de vérifier que $\Lambda^{et} \neq 0$ si et seulement si $\widetilde{\Ha}_i$ est inversible. On  a la suite exacte,
\[ 0 \fleche \im\left( pM_i \otimes \bigwedge^{q_i-1} VM_{i-1} \fleche \bigwedge^{q_i}VM_{i-1}\right)
\fleche \underbrace{\bigwedge^{q_i} VM_{i-1}}_{\Lambda} \overset{\pi}{\fleche} \bigwedge^{q_i} \quotient{VM_{i-1}}{pM_i} \fleche 0.\]
Notons $N = \Ker \pi$ le premier module. Le module quotient est de dimension 1 sur $k$, et le morphisme semi-linéaire $p^{-k_i}V^f = \phi$ qui agit dessus n'est autre que 
$\widetilde\Ha_i$. Or on sait que (par exemple par la remarque \ref{remfort}, mais on peut le déduire facilement du lemme \ref{lem24}),
\[ \phi(N) \subset p\Lambda,\]
et donc $N \subset \Lambda^{nilp}$, et $\pi(\Lambda^{et})$ est un sous-module de \[\bigwedge^{q_i} \quotient{VM_{i-1}}{pM_i},\]
sur lequel $\phi = \widetilde\Ha_i$ est inversible. Or ce dernier espace est un $k$-espace vectoriel de dimension 1, donc $\widetilde\Ha_i$ est inversible si et seulement si $\Lambda^{et} \neq 0$.
\edem

\section{Construction dans le cas d'une base lisse}
\label{sect4}

Soit comme dans l'introduction $\mathcal O = \mathcal O_F$ non ramifié, et $G/S$ un $\mathcal O$-module $p$-divisible tronqué d'échelon $r$, on suppose cette fois de plus que $S$ est un schéma 
lisse sur $\Spec(\kappa_F)$, de caractéristique $p$, tel que les dimensions des modules localement libres $\omega_{G,\tau}$ et $\omega_{G^D,\tau}$, pour tout 
$\tau \dans \mathcal I$, soient constantes sur $S$. 

On se place dans le topos cristallin (Zariski) $(S/\Sigma)_{Cris}$, où $\Sigma = \Spec(W(\kappa_F))$. Soit 
$\mathcal E = \mathcal Ext^1(\underline{G},\mathcal{O}_{S/\Sigma})$ le cristal de Dieudonné de $G$, c'est un $\mathcal O_{S/\Sigma}/p^r$-module localement libre de rang fini (cf. \cite{BBM} Théorème 3.3.10 ).

\defi
Les faisceaux cristallins $\mathcal O_{S/\Sigma}$ et $J_{S/\Sigma}$ vérifient, pour tout objet $(U,T,\gamma)$ de $Cris(S/\Sigma)$,
\[ \mathcal O_{S/\Sigma}(U,T,\gamma) = \mathcal O_T, \quad \text{et} \quad J_{S/\Sigma}(U,T,\gamma) = \Ker(\mathcal O_T \twoheadrightarrow O_U).\]
On a un morphisme de topos $i_{S/\Sigma} : S_{Zar} \fleche (S/\Sigma)_{Cris}$ où $S_{Zar}$ est le topos Zariski de $S$, défini par,
pour tout $\mathcal F \dans (S/\Sigma)_{cris}$, pour tout $\mathcal G \dans S_{Zar}$ et pour tous ouverts $U \subset S, (U,T,\gamma) \dans Cris(S/\Sigma)$,
\[ i_{S/\Sigma}^*\mathcal F(U) = \mathcal F(U,U,0) \quad \text{et} \quad i_{S/\Sigma,*}\mathcal G(U,T,\gamma) = \mathcal G(U).\]
\edefi 

On a alors la suite exacte,
\[ 0 \fleche J_{S/\Sigma} \fleche \mathcal O_{S/\Sigma} \fleche i_{S/\Sigma*} \mathcal O_S \fleche 0.\]
On obtient donc une suite exacte,
\[0 \fleche J_{S/\Sigma}\mathcal E \fleche \mathcal E \fleche i_{S/\Sigma*} 
\left(i_{S/\Sigma}^*\mathcal E\right) \fleche 0.\]
Comme $\mathcal E$ est un cristal, on a \[i_{S/\Sigma*} 
\left(i_{S/\Sigma}^*\mathcal E\right) \simeq \mathcal E/J_{S/\Sigma}\mathcal E.\]
Mais maintenant, on a une filtration (cf \cite{BBM} Corolaire 3.3.5) de $\mathcal O_{S}$-modules,
\[ 0 \fleche \omega_{G^D} \fleche \mathcal E_{(S\overset{\id}{\rightarrow} S)} 
\fleche \omega_{G}^\vee\fleche 0.\]
On en déduit donc une filtration (cf \cite{BBM}, Corollaire 3.3.5),
\[ 0 \fleche i_{S/\Sigma*}\omega_{G^D} \fleche 
i_{S/\Sigma*}\left(i_{S/\Sigma}^*\mathcal E\right) \fleche 
i_{S/\Sigma*}\omega_{G}^\vee\fleche 0.\]

\defi
\label{def32}
On définit la filtration de $\mathcal E$ comme étant le sous-faisceau sur le site cristallin, 
\[ \Fil \mathcal E = \Ker(\mathcal E \fleche \mathcal E/J_{S/\Sigma}\mathcal E 
\fleche i_{S/\Sigma*}\omega_{G}^\vee).\]
\edefi

Cette filtration vérifie que,
\[\left(\mathcal E/\Fil \mathcal E\right)_{(S\rightarrow S)} = \omega_G^\vee \quad \text{et} \quad 
\left(\Fil \mathcal E\right)_{(S\rightarrow S)} = \omega_{G^D}.\]

On a aussi grâce à l'action de $\mathcal O$ sur $G$, des filtrations induites, 
\[\Fil \mathcal E_\tau \subset \mathcal E_\tau\]
qui vérifient,
\[\left(\mathcal E_\tau/\Fil \mathcal E_\tau\right)_{(S\rightarrow S)} = \omega_{G,\tau}^\vee 
\quad \text{et} \quad 
\left(\Fil \mathcal E_\tau\right)_{(S\rightarrow S)} = \omega_{G^D,\tau}.\]

On va construire les invariants de Hasse partiels comme sections de certains fibrés en droites, de manière analogue à la construction de l'invariant de Hasse classique d'un groupe $p$-divisible, telle qu'elle est par exemple décrite dans \cite{Far}, section 2.2. On renvoie également à \cite{GN} pour une construction analogue dans le cas $\mu$-ordinaires sur les variétés de Shimura.
On s'intéresse donc à l'application $\bigwedge^{q_\tau} V^f$, que l'on note plutôt $V^f$,
\[ V^f : \bigwedge^{q_\tau} \mathcal E_\tau \fleche \bigwedge^{q_\tau} \mathcal E_\tau^{(p^f)}.\]
Dans le but de relier cela à la construction décrite dans la section \ref{sect22}, introduisons 
les sous-faisceaux cristallins,
\[ \Fil\left(\bigwedge^{q_\tau} \mathcal E_\tau\right) = \im\left(\bigwedge^{q_\tau}\Fil\mathcal E_\tau \fleche 
\bigwedge^{q_\tau}\mathcal E_\tau\right),\] 
\[\Fil\left(\bigwedge^{q_\tau} \mathcal E_\tau^{(p^f)}\right) = \im\left(\Fil\left(\bigwedge^{q_\tau}\mathcal E_\tau\right)^{(p^f)} \fleche 
\bigwedge^{q_\tau}\mathcal E_\tau^{(p^f)}\right).\]

\lem
\label{pro2}
Soit $\tau \dans \mathcal I$. L'application \[V : \bigwedge^{q_{\tau_0}} \mathcal E_\tau \fleche \bigwedge^{q_{\tau_0}} \mathcal E_{\sigma^{-1}\tau}^{(p)},\] se factorise en,
\begin{center}
\begin{tikzpicture}[description/.style={fill=white,inner sep=2pt}] 
\matrix (m) [matrix of math nodes, row sep=3em, column sep=2.5em, text height=1.5ex, text depth=0.25ex] at (0,0)
{ 
\bigwedge^{q_{\tau_0}} \mathcal E_\tau & & \bigwedge^{q_{\tau_0}} \mathcal 
E_{\sigma^{-1}\tau}^{(p)} \\
& \im\left(\bigwedge^{q_{\tau_0}} (\Fil\mathcal E_{\sigma^{-1}\tau})^{(p)} \fleche \bigwedge^{q_{\tau_0}} (\mathcal E_{\sigma^{-1}\tau})^{(p)} \right)  & \\
 };

\path[->,font=\scriptsize] 
(m-1-1) edge node[auto] {$V$} (m-1-3)
(m-1-1) edge node[auto,left] {} (m-2-2)
(m-2-2) edge node[auto,right] {} (m-1-3);
\end{tikzpicture}
\end{center}
\elem
%
%

\defi
Soit $S,R$ deux anneaux et $S \twoheadrightarrow R$ un épaississement à puissances divisés. Supposons que $pR = 0$.
Soit $\phi : S \fleche S$ un morphisme. On dit que $\phi$ est un \textit{relèvement fort} de Frobenius si $\phi(s) \equiv s^p \pmod p$.
En particulier $\phi$ relève (faiblement) le Frobenius de $R$, au sens où le diagramme suivant commute,
\begin{center}
\begin{tikzpicture}[description/.style={fill=white,inner sep=2pt}] 
\matrix (m) [matrix of math nodes, row sep=3em, column sep=2.5em, text height=1.5ex, text depth=0.25ex] at (0,0)
{ 
S & &S \\
R& &R\\
 };

\path[->,font=\scriptsize] 
(m-1-1) edge node[auto] {$\phi$} (m-1-3)
(m-1-1) edge node[auto,left] {} (m-2-1)
(m-1-3) edge node[auto,left] {} (m-2-3)
(m-2-1) edge node[auto] {$\Frob_R$} (m-2-3);
\end{tikzpicture}
\end{center}
\edefi

\exe
Si $R = k$ un corps parfait de caractéristique $p$, et $S = W(k)$, d'ideal $J_{S/R} = (p)$, le Frobenius $\sigma$ de $S$ est un relèvement fort de Frobenius.

Si $C = \widehat{\overline{\QQ_p}}$, $R = \mathcal O_C/p$ et $S = A_{cris}(\mathcal O_C)$, $\theta : A_{cris} \fleche \mathcal O_C/p$ d'idéal $\ker \theta$, qui n'est pas $p$-adique, mais le Frobenius $\phi$ de $A_{cris}$ est un relèvement fort de Frobenius.

Si $R = k$ un corps parfait de caractéristique $p$, et $S = k[T]/(T^p)$ muni de ses puissances divisées canoniques, alors $\phi : S \fleche S$ tel que
$\phi(x) = x^p$ pour $x \dans k$ et $\phi(T) = T$, n'est pas un relèvement fort de Frobenius sur $S$.
\eexe

À partir de maintenant on fixe un $\tau_0 \dans \mathcal I$.

\pro
\label{pro34}
Si $q_{\tau} < q_{\tau_0}$, alors,
\[\im\left(\bigwedge^{q_{\tau_0}} \Fil\mathcal E_{\tau} \fleche \bigwedge^{q_{\tau_0}} \mathcal E_{\tau} \right) \subset 
J_{S/\Sigma}^{q_{\tau_0} - q_\tau}\bigwedge^{q_{\tau_0}} \mathcal E_{\tau}. \]
Si de plus $(U,T,\gamma)$ est un ouvert qui est un épaississement $p$-adique, ou bien tel qu'il existe $\phi : T\fleche T$ un relèvement fort de Frobenius, alors,
\[\im\left(\bigwedge^{q_{\tau_0}} (\Fil\mathcal E_{\tau})^{(p)}_{(U,T,\gamma)} \fleche \left(\bigwedge^{q_{\tau_0}} \mathcal E_{\tau}\right)_{(U,T,\gamma)}^{(p)} \right) \subset p^{q_{\tau_0} - q_\tau}\left(\bigwedge^{q_{\tau_0}} \mathcal E_{\tau}\right)_{(U,T,\gamma)}^{(p)}. \]
\epro

\dem
On utilise le fait précédent, \ref{fact1}, qui se transpose aux faisceaux cristallins :
soit $(U \hookrightarrow T) \dans Cris(S/\Sigma)$, alors on a un morphisme surjectif,
\[ \bigwedge^i J_{S/\Sigma}\mathcal {E_\tau}_{(U\hookrightarrow T)} \otimes_{\mathcal O_T} 
\bigwedge^{n-i} \Fil \mathcal {E_\tau}_{(U\hookrightarrow T)}/J_{S/\Sigma}\mathcal {E_\tau}_{(U\hookrightarrow T)} \twoheadrightarrow \Gr^i\left( \bigwedge^n 
(\Fil \mathcal E_\tau)_{(U\hookrightarrow T)}\right)\]
par le fait \ref{fact1}. Or le foncteur faisceau associé est un foncteur exact, \cite{SGA4}, 4.4.1, donc, dans $\mathcal O_{S/\Sigma}$-Mod, on a un épimorphisme,
\[ \bigwedge^i J_{S/\Sigma}\mathcal {E_\tau} \otimes_{\mathcal O_{S/\Sigma}}
\bigwedge^{n-i} \Fil \mathcal {E_\tau}/J_{S/\Sigma}\mathcal {E_\tau} \twoheadrightarrow \Gr^i\left( \bigwedge^n (\Fil \mathcal E_\tau)\right).\]
Appliquons alors cela avec $n = q_{\tau_0}$ et $0 \leq i < q_{\tau_0} - q_\tau$ (i.e. $q_{\tau_0} \geq q_{\tau_0} - i > q_\tau$).
On a l'isomorphisme, \[\Fil \mathcal E_\tau/J_{S/\Sigma}\mathcal E_\tau \simeq i_{S/\Sigma*}\omega_{G^D,\tau}.\]
On en déduit donc,
\[ \bigwedge^{q_{\tau_0}-i} \Fil \mathcal E_\tau/J_{S/\Sigma}\mathcal E_\tau  \simeq \bigwedge^{q_{\tau_0}-i}i_{S/\Sigma*}\omega_{G^D,\tau} 
\simeq i_{S/\Sigma*}\bigwedge^{q_{\tau_0}-i}\omega_{G^D,\tau} = 0,\]
car $q_{\tau_0} - i > q_\tau = \dim_{\mathcal O_{S}} \omega_{G^D,\tau}$.
Tous les gradués considérés sont donc nuls, et on en déduit,
\[ \bigwedge^{q_{\tau_0}} \Fil \mathcal E_\tau = \im\left( \bigwedge^{q_{\tau_0}-q_\tau} J_{S/\Sigma}\mathcal {E_\tau} \otimes_{\mathcal O_{S/\Sigma}} 
\bigwedge^{q_\tau} \Fil \mathcal {E_\tau} \fleche \bigwedge^{q_{\tau_0}} \Fil \mathcal E_\tau\right)\]
et donc, en prenant l'image dans $\bigwedge^{q_{\tau_0}} \mathcal E_\tau$,
\[ \im\left( \bigwedge^{q_{\tau_0}} \Fil \mathcal E_\tau \fleche \bigwedge^{q_{\tau_0}} \mathcal E_\tau\right) = 
 \im\left( \bigwedge^{q_{\tau_0}-q_\tau} J_{S/\Sigma}\mathcal {E_\tau} \otimes_{\mathcal O_{S/\Sigma}} 
\bigwedge^{q_\tau} \Fil \mathcal {E_\tau} \fleche\bigwedge^{q_{\tau_0}} \mathcal E_\tau\right) 
\subset J_{S/\Sigma}^{q_{\tau_0} - q_\tau} \bigwedge^{q_{\tau_0}} \mathcal E_\tau.\]
Soit $(U,T,\gamma)$ un ouvert de $Cris(S/\Sigma)$. Si l'épaississement est $p$-adique, alors $J_{S/\Sigma}(U,T,\gamma)= p\mathcal O_{T}$ 
et donc la seconde assertion découle de précédemment (même sans tordre par le Frobenius de $S$). 
Supposons donc qu'il existe $\phi$ sur $\mathcal O_T$ qui relève fortement le Frobenius, c'est-à-dire que pour tout 
$x \dans \mathcal O_T$, $\phi(x) \equiv x^p \pmod {p\mathcal O_T}$. Or si $x \dans J_{S/\Sigma}(U,T,\gamma)$, alors $x^p = p!\gamma_p(x)$. 
Donc si $x \dans J_{S/\Sigma}(U,T,\gamma)$, alors $\phi(x) \dans p\mathcal O_T$, on en déduit donc le résultat en tordant par $\phi$ la première assertion.
\edem
%

\cor
\label{cor33}
L'application
\[V^f : \bigwedge^{q_{\tau_0}} \mathcal E_{\tau_0} \fleche \Fil\left(\bigwedge^{q_{\tau_0}} \mathcal E_{\tau_0}^{(p^f)}\right)\]
est divisible par $p^{k_{\tau_0}}$ sur les épaississements $(U,T,\delta)$, $p$-adiques ou munis d'un relèvement fort de Frobenius au sens de la proposition précédente, où,
\[ k_{\tau_0} = \sum_{q_\tau < q_{\tau_0}} q_{\tau_0} - q_\tau.\]
C'est-à-dire que son image est incluse dans
\[p^{k_{\tau_0}}\Fil\left(\bigwedge^{q_{\tau_0}} \mathcal E_{\tau_0}^{(p^f)}\right)_{(U,T,\delta)}.\]
\ecor

\dem
Considérons la décomposition suivante du morphisme $V^f$,
\begin{center}
\begin{tikzpicture}[description/.style={fill=white,inner sep=2pt}] 
\matrix (m) [matrix of math nodes, row sep=3em, column sep=2.5em, text height=1.5ex, text depth=0.25ex] at (0,0)
{ 
\bigwedge^{q_{\tau_0}} \Fil \mathcal E_{\tau_0} &
\bigwedge^{q_{\tau_0}} \Fil \mathcal E_{\sigma^{-1}\tau_0}^{(p)} &
{\dots}
& 
\bigwedge^{q_{\tau_0}} \Fil \mathcal E_{\sigma^{-f+1}\tau_0}^{(p^{f-1})} &
\bigwedge^{q_{\tau_0}} \Fil \mathcal E_{\tau_0}^{(p^f)}\\
\bigwedge^{q_{\tau_0}} \mathcal E_{\tau_0} & \bigwedge^{q_{\tau_0}}  \mathcal E_{\sigma^{-1}\tau_0}^{(p)} & 
{\dots} 
&
\bigwedge^{q_{\tau_0}} \mathcal E_{\sigma^{-f+1}\tau_0}^{(p^{f-1})} & \bigwedge^{q_{\tau_0}} \mathcal E_{\tau_0}^{(p^f)}  \\ 
};

\path[->,font=\scriptsize] 
(m-1-1) edge node[auto] {$V$} (m-1-2)
(m-2-1) edge node[auto] {$V$} (m-1-2)
(m-1-2) edge node[auto] {$V$} (m-1-3)
(m-2-2) edge node[auto] {$V$} (m-1-3)
(m-1-4) edge node[auto] {$V$} (m-1-5)
(m-2-4) edge node[auto] {$V$} (m-1-5)
(m-1-1) edge node[auto] {} (m-2-1)
(m-1-2) edge node[auto] {} (m-2-2)
(m-1-4) edge node[auto] {} (m-2-4)
(m-1-5) edge node[auto] {} (m-2-5)
;
\end{tikzpicture}
\end{center}
De plus, d'après la proposition précédente, pour tout $\tau \dans \mathcal I$ tel que $q_\tau < q_{\tau_0}$, la flèche verticale se factorise en fait sur les ouverts considérés, par,
\[\bigwedge^{q_{\tau_0}} \Fil \mathcal E_\tau^{(p^r)} \fleche p^{q_{\tau_0}-q_\tau} \bigwedge^{q_{\tau_0}} \mathcal E_\tau^{(p^r)}.\]
La commutativité des triangles est donnée par la proposition \ref{pro2}.
On en déduit donc que l'application $V^f$ sur $\bigwedge^{q_{\tau_0}} \mathcal E_{\tau_0}$, se factorise sur les ouverts $(U,T,\delta)$ comme dans l'énoncé, par
\[ V^f : \bigwedge^{q_{\tau_0}} \mathcal E_{\tau_0} \fleche p^{k_{\tau_0}} \Fil \bigwedge^{q_{\tau_0}} \mathcal E_{\tau_0}^{(p^f)}.\qedhere\]
\edem

\rem
\label{remfort}
En fait on a montré un peu plus fort, c'est-à-dire que sur ces ouverts l'application
\[V^{f-1} : \bigwedge^{q_{\tau_0}} \mathcal E_{\tau_0} \fleche \left(\bigwedge^{q_{\tau_0}} \mathcal E_{\tau_0}^{(p^{f-1})}\right),\]
est divisible par $p^{k_{\tau_0}}$.
\erem

\thr
\label{thruni}
Supposons $S$ lisse et $r> k_{\tau_0}$.
Il existe une unique application entre faisceaux cristallins,
\[\phi_{\tau_0} :\displaystyle\left(\bigwedge^{q_{\tau_0}} \mathcal E_{\tau_0}\right) \otimes \mathcal O_{S/\Sigma}/p^{r - k_{\tau_0}}\fleche \displaystyle\Fil\left(\bigwedge^{q_{\tau_0}} \mathcal E_{\tau_0}^{(p^f)}\right)\otimes \mathcal O_{S/\Sigma}/p^{r - k_{\tau_0}},\]
telle que $p^{k_{\tau_0}}\phi_{\tau_0}$ s'étende en l'application $V^f$, modulo $p^r$, de la remarque précédente.
\ethr

\dem
Il suffit de prouver l'existence et l'unicité de $\phi_{\tau_0}$ localement, l'unicité 
locale impliquera alors automatiquement que la construction se recolle, et donc l'existence 
globale.
Supposons $S$ affine (et lisse) sur $\Sigma_1$, alors 
il existe un carré cartésien,
\begin{center}
\begin{tikzpicture}[description/.style={fill=white,inner sep=2pt}] 
\matrix (m) [matrix of math nodes, row sep=3em, column sep=2.5em, text height=1.5ex, text depth=0.25ex] at (0,0)
{ 
S & & S_\infty \\
\Spec(\kappa_F) & & \Sigma\\
 };

\path[->,font=\scriptsize] 
(m-1-1) edge node[auto] {} (m-1-3)
(m-1-1) edge node[auto,left] {} (m-2-1)
(m-1-3) edge node[auto,left] {lisse} (m-2-3)
(m-2-1) edge node[auto,right] {} (m-2-3);
\end{tikzpicture}
\end{center}
où $S_\infty \fleche \Sigma = \Spec(W(\kappa_F))$ est lisse. Notons $S_r = S_\infty \times \Spec(W(\kappa_F)/p^r)$. Par lissité de $S_r$, $p^{r-1}\mathcal O_{S_r} \neq 0$, 
et l'épaississement $(S \hookrightarrow S_\infty)$ est $p$-adique.

D'après \cite{BerO} Théorème 6.6, considérons $(M,\nabla)$ le $\mathcal O_{S_\infty}$-module à connexion associé à 
$\bigwedge^{q_{\tau_0}} \mathcal E_{\tau_0}$, c'est à dire,
\[ M = \left(\bigwedge^{q_{\tau_0}} \mathcal E_{\tau_0}\right)_{(S\hookrightarrow S_\infty)}.\]
C'est un $\mathcal O_{S_r}$-module localement libre, car $G$ est un groupe de Barsotti-Tate tronqué de rang $r$ (cf. \cite{BBM} Théorème 3.3.10). Le morphisme \[V^f : M \fleche M^{(p^f)}\] 
est divisible par $p^{k_{\tau_0}}$, d'après le corollaire \ref{cor33}, donc il existe 
\[\psi : M \fleche M^{(p^f)} \quad \text{tel que } p^{k_{\tau_0}}\psi = V^f.\]
De plus, $\psi$ est unique modulo $p^{r - k_{\tau_0}}$ car $M^{(p^f)}$ est localement libre sur $\mathcal O_{S_r}$.
Et pour que $\psi$ donne un morphisme de cristaux, il faut qu'il soit compatible aux 
connexions de $M$ et $M^{(p^f)}$.
Or il l'est après multiplication par $p^{k_{\tau_0}}$, et $M$ est localement libre sur 
$\mathcal O_{S_r}$, donc $\psi$ est compatible aux connexions au moins après réduction
 modulo $p^{r - k_{\tau_0}}$, c'est-à-dire que,
\[ \overline{\psi} : M/p^{r - k_{\tau_0}}M \fleche M^{(p^f)}/p^{r - k_{\tau_0}}M^{(p^f)}\]
est compatible aux connexions. On en déduit donc, toujours par \cite{BerO} Théorème 6.6, un unique morphisme de cristaux,
\[ \overline{\psi} : \bigwedge^{q_{\tau_0}} \mathcal E_{\tau_0}/p^{r - k_{\tau_0}} \fleche 
\bigwedge^{q_{\tau_0}} \mathcal E^{(p^f)}_{\tau_0}/p^{r - k_{\tau_0}}.\]

\lem
Le morphisme précédemment construit,
\[ \overline{\psi} : \bigwedge^{q_{\tau_0}} \mathcal E_{\tau_0}/p^{r - k_{\tau_0}} \fleche 
\bigwedge^{q_{\tau_0}} \mathcal E^{(p^f)}_{\tau_0}/p^{r - k_{\tau_0}},\]
se factorise (à l'arrivée) par,
\[\phi_{\tau_0}  : \bigwedge^{q_{\tau_0}} \mathcal E_{\tau_0}/p^{r - k_{\tau_0}} \fleche \Fil\left(\bigwedge^{q_{\tau_0}} \mathcal E^{(p^f)}_{\tau_0}\right)/p^{r - k_{\tau_0}}.\]
\elem

\dem[Lemme] 
Malheureusement le faisceau cristallin par lequel on veut factoriser n'est pas un cristal, on va donc montrer la factorisation sur chaque épaississement $(U \hookrightarrow T, \delta)$ de $Cris(S/\Sigma)$. La factorisation est vraie pour $V^f$, d'après le corollaire 3.3. De plus pour l'épaississement à puissances divisées $(S \hookrightarrow S_\infty)$ d'idéal $p\mathcal O_{S_\infty}$, 
 \[\left(\Fil \bigwedge^{q_{\tau_0}} \mathcal E_{\tau_0}\right)_{(S\hookrightarrow S_\infty)} \quad
 \text{et} \quad \left(\Fil \bigwedge^{q_{\tau_0}} \mathcal E^{(p^f)}
 _{\tau_0}\right)_{(S\hookrightarrow S_\infty)}, \]
ne sont pas des $\mathcal O_{S_r}$-modules localement libres, 
mais on peut diviser comme précédemment, d'après le corollaire \ref{cor33} et la remarque \ref{remfort}, le morphisme de cristaux,
\[ V^{f-1} : \left(\bigwedge^{q_{\tau_0}} \mathcal E_{\tau_0}\right)_{(S\hookrightarrow S_r)} \fleche
 \left(\bigwedge^{q_{\tau_0}} \mathcal E^{(p^{f-1})}
 _{\tau_0}\right)_{(S\hookrightarrow S_\infty)}, \]
par $p^{k_{\tau_0}}$, car $ \left(\bigwedge^{q_{\tau_0}} \mathcal E^{(p^{f-1})}
 _{\tau_0}\right)$ est un cristal, et en composant avec la flèche, cf proposition \ref{pro2},
\[ V : \left(\bigwedge^{q_{\tau_0}} \mathcal E^{(p^{f-1})}
 _{\tau_0}\right)_{(S\hookrightarrow S_\infty)} \fleche  \left(\Fil \bigwedge^{q_{\tau_0}} \mathcal E^{(p^f)}_{\tau_0}\right)_{(S \hookrightarrow S_\infty)},\]
on obtient un morphisme entre les espaces désirés,
\[(\phi_{\tau_0})_{(S \hookrightarrow S_\infty)}  : \left(\bigwedge^{q_{\tau_0}} \mathcal E_{\tau_0}/p^{r - k_{\tau_0}}\right)_{(S \hookrightarrow S_\infty)} \fleche \left(\Fil\left(\bigwedge^{q_{\tau_0}} \mathcal E^{(p^f)}_{\tau_0}\right)/p^{r - k_{\tau_0}}\right)_{(S \hookrightarrow S_\infty)}.\]
 De plus si on compose ce dernier morphisme par,
\[  \left(\Fil \bigwedge^{q_{\tau_0}} \mathcal E^{(p^{f})}
 _{\tau_0}\right)_{(S\hookrightarrow S_\infty)} \subset  \left(\bigwedge^{q_{\tau_0}} \mathcal E^{(p^{f})}
 _{\tau_0}\right)_{(S\hookrightarrow S_\infty)},\]
et par l'unicité de $\overline{\psi}$, on en déduit que $\overline{\psi}_{(S \hookrightarrow S_\infty)}$ se factorise comme voulu.

 Mais maintenant soit $(U \hookrightarrow T,\delta)$ un épaississement quelconque de $Cris(S/\Sigma)$ avec $U$ et $T$ affines. Comme $S_\infty$ est lisse, il existe un morphisme $u$ s'insérant dans le diagramme suivant,
\begin{center}
\begin{tikzpicture}[description/.style={fill=white,inner sep=2pt}] 
\matrix (m) [matrix of math nodes, row sep=3em, column sep=2.5em, text height=1.5ex, text depth=0.25ex] at (0,0)
{ 
U  & & T \\
S & & S_\infty \\
 & & \Sigma\\
 };

\path[->,font=\scriptsize] 
(m-2-1) edge node[auto] {} (m-3-3)
(m-2-3) edge node[auto] {} (m-3-3)
(m-2-1) edge node[auto,left] {} (m-2-3)
(m-1-1) edge node[auto,right] {} (m-1-3)
(m-1-1) edge node[auto,right] {} (m-2-1);
\path[dashed,->,font=\scriptsize] 
(m-1-3) edge node[auto,right] {$u$} (m-2-3);
\end{tikzpicture}
\end{center}

Par propriété universelle de l'enveloppe à puissances divisées (qui est ici $S_\infty$ puisque cet espace est lisse sur $\Sigma$), on a que $u$ 
est un morphisme dans $Cris(S/\Sigma)$.
Comme $\overline{\psi}$ est un morphisme de cristal, il vérifie que le diagramme suivant commute, 
 \begin{center}
\begin{tikzpicture}[description/.style={fill=white,inner sep=2pt}] 
\matrix (m) [matrix of math nodes, row sep=3em, column sep=2.5em, text height=1.5ex, text depth=0.25ex] at (0,0)
{ 
 \left(\bigwedge^{q_{\tau_0}} \mathcal E_{\tau_0}\right)_{(U\hookrightarrow T)}/p^{r - k_{\tau_0}} & &  \left(\bigwedge^{q_{\tau_0}} \mathcal E^{(p^f)}_{\tau_0}\right)_{(U\hookrightarrow T)}/p^{r - k_{\tau_0}} \\
u^*\left(\bigwedge^{q_{\tau_0}} \mathcal E_{\tau_0}\right)_{(S\hookrightarrow S_\infty)}/p^{r - k_{\tau_0}} & &  
u^*\left(\bigwedge^{q_{\tau_0}} \mathcal E^{(p^f)}_{\tau_0}\right)_{(S\hookrightarrow S_\infty)}/p^{r - k_{\tau_0}}
\\
 };

\path[->,font=\scriptsize] 
(m-2-1) edge node[auto] {$\rho_u'$} (m-1-1)
(m-2-1) edge node[auto,right] {$\simeq$} (m-1-1)
(m-2-1) edge node[auto] {$u^*{\overline{\psi}}$} (m-2-3)
(m-2-3) edge node[auto] {$\rho_u'$} (m-1-3)
(m-2-3) edge node[auto,right] {$\simeq$} (m-1-3);
\path[dashed,->,font=\scriptsize] 
(m-1-1) edge node[auto] {${\overline{\psi}}_{(U\hookrightarrow T)}$} (m-1-3);
\end{tikzpicture}
\end{center}

Mais on vient de montrer qu'on peut factoriser la flèche du bas, et donc on en déduit une factorisation,
\begin{center}
\begin{tikzpicture}[description/.style={fill=white,inner sep=2pt}] 
\matrix (m) [matrix of math nodes, row sep=3em, column sep=2.5em, text height=1.5ex, text depth=0.25ex] at (0,0)
{ 
 \left(\bigwedge^{q_{\tau_0}} \mathcal E_{\tau_0}\right)_{(U\hookrightarrow T)} & & 
 \left(\Fil \bigwedge^{q_{\tau_0}} \mathcal E^{(p^f)}
 _{\tau_0}\right)_{(U\hookrightarrow T)}  & &  \left(\bigwedge^{q_{\tau_0}} \mathcal E^{(p^f)}_{\tau_0}\right)_{(U\hookrightarrow T)} \\
u^*\left(\bigwedge^{q_{\tau_0}} \mathcal E_{\tau_0}\right)_{(S\hookrightarrow S_\infty)} & &  
u^*\left(\Fil \bigwedge^{q_{\tau_0}} \mathcal E^{(p^f)} _{\tau_0}\right)_{(S\hookrightarrow S_\infty)} & & 
u^*\left(\bigwedge^{q_{\tau_0}} \mathcal E^{(p^f)}_{\tau_0}\right)_{(S\hookrightarrow S_\infty)} 
\\
 };

\path[->,font=\scriptsize] 
(m-2-1) edge node[auto] {$\rho_u'$} (m-1-1)
(m-2-1) edge node[auto,right] {$\simeq$} (m-1-1)
(m-2-3) edge node[auto] {$\rho_u$} (m-1-3)
(m-2-1) edge node[auto] {$u^*{\phi_{\tau_0}}_{(S \hookrightarrow S_\infty)}$} (m-2-3)
(m-2-5) edge node[auto] {$\rho_u'$} (m-1-5)
(m-2-5) edge node[auto,right] {$\simeq$} (m-1-5)
(m-2-3) edge node[auto] {$u^*i_{S_\infty}$} (m-2-5)
(m-1-3) edge node[auto] {$i_T$} (m-1-5)
;
\path[dashed,->,font=\scriptsize] 
(m-1-1) edge node[auto] {${\phi_{\tau_0}}_{(U\hookrightarrow T),u}$} (m-1-3);
\end{tikzpicture}
\end{center}
telle que $i_T\circ {\phi_{\tau_0}}_{(U\hookrightarrow T),u} = \overline{\psi}_{(U\hookrightarrow T)}$ 
ne dépends pas de $u$. Or $i_T$ est un monomorphisme, donc ${\phi_{\tau_0}}_{(U\hookrightarrow T),u}$ ne dépend pas du choix du morphisme 
\[ u : T \fleche S_r,\] 
et donc factorise $\overline{\psi}_{(U,T,\delta)}$.
\edem

Il reste donc à montrer que pour tout morphisme 
\[ f : \quotient{\displaystyle\bigwedge^{q_{\tau_0}} \mathcal E_{\tau_0}}{p^{r - k_{\tau_0}}}\fleche \quotient{\displaystyle\Fil\left(\bigwedge^{q_{\tau_0}} \mathcal E_{\tau_0}^{(p^f)}\right)}{p^{r - k_{\tau_0}}},\]
tel que le relevé de $p^{k_{\tau_0}}f$ modulo $p^r$ soit $V^f$, alors après composition avec,
\[ i :  \left(\Fil \bigwedge^{q_{\tau_0}} \mathcal E^{(p^f)}
 _{\tau_0}\right)/p^{r-k_{\tau_0}}\subset \left(\bigwedge^{q_{\tau_0}} \mathcal E^{(p^f)}_{\tau_0}/p^{r - k_{\tau_0}} \right),\]
$i \circ f$ est le morphisme $\overline{\psi}$. Mais par hypothèse $i \circ p^{k_{\tau_0}}f = p^{k_{\tau_0}} i \circ f$ est un morphisme de cristaux, donc, 
\[i \circ f :  \quotient{\displaystyle\bigwedge^{q_{\tau_0}} \mathcal E_{\tau_0}}{p^{r - k_{\tau_0}}}\fleche \quotient{\displaystyle{\bigwedge^{q_{\tau_0}} \mathcal E_{\tau_0}^{(p^f)}}}{p^{r - k_{\tau_0}}},\]
est un morphisme de cristaux (par \cite{BerO} Théorème 6.6, comme précédemment) qui divise $V^f$. C'est donc $\overline{\psi}$. Comme $i$ est un monomorphisme, cela conclut.
\edem
 
\rem L'unicité va être centrale dans la suite, elle nous permettra de voir que $\phi_{\tau_0}$ (et donc dans la suite $\widetilde\Ha_\tau$) vérifie toutes les compatibilités désirées.

Dans le cas où $G/S$ est un groupe $p$-divisible (non tronqué), tous les modules sont des $\ZZ_p$-modules 
libres, et on a bien existence et unicité d'une application,
\[ \phi_{\tau_0} :  \left(\bigwedge^{q_{\tau_0}} \mathcal E_{\tau_0}\right)\fleche 
\Fil\left(\bigwedge^{q_{\tau_0}} \mathcal E_{\tau_0}^{(p^f)}\right),\]
vérifiant que $p^{k_{\tau_0}}\phi_{\tau_0} = V^f$.
\erem

On voudrait voir $\phi_{\tau_0}$ comme une application sur 
$\bigwedge^{q_{\tau_0}} \omega_{G^D,\tau_0}$, pour cela notons $\zeta_{\tau_0}$ la restriction de $\phi_{\tau_0}$ a $\Fil\left(\bigwedge^{q_{\tau_0}} \mathcal E_{\tau_0}\right),$ 
et considérons,

\[ \Fil_2^{\tau_0} = \Ker\left( \Fil \left( \bigwedge^{q_{\tau_0}} \mathcal E_{\tau_0}\right)/p^{r - k_{\tau_0}}
\overset{\pi}{\fleche} \bigwedge^{q_{\tau_0}} \mathcal E_{\tau_0}/J_{S/\Sigma}\mathcal E_{\tau_0}\right),\]
et notons 
\[\Fil_2^{\tau_0,(p^f)} = \Ker\left( \Fil \left( \bigwedge^{q_{\tau_0}} \mathcal E^{(p^f)}_{\tau_0}\right) /p^{r - k_{\tau_0}}
\overset{\pi}{\fleche} 
\bigwedge^{q_{\tau_0}} \mathcal E^{(p^f)}_{\tau_0}/J_{S/\Sigma}\mathcal E^{(p^f)}_{\tau_0}\right).\]

Vérifions tout d'abord que $\zeta_{\tau_0}$ passe bien au quotient par ces sous-modules :

\pro
Le morphisme cristallin défini précédemment,
\[\zeta_{\tau_0} : \quotient{\displaystyle\Fil \left(\bigwedge^{q_{\tau_0}} \mathcal E_{\tau_0}\right)}{p^{r - k_{\tau_0}}}\fleche \quotient{\displaystyle\Fil\left(\bigwedge^{q_{\tau_0}} \mathcal E_{\tau_0}^{(p^f)}\right)}{p^{r - k_{\tau_0}}},\]
passe au quotient par les sous modules $\Fil_2^{q_{\tau_0}}$ et $\Fil_2^{q_{\tau_0,(p^f)}}$, 
c'est à dire que,
\[\zeta_{\tau_0}(\Fil_2^{\tau_0}) \subset \Fil_2^{\tau_0,(p^f)}.\]
\epro

\dem
Par le théorème précédent, il existe un morphisme cristallin,
\[ \phi_{\tau_0} : \quotient{\displaystyle\bigwedge^{q_{\tau_0}} \mathcal E_{\tau_0}}{p^{r - k_{\tau_0}}}\fleche \quotient{\displaystyle\Fil\left(\bigwedge^{q_{\tau_0}} \mathcal E_{\tau_0}^{(p^f)}\right)}{p^{r - k_{\tau_0}}},\]
qui prolonge $\zeta_{\tau_0}$ et tel que $p^{k_{\tau_0}}\phi_{\tau_0} = V^f$.
Regardons le sous-module,
\[\Ker\left( \left( \bigwedge^{q_{\tau_0}} \mathcal E_{\tau_0}\right)/p^{r - k_{\tau_0}}
\overset{\pi}{\fleche} \bigwedge^{q_{\tau_0}} \mathcal E_{\tau_0}/J_{S/\Sigma}\mathcal E_{\tau_0}\right) = 
\im\left( J_{S/\Sigma}\mathcal E_{\tau_0}\otimes \bigwedge^{q_{\tau_0}-1} \mathcal E_{\tau_0}  \fleche 
\left(\bigwedge^{q_{\tau_0}} \mathcal E_{\tau_0}\right)/p^{r - k_{\tau_0}}\right).\]
%
L'application $V^f$ est divisible par $p^{k_{\tau_0}}$, on a donc que l'application $\phi_{\tau_0}$ restreinte à ce noyau est d'image dans,
\[J_{S/\Sigma}\left(\bigwedge^{q_{\tau_0}} \mathcal E^{(p^f)}_{\tau_0}\right)/p^{r - k_{\tau_0}}  =
\Ker\left( \left( \bigwedge^{q_{\tau_0}} \mathcal E^{(p^f)}_{\tau_0}\right)/p^{r - k_{\tau_0}}\overset{\pi}{\fleche} \bigwedge^{q_{\tau_0}} \mathcal E^{(p^f)}_{\tau_0}/J_{S/\Sigma}\mathcal E^{(p^f)}_{\tau_0}\right),\]
donc nulle après projection dans $\bigwedge^{q_{\tau_0}} \mathcal E_{\tau_0}^{(p^f)}/J_{S/\Sigma}\mathcal E_{\tau_0}^{(p^f)}$. Donc 
%
%
\[ \zeta_{\tau_0}(\Fil_2^{\tau_0}) \subset \Fil_2^{\tau_0,(p^f)}.\qedhere\]
\edem

Avant de conclure, réécrivons plus simplement ces quotients.

\lem
On a un isomorphisme,
\[ \quotient{\displaystyle\Fil\left(\bigwedge^{q_{\tau_0}} \mathcal E_{\tau_0}\right)}{\Fil_2^{\tau_0}} \simeq
\bigwedge^{q_{\tau_0}} i_{S/\Sigma*}\omega_{G^D,\tau_0}.\]
\elem

\rem
On a la même chose pour $\mathcal E_{\tau_0}^{(p^f)}$ en appliquant le lemme à $G^{(p^f)}$.
\erem

\dem 
On a une suite exacte,
\[ 0 \fleche J_{S/\Sigma}\mathcal E_{\tau_0} \fleche \Fil \mathcal E_{\tau_0} \fleche i_{S/\Sigma*}\omega_{G^D,\tau_0} \fleche 0,\]
on en déduit donc une flèche,
\[ \bigwedge^{q_{\tau_0}} \Fil \mathcal E_{\tau_0} \fleche \bigwedge^{q_{\tau_0}} i_{S/\Sigma*}\omega_{G^D,\tau_0} \fleche 0.\]
Il faut d'abord voir que cette flèche se factorise par,
\[   \Fil\bigwedge^{q_{\tau_0}} \mathcal E_{\tau_0}.\]
Or on a une suite exacte,
\[ 0 \fleche  \bigwedge^{q_{\tau_0}} i_{S/\Sigma*} \omega_{G^D,\tau_0} \fleche  \bigwedge^{q_{\tau_0}} \mathcal E_{\tau_0}/J_{S/\Sigma*}\mathcal E_{\tau_0} \fleche  \bigwedge^{q_{\tau_0}} i_{S/\Sigma*} \omega_{G,\tau_0}^\vee \fleche 0.\]
En effet, cette suite est exacte pour tout épaississement $(U\hookrightarrow T) \dans Cris(S/\Sigma)$ car ${\omega_{G^D,\tau_0}\otimes \mathcal O_U}$ est localement facteur direct de $\mathcal E_{\tau_0,(U \overset{id}{\rightarrow} U)}$.
Et par définition de $\Fil \mathcal E_{\tau_0}$ on a un carré commutatif,
\begin{center}
\begin{tikzpicture}[description/.style={fill=white,inner sep=2pt}] 
\matrix (m) [matrix of math nodes, row sep=3em, column sep=2.5em, text height=1.5ex, text depth=0.25ex] at (0,0)
{ 
 \bigwedge^{q_{\tau_0}}   i_{S/\Sigma*} \omega_{G^D,\tau_0}& 
 \bigwedge^{q_{\tau_0}} \left(\mathcal E_{\tau_0}/J_{S/\Sigma}\mathcal E_{\tau_0}\right)
 \\
\bigwedge^{q_{\tau_0}}\Fil \mathcal E_{\tau_0}& 
\bigwedge^{q_{\tau_0}} \mathcal E_{\tau_0}
\\
 };

\path[->,font=\scriptsize] 
(m-2-1) edge node[auto] {} (m-2-2)
(m-2-1) edge node[auto] {} (m-1-1)
(m-2-2) edge node[auto] {} (m-1-2)
;
\path[ right hook->,font=\scriptsize] 
(m-1-1) edge node[auto] {} (m-1-2)
;
\end{tikzpicture}
\end{center}
Donc l'application cherchée se factorise bien,
\begin{center}
\begin{tikzpicture}[description/.style={fill=white,inner sep=2pt}] 
\matrix (m) [matrix of math nodes, row sep=3em, column sep=2.5em, text height=1.5ex, text depth=0.25ex] at (0,0)
{ 
 \bigwedge^{q_{\tau_0}}\Fil \mathcal E_{\tau_0} &  & \bigwedge^{q_{\tau_0}}   i_{S/\Sigma*} \omega_{G^D,\tau_0} \\

 & \Fil \left(\bigwedge^{q_{\tau_0}}\mathcal E_{\tau_0}\right) & 
\\
 };

\path[->,font=\scriptsize] 
(m-1-1) edge node[auto] {} (m-1-3)
(m-2-2) edge node[auto] {$\pi$} (m-1-3)
(m-1-1) edge node[auto] {} (m-2-2)
;
\end{tikzpicture}
\end{center}
Il reste à déterminer le noyau de $\pi$, mais le carré commutatif précédent nous dit que le noyau de $\pi$ est exactement,
\[ \Ker\left( \Fil \bigwedge^{q_{\tau_0}} \mathcal E_{\tau_0} \fleche 
 \bigwedge^{q_{\tau_0}} \left(\mathcal E_{\tau_0}/J_{S/\Sigma}\mathcal E_{\tau_0}\right)\right),
 \]
 c'est à dire $\Fil_2^{\tau_0}$.
\edem

Grâce à tout ce qui précède, on a donc construit une application "de Hasse" associée à 
$\tau_0 \dans \mathcal I$, qui provient de $\zeta_{\tau_0}$ par passage au quotient,
\[ \widetilde{\Ha_{\tau_0}}(G) : i_{S/\Sigma*} \bigwedge^{q_{\tau_0}} \omega_{G^D,\tau_0} \fleche 
i_{S/\Sigma*} \bigwedge^{q_{\tau_0}} \omega_{G^D,\tau_0}^{(p^f)}.\]

Or le foncteur $i_{S/\Sigma*}$ est pleinement fidèle, cf \cite{BerO} 5.19, on en déduit donc un morphisme de $\mathcal O_{S}$-modules,
\[ \widetilde{\Ha_{\tau_0}}(G) : \bigwedge^{q_{\tau_0}} \omega_{G^D,\tau_0} \fleche 
\bigwedge^{q_{\tau_0}} \omega_{G^D,\tau_0}^{(p^f)}.\]

\pro
\label{proechelon}
Soit $G/S$ un groupe de Barsotti-Tate tronqué d'échelon $r$ (et supposons toujours $S$ lisse). Soit $\tau \dans \mathcal I$. Supposons donné,
\[ r > s > k_\tau.\]
Alors $G[p^s]$ est un groupe de Barsotti-Tate tronqué d'échelon $s<r$, et 
\[ \widetilde{\Ha_\tau}(G[p^s]) =\widetilde{ \Ha_\tau}(G).\]
\epro

\dem
Notons $G' := G[p^s]$, et d'après \cite{BBM} Théorème 3.3.3, on a,
\[ \mathcal Ext^1_{S/\Sigma}(G',\mathcal O_{S/\Sigma}) = \mathcal Ext^1_{S/\Sigma}(G,\mathcal O_{S/\Sigma}) \otimes \mathcal O_{S/\Sigma}/p^s.\]
C'est donc un $\mathcal O_{S/\Sigma}/p^s$-module localement libre, et le carré de cristaux suivant est commutatif,
\begin{center}
\begin{tikzpicture}[description/.style={fill=white,inner sep=2pt}] 
\matrix (m) [matrix of math nodes, row sep=3em, column sep=2.5em, text height=1.5ex, text depth=0.25ex] at (0,0)
{ 
 \bigwedge^{q_{\tau}}\mathcal Ext^1_{S/\Sigma}(G,\mathcal O_{S/\Sigma})_\tau &  &   \bigwedge^{q_{\tau}}\mathcal Ext^1_{S/\Sigma}(G,\mathcal O_{S/\Sigma})_\tau^{(p^f)}\\
 \bigwedge^{q_{\tau}}\mathcal Ext^1_{S/\Sigma}(G',\mathcal O_{S/\Sigma})_\tau &  &   \bigwedge^{q_{\tau}}\mathcal Ext^1_{S/\Sigma}(G',\mathcal O_{S/\Sigma})_\tau^{(p^f)}\\
 };

\path[->,font=\scriptsize] 
(m-1-1) edge node[auto] {$V^f$} (m-1-3)
(m-2-1) edge node[auto] {$V^f$} (m-2-3)
;
\path[->>,font=\scriptsize] 
(m-1-1) edge node[auto] {$$} (m-2-1)
(m-1-3) edge node[auto] {$$} (m-2-3)
;
\end{tikzpicture}
\end{center}
Remarquons aussi que $S$ étant de caractéristique $p$, canoniquement,
\[ \omega_{G^D} = \omega_{G'^D}.\]
Mais d'après le théorème \ref{thruni}, il existe une unique application,
\[\bigwedge^{q_{\tau}}\mathcal Ext^1_{S/\Sigma}(G',\mathcal O_{S/\Sigma})_\tau \pmod{p^s}  \overset{ \phi_\tau(G')}{\fleche} \bigwedge^{q_{\tau}}\mathcal Ext^1_{S/\Sigma}(G',\mathcal O_{S/\Sigma})_\tau^{(p^f)}\pmod{p^{s-k_\tau}},\]
qui divise $V^f$. Donc $\phi_\tau(G') = \phi_\tau(G) \pmod{p^{s-k_\tau}}.$
Il suffit ensuite de remarquer que la flèche,
\[\bigwedge^{q_{\tau}}\mathcal Ext^1_{S/\Sigma}(G,\mathcal O_{S/\Sigma})_\tau \twoheadrightarrow i_{S/\Sigma*} \bigwedge^{q_{\tau}} \omega_{G^D,\tau},\]
ainsi que ses variantes tordues par Frobenius, se factorise par,
\[ \bigwedge^{q_{\tau}}\mathcal Ext^1_{S/\Sigma}(G,\mathcal O_{S/\Sigma})_\tau \twoheadrightarrow \bigwedge^{q_{\tau}}\mathcal Ext^1_{S/\Sigma}(G',\mathcal O_{S/\Sigma})_\tau.\qedhere\]
\edem

\defi
Soit $G/S$ un groupe de Barsotti-Tate tronqué d'échelon \[r > \max_{\tau \in \mathcal I} k_\tau.\] Si $\tau \dans \mathcal I$ est tel que 
\[ q_\tau = \dim_{S} \omega_{G^D,\tau} = 0,.\]
on pose alors,
\[ \det(\omega_{G^D,\tau}) = \mathcal O_{S}, \quad \text{et} \quad \widetilde{\Ha_\tau}(G) = \id.\]

On définit l'\textit{invariant de Hasse partiel}, associé à $\tau \dans \mathcal I$, de 
$G/S$, par,
\[\widetilde\Ha_{\tau}(G) \dans \Gamma(S,\det\left(\omega_{G^D,\tau}\right)^{\otimes(p^f-1)}),\]
donné par le morphisme $\widetilde{\Ha_{\tau}}(G)$ précédent, sous l'identification,
\[\det\left(\omega_{G^D,\tau}^{(p^f)}\right)= \det\left(\omega_{G^D,\tau}\right)^{\otimes p^f}.\]
On appelle $\mu$-\textit{invariant de Hasse} la section produit des sections $\Ha_{\tau}(G)$, pour $\tau \dans \mathcal I$, c'est à dire,
\[ \widetilde{^\mu\Ha}(G) = \bigotimes_{\tau \in \mathcal I} \widetilde{\Ha_{\tau}}(G) \dans \Gamma(S, \det(\omega_{G^D})^{\otimes(p^f-1)}).\]
D'après la proposition \ref{proechelon}, cette définition est bien cohérente, et ne dépends que de la $p^s$-torsion de $G$, où,
\[s = \max_{\tau \in \mathcal I} k_\tau +1.\]
\edefi

\rem
Dans le cas où la base $S$ est une variété de Hilbert-Siegel (modulo un $p$ de bonne réduction), et $A$ est la variété abélienne universelle, on a aussi des invariants de Hasse partiels, donnés par les sections (cf. \cite{AndGo} par exemple),
\[ V_\tau \dans \Gamma(S, \det(\omega_{A,\sigma^{-1}\tau}^{\otimes p}\otimes\omega_{A,\tau}^{\otimes(-1)}))\]
car les $q_{\tau}$ sont tous égaux. Mais ces invariants de Hasse partiels sont différents de ceux défini dans cet article : on a la relation 
(comme tous les $q_\tau$ sont égaux, $k_\tau =0$) ,
\[ \widetilde\Ha_\tau(A[p^\infty]) = V_\tau\otimes V_{\sigma^{-1}\tau}^{\otimes p}\otimes \dots \otimes V_{\sigma\tau}^{\otimes p^{f-1}}.\]
En particulier l'invariant de Hasse partiel $\widetilde\Ha_\tau$ est "moins précis" que $V_\tau$, et il n'est pas réduit dans ce cas.
\erem

\section{Descente au corps réflexe}
\label{sect5}
Considérons maintenant un schéma $S/\kappa_F$ de caractéristique $p$, et $G/S$ un groupe de Barsotti-Tate tronqué d'échelon $r$ muni d'une action de $\mathcal O = \mathcal O_F$, 
et de signature $(p_\tau,q_\tau)_{\tau \in \mathcal I}$ donnée. Sans donner plus de précision, supposons $r$ assez grand (au sens de la section précédente) pour les constructions 
utilisées aient du sens.

\defi
Soit $\kappa_E$ le plus petit sous-corps de $\kappa_F$ qui laisse invariante la signature, c'est-à-dire tel que,
$\forall \theta \in \Gal(\kappa_F/\kappa_E), \forall \tau \in \mathcal I,$
\[ q_{\tau\circ\theta} = q_{\tau}.\]
\edefi

Cela nous fixe un plongement, $\kappa_E \fleche \kappa_F$.
On suppose qu'en plus, $S$ et $G$ descendent à $\kappa_E$, c'est à dire,
qu'il existe $S_0/\kappa_E$, et un groupe $p$-divisible (tronqué) muni d'une action de $\mathcal O$, noté $G_0$ sur $S_0$, tel que,
\[S = S_0 \otimes_{\kappa_E} \kappa_F \quad \text{et} \quad G_0 \times_{S_0} S = G,\]
et l'action sur $G$ est donnée par extension des scalaires par celle de $G_0/S_0$.

Sur $S$, c'est à dire associé à $G$, on a construit $\widetilde\Ha_\tau(G)$ pour tout $\tau \dans \mathcal I$.
On voudrait montrer qu'en fait on peut redescendre la construction de l'invariant de Hasse 
$\mu$-ordinaire à $S_0$. Soit alors $\theta \dans \Gal(\kappa_F/\kappa_E)$, et notons,
\[ G^{(\theta)} = G \otimes_{\kappa_F,\theta} \kappa_F = G \times_{S} (S \otimes_{\kappa_F,\theta} \kappa_F).\]
On a naturellement des isomorphismes,
\[ \rho_\theta : G^{(\theta)} \overset{\simeq}{\fleche} G,\]
qui induisent des isomorphismes,
\[ \mathcal E_{\theta} := \mathcal{E}xt^1_{S^{(\theta)}/\sigma}(G^{(\theta),D},\mathcal 
O_{S^{(\theta)}/\sigma}) \simeq  \theta^*\mathcal{E}xt^1_{S/\sigma}(G^D,\mathcal 
O_{S/\sigma}) = \theta^*\mathcal E.\]

ainsi que,
\[ \omega_{G^{(\theta)}} \simeq \theta^*\omega_{G} \quad \text{et} \quad \omega_{G^{(\theta),D}} 
= \theta^* \omega_{G^D}.\]

Et donc grâce à l'isomorphisme $\rho_\theta$, on en déduit des isomorphismes,
\[ \theta^* \mathcal E \simeq \mathcal E, \quad \theta^*\omega_G \simeq \omega_G \quad \text{et} \quad \theta^*\omega_{G^D} \simeq \omega_{G^D}.\]

Malheureusement, a priori, les composantes correspondant à $\tau \dans \mathcal I$ ne sont pas 
préservées par torsion par $\theta$ (et donc le fibré sur $S$, $\omega_{G^D,\tau}$, ne redescend 
pas a priori à $S_0$). Mais on a des isomorphismes,
\[\left(\omega_{G^{(\theta)}}\right)_\tau = \left( \theta^*\omega_G\right)_\tau \simeq \theta^*\left( \omega_{G,\tau\circ\theta}\right) \quad 
\text{et} \quad \left(\omega_{G^{(\theta),D}}\right)_\tau =\left( \theta^*\omega_{G^D}\right)_\tau \simeq 
\theta^*\left( \omega_{G^D,\tau\circ\theta}\right).\]
Soit $a,b \dans \NN$ tels que $a+b = h$. Il est alors naturel de considérer,
\[ \omega_{G,a} = \bigoplus_{\tau \in \mathcal I, q_\tau = b} \omega_{G,\tau}
\quad \text{et} \quad  \omega_{G^D,b} = \bigoplus_{\tau \in \mathcal I, q_\tau = b} 
\omega_{G^D,\tau}.\]
Les isomorphismes précédents induisent alors des isomorphismes,
\[ \theta^*(\omega_{G,a}) \simeq \omega_{G,a} \quad \text{et} \quad  r_\theta : \theta^*(\omega_{G^D,b}) \simeq \omega_{G^D,b}.\]
De plus ceux-ci sont compatibles à la condition de cocycle. En effet, c'est le cas de $\omega_G$ et $\omega_{G^D}$, puisque ces derniers proviennent de $S_0$, et donc c'est le 
cas de $\omega_{G,a}$ et $\omega_{G^D,b}$ par restriction. On en déduit par descente galoisienne qu'ils redescendent en des modules localement libres 
$\omega_{G_0,a}$ et $\omega_{G_0^D,b}$ sur $S_0$ !

\pro
Les morphismes de faisceaux sur $S$, pour tout $q \dans \{q_\tau, \tau \dans \mathcal I\}$,
\[ \bigotimes_{\tau \in \mathcal I, q_\tau =b} \widetilde\Ha_\tau(G) : \bigotimes_{\tau \in \mathcal I, q_\tau =b} \det\left(\omega_{G^D,\tau}\right) \fleche \bigotimes_{\tau \in \mathcal I, q_\tau =b} \det\left(\omega_{G^D,\tau}\right)^{\otimes p^f},\]
proviennent canoniquement de morphismes sur $S_0$, notés,
\[ \widetilde\Ha_b(G_0) : \det\left(\omega_{G_0^D,b}\right) \fleche  \det\left(\omega_{G_0^D,b}\right)^{\otimes p^f}.\]
En particulier l'invariant de Hasse $\mu$-ordinaire provient d'un morphisme sur $S_0$,
\[ \widetilde{^\mu\Ha}(G_0) : \det\left(\omega_{G_0^D}\right) \fleche  \det\left(\omega_{G_0^D}\right)^{\otimes p^f}.\]
\epro

\dem
Il suffit de montrer que $\widetilde\Ha_\tau(G^{(\theta)}) = \theta^*\widetilde\Ha_{\theta\circ\tau}(G)$ comme morphismes sur $S$, c'est à dire, que le carré suivant commute :
\begin{center}
\begin{tikzpicture}[description/.style={fill=white,inner sep=2pt}] 
\matrix (m) [matrix of math nodes, row sep=3em, column sep=2.5em, text height=1.5ex, text depth=0.25ex] at (0,0)
{ 
\det\left(\omega_{G^{D},\tau}\right)
 &  & \det\left(\omega_{G^{D},\tau}\right)^{\otimes p^f} \\
\theta^*\det\left(\omega_{G^{D},\tau\theta}\right)
 &  & \theta^*\det\left(\omega_{G^{D},\tau\theta}\right)^{\otimes p^f} \\
 };

\path[->,font=\scriptsize] 
(m-1-1) edge node[auto] {$\widetilde\Ha_\tau(G^{(\theta)})$} (m-1-3)
(m-2-1) edge node[auto] {$\theta^*\widetilde\Ha_{\theta\circ\tau}(G)$} (m-2-3)
(m-2-1) edge node[auto] {$\overset{r_\theta}{\simeq}$} (m-1-1)
(m-2-3) edge node[auto] {$\overset{r_\theta^{\otimes p^f}}{\simeq}$} (m-1-3)
;
\end{tikzpicture}
\end{center}
Or on a un isomorphisme $r_\theta : \theta^*\mathcal E \simeq \mathcal E$, dont on déduit aussi que,
\[ \Fil \theta^*\mathcal E = \theta^*(\Fil \mathcal E),\]
Et idem avec application des opérateurs $\bigwedge$, quotients, etc... De même,
\[ (\theta^*\mathcal E)_{\tau} = \theta^*\left(\mathcal E_{\theta\circ\tau}\right),\]
et idem avec application de $\Fil$, $\bigwedge$, et/ou quotients, etc...
De plus, si $\sigma$ dénote le Frobenius de $F$, on a évidement que $\theta \circ \sigma = \sigma \circ \theta$, on en déduit que le carré,
\begin{center}
\begin{tikzpicture}[description/.style={fill=white,inner sep=2pt}] 
\matrix (m) [matrix of math nodes, row sep=3em, column sep=2.5em, text height=1.5ex, text depth=0.25ex] at (0,0)
{ 
 \mathcal E_{\tau}  &  & \mathcal E_{\sigma^{-1}\circ\tau}^{(p)} & \\
\theta^*\left(\mathcal E_{\tau\theta}\right) 
 &  & \theta^*\left(\mathcal E_{\tau\theta\sigma^{-1}}\right)^{(p)}  &
=  \theta^*\left(\mathcal E_{\tau\sigma^{-1}\theta}\right)^{(p)}\\
 };

\path[->,font=\scriptsize] 
(m-1-1) edge node[auto] {$V$} (m-1-3)
(m-2-1) edge node[auto] {$\theta^*V$} (m-2-3)
(m-2-1) edge node[auto] {$\simeq$} (m-1-1)
(m-2-3) edge node[auto] {$\simeq$} (m-1-3)
;
\end{tikzpicture}
\end{center}
commute, et de même on a le carré commutatif,
\begin{center}
\begin{tikzpicture}[description/.style={fill=white,inner sep=2pt}] 
\matrix (m) [matrix of math nodes, row sep=3em, column sep=2.5em, text height=1.5ex, text depth=0.25ex] at (0,0)
{ 
\bigwedge^{q_\tau} \mathcal E_{\tau}  &  & \bigwedge^{q_\tau}\mathcal E_{\tau}^{(p^f)} & \\
\theta^*\left(\bigwedge^{q_\tau}\mathcal E_{\tau\theta}\right) 
 &  & \theta^*\left(\bigwedge^{q_\tau}\mathcal E_{\tau\theta}\right)^{(p^f)}  \\
 };

\path[->,font=\scriptsize] 
(m-1-1) edge node[auto] {$V^f$} (m-1-3)
(m-2-1) edge node[auto] {$\theta^*V^f$} (m-2-3)
(m-2-1) edge node[auto] {$\simeq$} (m-1-1)
(m-2-3) edge node[auto] {$\simeq$} (m-1-3)
;
\end{tikzpicture}
\end{center}
Mais par la proposition (\ref{thruni}) il existe une unique application,
\[\phi_{\tau} : \quotient{\displaystyle\bigwedge^{q_{\tau}} \mathcal E_{\tau}}{p^{r - k_{\tau}}}\fleche \quotient{\displaystyle\Fil\left(\bigwedge^{q_{\tau}} \mathcal E_{\tau}^{(p^f)}\right)}{p^{r - k_{\tau}}},\]
qui, multipliée par $p^{k_\tau}$, se relève en $V^f$. Or c'est aussi le cas de $\theta^*\phi_{\tau\theta}$, où $\phi_{\tau\theta}$ est l'application associée à $G$ et au plongement $\tau\theta$, donc on en déduit que,
\[ \phi_{\tau} = \theta^*\phi_{\tau\theta}.\]
Étant donné que $\widetilde\Ha_{\tau}$ s'en déduit par restriction et quotient (et que ces opérations commutent avec $\theta^*$) on en déduit que le premier diagramme commute, et donc que $\widetilde\Ha_q(G)$ descend comme une section sur $S_0$, 
\[\widetilde\Ha_b(G_0) \dans H^0\left(S_0,\det\left(\omega_{G_0^D,b}\right)^{\otimes(p^f-1)}\right).\qedhere\] 
\edem

\section{Changement de Base}
\label{sect6}

Soit $S,S'$ deux schémas lisses sur $\Spec(\kappa_E)$, et un morphisme,
\[ \pi : S' \fleche S.\]
Soit $G/S$ un groupe de Barsotti-Tate, muni d'une action de $\mathcal O$, de signature 
$(p_\tau,q_\tau)_\tau$. 
On a défini son invariant de Hasse partiel (associé à une abscisse $q$) comme un morphisme,
\[ \widetilde\Ha_b(G) : \det(\omega_{G^D,b}) \fleche \det(\omega_{G^D,b})^{\otimes p^f},\]
ou de manière équivalente comme une section sur $S$ de 
$\det(\omega_{G^D,b})^{\otimes (p^f-1)}$.
Maintenant on peut se demander si cette construction est compatible au changement de base, c'est-à-dire si on note $G' = G \times_S S'$, peut-on relier son invariant de Hasse à celui de $G$ ?

\thr
\label{chgtbase}
La construction de l'invariant de Hasse partiel est fonctorielle (ou encore compatible au changement de base), c'est-à-dire que l'on a l'égalité,
\[ \widetilde\Ha_b(G') = \pi^*\widetilde\Ha_b(G).\]
\ethr

\dem
En fait cela découle  simplement de la compatibilité au changement de base du cristal de 
Dieudonné,
\[ \mathbb D(G\times_S S') = \pi_{CRIS}^*\mathbb D(G).\]
On peut supposer (par la propriété de descente du chapitre si dessus) que quitte à étendre les 
scalaires à $\kappa_F$, $S$ et $S'$ sont des $\Spec(\kappa_F)$-schémas.
De même, il est formel que $V_{G'},F_{G'}$ sur le cristal précédent soit les tirées en arrière par $\pi_{CRIS}$ de $V,F$ sur 
$\mathbb D(G)$, et l'unicité de la décomposition suivant l'action de $\mathcal O$ assure que,
\[ \mathbb D(G')_\tau = \pi_{CRIS}^*\mathbb D(G)_\tau, \quad \forall \tau \in \mathcal I.\]
De plus, la formation de $\omega_{G}$ est compatible au changements de base (\cite{BBM} p134), et donc comme $\Fil \bigwedge^{q_\tau} \mathcal E_\tau$ est défini par son image dans $\bigwedge^{q_\tau} \mathcal E_\tau$, celui-ci est compatible au changement de base (bien que ce ne soit pas le cas pour $\Fil \mathcal E$ à priori, mais $\im(\pi_{CRIS}^*\Fil \mathcal E \fleche \pi^*_{CRIS}\mathcal E) = \Fil \pi_{CRIS}^*\mathcal E$). Dès lors, par l'unicité du théorème \ref{thruni} (avec les mêmes notations), on a que $\phi_\tau(G') = \pi_{CRIS}^*\phi_\tau(G)$, et donc
\[\widetilde{\Ha_\tau(G')} = \pi^*\widetilde{\Ha_\tau(G)}, \quad \forall \tau \in \mathcal I. \qedhere\]
\edem

\rem
Cette fonctorialité est encore vraie pour $\Ha_\tau$, pour tout $\tau \dans \mathcal I$, dès lors que $S,S'$ sont au-dessus de $\kappa_F$,
 comme le montre la démonstration.
\erem

On en déduit alors la proposition suivante,

\pro
Soit $S$ un $\kappa_F$-schéma lisse. 
Soit $G/S$ un groupe de Barsotti-Tate tronqué d'échelon $r$.
Si $\tau \in \mathcal I$ est un plongement tel que $q_\tau = \Ht_{\mathcal O}(G)$, et 
\[ r > k_\tau,\] alors $\widetilde\Ha_\tau(G)$ est une section inversible.
\epro

\dem
Supposons $q_\tau = \Ht_{\mathcal O}(G)$, c'est à dire
\[ i_{S/\Sigma*}\omega_{G^D,\tau} \overset{\simeq}{\fleche} 
i_{S/\Sigma*}\left(i_{S/\Sigma}^*\mathcal E_\tau\right) = \mathcal E_\tau/J_{S/\Sigma*}\mathcal E_\tau.\]
Il suffit de vérifier que $\widetilde\Ha_\tau(G)$ est inversible en chaque point géométrique de $S$, par la proposition précédente. Supposons donc $k$ algébriquement clos et $G/k$. Relevons alors le $\mathcal O$-module $p$-divisible tronqué $G$ en un $\mathcal O$-module $p$-divisible (cf. \cite{Wed2} Proposition 3.2), noté encore $G$. Comme $\Ha_\tau(G)$ ne dépend que de sa $p^{k_\tau +1}$ torsion, on peut relier $\Ha_\tau(G[p^r])$ à la section \ref{sect3}.
Mais alors on utilise l'isomorphisme de \cite{BBM} Chapitre 4 (rappelé en \ref{proBBMFON}),
\[ \mathbb D(G)_{(W(k) \fleche k)} \simeq M(G)^{(\sigma)}.\]
Il suffit donc de travailler avec le module de Dieudonné $M(G)$, qui est un $W_r(k) = W(k)/p^r$ module libre.
Mais un calcul direct de l'indice entre les deux réseaux, donne, comme $q_\tau = h$,
\begin{IEEEeqnarray*}{ccc} 
[M(G)_\tau:V^f(M(G)_\tau)] &= &\sum_{i = 1}^f [V^i(M(G)_{\sigma^{-i}\tau}):V^{i+1}(M(G)_{\sigma^{-i-1}\tau})] \\
& = & \sum_{i = 1}^f p_{\sigma^{-i}\tau} \\
& = & \sum_{i = 1}^f h - q_{\sigma^{-i}\tau}\\
& = & k_\tau
\end{IEEEeqnarray*}
Donc, 
\[ [\bigwedge^hM(G)_\tau: V^f(\bigwedge^h M(G)_\tau)] = k_\tau,\]
et donc si $\phi_\tau = \frac{1}{p^{k_\tau}}V^f$, on obtient que $\phi_\tau$ est surjective. Donc,
\[ \phi_\tau : \bigwedge^h \omega_{G^D,\tau} \overset{\simeq}{\fleche} \bigwedge^h \omega_{G^D,\tau}^{(p^f)}.\qedhere\]
\edem

\section{Compatibilité au produit}
\label{sect7}

En général on ne peut pas espérer que les invariants de Hasse partiels précédemment définis soient compatibles au produit. C'est-à-dire que l'on n'a pas en général, comme le montre l'exemple suivant, pour $G,G'$ deux groupes de Barsotti-Tate (tronqués ou non) avec action de $\mathcal O$, que,
\[ \widetilde\Ha_\tau(G\times G') = \widetilde\Ha_\tau(G) \otimes \widetilde\Ha_\tau(G').\]

\exe
\label{exeprod}

Soit $k$ corps algébriquement clos de caractéristique $p$. Soit $\mathcal O = \ZZ_{p^2}$. Notons $\tau'$ et 
$\tau$ les plongements de $\mathcal O$ dans $W(k)$. On décompose tous les cristaux sous la forme $M = M_{\tau'}\oplus M_{\tau}$ et on écrira les matrices en respectant cette décomposition.
Soit $G_1/k$ un groupe $p$-divisible donné par le cristal $(W(k)^4,V_1)$ :
\[
V_1 = \left(
\begin{array}{cc}
0   &   
\left(
\begin{array}{cc}
 1 & 0     \\
  0 &   1   
\end{array}
\right)
\\
\left(
\begin{array}{cc}
 1 & 0   \\
 0 & p  
\end{array}
\right)
 &  0
\end{array}
\right)
\]
C'est-à-dire que $G_1$ est $\mu$-ordinaire et  
\[q^1_\tau = 1, \quad q^1_{\tau'} = 2, \quad \text{et} \quad \Ht_{\mathcal O}(G_1) = 2.\]
Soit aussi $G_2$ un groupe $p$-divisible sur $k$ donné par le cristal $(W(k)^6,V_2)$,

\[
V_2 = \left(
\begin{array}{cc}
0   &   
\left(
\begin{array}{ccc}
 p & 0 & 0     \\
  0 &   p   & 0 \\
  0 & 0 & p
\end{array}
\right)
\\
\left(
\begin{array}{ccc}
 1 & 0 & 0   \\
 0 &1 & 0\\
 0 & 0 & p  
\end{array}
\right)
 &  0
\end{array}
\right)
\]
C'est-à-dire que $G_2$ est $\mu$-ordinaire et  
\[q^2_\tau = 2, \quad q^2_{\tau'} = 0, \quad \text{et} \quad \Ht_{\mathcal O}(G_2) = 3.\]

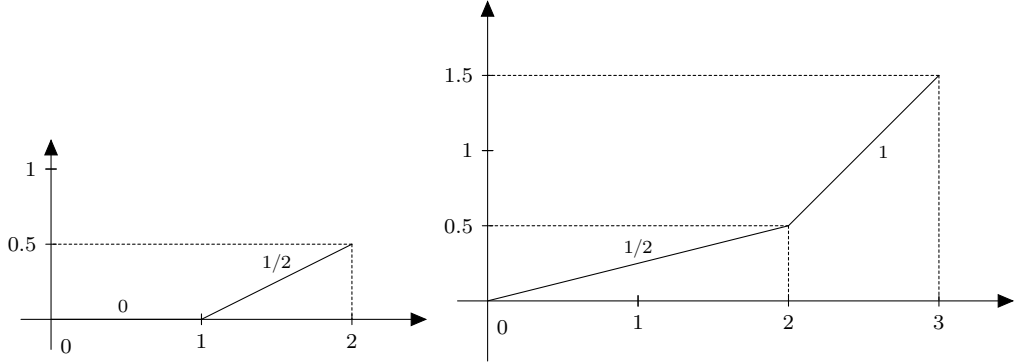
\begin{figure}[h]
\caption{$\mathcal O$-Polygones de Hodge et Newton de $G_1$ et $G_2$.}
\begin{tikzpicture}[line cap=round,line join=round,>=triangle 45,x=2.0cm,y=2.0cm]
\draw[->,color=black] (-0.2,0) -- (2.5,0);
\foreach \x in {,1,2}
\draw[shift={(\x,0)},color=black] (0pt,2pt) -- (0pt,-2pt) node[below] {\footnotesize $\x$};
\draw[->,color=black] (0,-0.2) -- (0,1.2);
\foreach \y in {,0.5,1}
\draw[shift={(0,\y)},color=black] (2pt,0pt) -- (-2pt,0pt) node[left] {\footnotesize $\y$};
\draw[color=black] (0pt,-10pt) node[right] {\footnotesize $0$};
\clip(-0.2,-0.2) rectangle (2.5,1.2);
\draw (0,0)-- (1,0);
\draw (1,0)-- (2,0.5);
\draw [dash pattern=on 1pt off 1pt] (2,0.5)-- (2,0);
\draw [dash pattern=on 1pt off 1pt] (2,0.5)-- (0,0.5);
\begin{scriptsize}
\draw[color=black] (0.48,0.09) node {$0$};
\draw[color=black] (1.5,0.37) node {$1/2$};
\end{scriptsize}
\end{tikzpicture}
\begin{tikzpicture}[line cap=round,line join=round,>=triangle 45,x=2.0cm,y=2.0cm]
\draw[->,color=black] (-0.2,0) -- (3.5,0);
\foreach \x in {,1,2,3}
\draw[shift={(\x,0)},color=black] (0pt,2pt) -- (0pt,-2pt) node[below] {\footnotesize $\x$};
\draw[->,color=black] (0,-0.4) -- (0,2);
\foreach \y in {,0.5,1,1.5}
\draw[shift={(0,\y)},color=black] (2pt,0pt) -- (-2pt,0pt) node[left] {\footnotesize $\y$};
\draw[color=black] (0pt,-10pt) node[right] {\footnotesize $0$};
\clip(-0.2,-0.4) rectangle (3.5,2);
\draw (0,0)-- (2,0.5);
\draw (2,0.5)-- (3,1.5);
\draw [dash pattern=on 1pt off 1pt] (2,0.5)-- (0,0.5);
\draw [dash pattern=on 1pt off 1pt] (2,0.5)-- (2,0);
\draw [dash pattern=on 1pt off 1pt] (3,1.5)-- (0,1.5);
\draw [dash pattern=on 1pt off 1pt] (3,1.5)-- (3,0);
\begin{scriptsize}
\draw[color=black] (2.63,0.99) node {$1$};
\draw[color=black] (1,0.35) node {$1/2$};
\end{scriptsize}
\end{tikzpicture}
\end{figure}

On peut alors donner la matrice de Frobenius du cristal de $G_1\times G_2$ : 
\[
V = \left(
\begin{array}{cc}
0   &   
\left(
\begin{array}{ccccc}
 1 & & &&0     \\
    &  1& && \\
   &  & p&& \\
  & & &p& \\
  0 & & && p
\end{array}
\right)
\\
\left(
\begin{array}{ccccc}
 1 & & &&0     \\
    &  p& && \\
   &  &1&& \\
  & & &1& \\
  0 & & & &p
\end{array}
\right)
 &  0
\end{array}
\right)
\]

\begin{figure}[h]
\begin{center}
\caption{$\mathcal O$-Polygones de Hodge et Newton de $G = G_1\times G_2$.}
\label{figHN5}
\definecolor{ffqqqq}{rgb}{1,0,0}
\begin{tikzpicture}[line cap=round,line join=round,>=triangle 45,x=2.0cm,y=2.0cm]
\draw[->,color=black] (-0.2,0) -- (5.2,0);
\foreach \x in {,1,2,3,4,5}
\draw[shift={(\x,0)},color=black] (0pt,2pt) -- (0pt,-2pt) node[below] {\footnotesize $\x$};
\draw[->,color=black] (0,-0.2) -- (0,3.7);
\foreach \y in {,0.5,1,1.5,2,2.5,3,3.5}
\draw[shift={(0,\y)},color=black] (2pt,0pt) -- (-2pt,0pt) node[left] {\footnotesize $\y$};
\draw[color=black] (0pt,-10pt) node[right] {\footnotesize $0$};
\clip(-0.2,-0.2) rectangle (5.2,3.7);
\draw [color=ffqqqq] (0,0)-- (2,0);
\draw [color=ffqqqq] (3,0.5)-- (2,0);
\draw [color=ffqqqq] (3,0.5)-- (5,2.5);
\draw (0,0)-- (1,0);
\draw (1,0)-- (4,1.5);
\draw (4,1.5)-- (5,2.5);
\draw (0.53,3.19) node[anchor=north west] {$Newt_{\mathcal O}(G)$};
\draw (0.53,2.98) node[anchor=north west] {$Hdg_{\mathcal O}(G)$};
\draw (2,3.03)-- (1.5,3.03);
\draw [color=ffqqqq] (2,2.82)-- (1.5,2.82);
\draw [dash pattern=on 2pt off 2pt] (4,1.5)-- (4,0);
\draw [dash pattern=on 2pt off 2pt] (4,1.5)-- (0,1.5);
\draw [dash pattern=on 2pt off 2pt] (3,0.5)-- (3,0);
\draw [dash pattern=on 2pt off 2pt] (3,0.5)-- (0,0.5);
\draw [dash pattern=on 2pt off 2pt] (5,2.5)-- (5,0);
\draw [dash pattern=on 2pt off 2pt] (5,2.5)-- (0,2.5);
\begin{scriptsize}
\draw[color=black] (4.5,2.2) node {$1$};
\draw[color=black] (2.5,1) node {$1/2$};
\draw[color=black] (0.48,0.09) node {$0$};
\end{scriptsize}
\end{tikzpicture}
\end{center}
\end{figure}
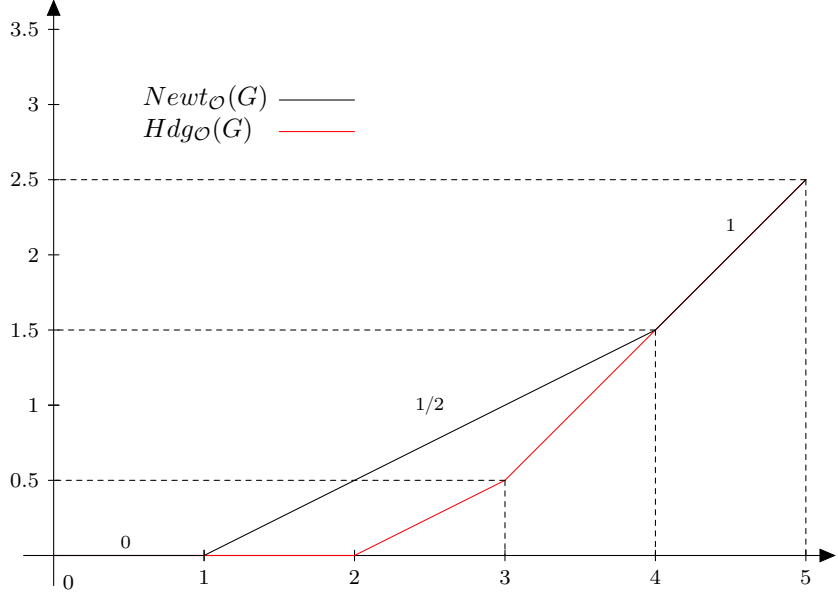
En particulier il n'est pas $\mu$-ordinaire (voir la caractérisation de Moonen, théorème \ref{thrmoo}). D'ailleurs on le voit sur le calcul de (la valuation de) l'invariant de Hasse associé à $\tau$ :
\[ \Ha_\tau(G_1\times G_2) = 1, \quad \Ha_\tau(G_1) = \Ha_\tau(G_2) = 0.\]
Mais on le voit aussi sur les polygones de Newton (figure \ref{figHN5}).
\eexe

Afin d'espérer avoir une égalité entre le $\mu$-invariant d'un produit et le produit des $\mu$-invariants, étant donné que ceux-ci impliquent tous deux de diviser par une bonne puissance de $p$, il est naturel de supposer que l'on divise par la même puissance. C'est à dire que si 
\[G = G_1 \times G_2,\]
et que l'on note $(q_\tau)$ la signature de $G$, $(q_\tau^1), (q_\tau^2)$ celles de $G_1$ et $G_2$, et pour tout $\tau \in \mathcal I$, $k_\tau$, (respectivement $k_\tau^1, k_\tau^2$) les puissances de $p$ par lequelles on divise dans la construction de l'invariant de Hasse associé à $\tau$ (cf théorème \ref{thruni}) de $G$ (respectivement $G_1,G_2$). Alors une hypothèse raisonnable est que
\[ k_\tau = k_\tau^1 + k_\tau^2.\]

\rem
\label{rem62}
Soit $G/\Spec(k)$ un groupe $p$-divisible avec $\mathcal O$-action. Si les $\mathcal O$-polygones de Hodge et Newton de $G$ se touchent en une abscisse $x$ qui est un point de rupture pour le $\mathcal O$-polygone de Newton, on a alors par le théorème \ref{thrHN} de Mantovan-Viehmann, une décomposition de Hodge-Newton,
\[ G = G_1 \times G_2.\]
Et d'après la remarque \ref{rem29}, on connaît les signatures de $G_1$ et $G_2$, on en déduit donc que, pour tout $\tau \in \mathcal I$,
\[ k_\tau = k_\tau^1 + k_\tau^2.\]
\erem

\pro
\label{proproduit}
Soit $S$ un schéma lisse sur $\Spec(\kappa_E)$. Soit $G_1,G_2$ deux groupes de Barsotti-Tate tronqués de rang $r$, munis d'une action de $\mathcal O$, et de signature 
$(p_\tau^1,q_\tau^1)_\tau$ et $(p_\tau^2,q_\tau^2)_\tau$. 
La signature de $G = G_1\times G_2$ est alors donnée par,
\[p_\tau = p_\tau^1 + p_\tau^2, \quad \text{et} \quad q_\tau = q_\tau^1 + q_\tau^2, \quad \forall \tau \in \mathcal I.\]
On fait l'hypothèse suivante,
\[ k_\tau := \sum_{\tau' | q_{\tau'} \leq q_\tau} q_\tau - q_{\tau'} = k_\tau^1 + k_\tau^2.\]
Supposons $r > k_\tau^1 + k_\tau^2$ .Notons $G = G_1 \times G_2$.
Alors, quitte à étendre les scalaires de $S$, pour tout $\tau \dans \mathcal I$, on a un isomorphisme,
\[ \det(\omega_{G^D,\tau}) \simeq \det(\omega_{G_1^D,\tau}) \otimes \det(\omega_{G_2^D,\tau}),\]
et sous cet isomorphisme,
\[\widetilde\Ha_\tau(G_1\times G_2) = \widetilde\Ha_\tau(G_1) \otimes \widetilde\Ha_\tau(G_2).\]
En particulier, si 
\[ k_\tau = k_\tau^1 + k_\tau^2, \quad \forall \tau \in \mathcal I,\]
alors,
\[ \widetilde{^\mu\Ha}(G_1\times G_2) = \widetilde{^\mu\Ha}(G_1) \otimes  \widetilde{^\mu\Ha}(G_2).\] 
\epro

\dem
Tout d'abord, quitte à faire un changement de base, supposons $S$ au-dessus de $\kappa_F$. L'égalité,
\[ \omega_G = \omega_{G_1} \oplus \omega_{G_2},\]
induit l'isomorphisme voulu,
\[ \det(\omega_{G,\tau}) \simeq \det(\omega_{G_1,\tau}) \otimes \det(\omega_{G_2,\tau}).\]
Notons $\mathcal E = \mathcal Ext^1(G^D,\mathcal O_{S/\Sigma})$ et 
$\mathcal E^1,\mathcal E^2$ les cristaux en $\mathcal O_{S/\Sigma}/p^r$-modules 
localement libres associés à $G_1$ et $G_2$. On a alors,
\[ \mathcal E = \mathcal E^1 \oplus \mathcal E^2, \quad \text{et} \quad V := V_{\mathcal E} = V_{\mathcal E^1} \oplus V_{\mathcal E^2}.\]
On a alors les isomorphismes,
\[ \displaystyle\bigwedge^{q_{\tau}} \mathcal E_{\tau} \pmod{p^{r-k_\tau}} = \bigoplus_{j=0}^{q_\tau} \left(\displaystyle\bigwedge^{j} \mathcal E^1_{\tau} \otimes 
\displaystyle\bigwedge^{q_\tau - j} \mathcal E^2_{\tau}\right) \pmod{p^{r-k_\tau}},\]
et \[ \Fil\displaystyle\bigwedge^{q_{\tau}} \mathcal E_{\tau} \pmod{p^{r-k_\tau}} = \bigoplus_{j=0}^{q_\tau} \left(\Fil\displaystyle\bigwedge^{j} \mathcal E^1_{\tau} \otimes \Fil\displaystyle\bigwedge^{q_\tau - j} \mathcal E^2_{\tau}\right) \pmod{p^{r-k_\tau}}.\]
De plus, on a un diagramme commutatif de morphismes de cristaux,
\begin{center}
\begin{tikzpicture}[description/.style={fill=white,inner sep=2pt}] 
\matrix (m) [matrix of math nodes, row sep=3em, column sep=2.5em, text height=1.5ex, text depth=0.25ex] at (0,0)
{ 
\displaystyle\bigwedge^{q_{\tau}} \mathcal E_{\tau} &  & \bigwedge^{q_{\tau}} \mathcal E_{\tau}^{(p^f)} \\
 \displaystyle\bigwedge^{q_\tau^1} \mathcal E^1_{\tau} \otimes \displaystyle\bigwedge^{q_\tau^2} \mathcal E^2_{\tau} & & \displaystyle\bigwedge^{q_\tau^1} \mathcal E^{1,(p^f)}_{\tau} \otimes \displaystyle\bigwedge^{q_\tau^2} \mathcal E^{2,(p^f)}_{\tau}  \\
 };

\path[->,font=\scriptsize] 
(m-1-1) edge node[auto] {$V^f$} (m-1-3)
(m-2-1) edge node[auto] {$V_1^f \wedge V_2^f$} (m-2-3)
;
\path[right hook->,font=\scriptsize] 
(m-2-1) edge node[auto] {} (m-1-1)
(m-2-3) edge node[auto] {} (m-1-3)
;
\end{tikzpicture}
\end{center}
dont les cristaux sont des cristaux en $\mathcal O_{S/\Sigma}$-modules localement libres, et dont les flèches horizontales sont divisibles par $p^{k_\tau} = p^{k_\tau^1}\cdot p^{k_\tau^2}$ par hypothèse, de manière unique modulo $p^{r-k_\tau}$. On en déduit la commutativité du diagramme,
 \begin{center}
\begin{tikzpicture}[description/.style={fill=white,inner sep=2pt}] 
\matrix (m) [matrix of math nodes, row sep=3em, column sep=2.5em, text height=1.5ex, text depth=0.25ex] at (0,0)
{ 
\Fil\displaystyle\bigwedge^{q_{\tau}} \mathcal E_{\tau} \pmod{p^{r-k_\tau}} &  & \Fil\bigwedge^{q_{\tau}} \mathcal E_{\tau}^{(p^f)}\pmod{p^{r-k_\tau}} \\
\left( \Fil\displaystyle\bigwedge^{q_\tau^1} \mathcal E^1_{\tau}\right) \otimes \left(\Fil\displaystyle\bigwedge^{q_\tau^2} \mathcal E^2_{\tau}\right)\pmod{p^{r-k_\tau}} & & \left(\Fil\displaystyle\bigwedge^{q_\tau^1} \mathcal E^{1,(p^f)}_{\tau}\right) \otimes\left(\Fil\displaystyle\bigwedge^{q_\tau^2} \mathcal E^{2,(p^f)}_{\tau}\right) \pmod{p^{r-k_\tau}} \\
 };

\path[->,font=\scriptsize] 
(m-1-1) edge node[auto] {$\phi_\tau$} (m-1-3)
(m-2-1) edge node[auto] {$\phi_\tau^1 \otimes \phi_\tau^2$} (m-2-3)
;
\path[right hook->,font=\scriptsize] 
(m-2-1) edge node[auto] {} (m-1-1)
(m-2-3) edge node[auto] {} (m-1-3)
;
\end{tikzpicture}
\end{center}
où les flèches $\phi_\tau, \phi_\tau^1,\phi_\tau^2$ sont données par le théorème $\ref{thruni}$ appliqué à $G$ (respectivement $G_1,G_2$).
Mais en passant au quotient à 
$\omega_{G^D,\tau} = \omega_{G_1^D,\tau} \oplus \omega_{G_2^D,\tau}$, 
les monomorphismes verticaux deviennent des isomorphismes, on en déduit donc que,
\[\widetilde\Ha_\tau(G) = \widetilde\Ha_\tau(G_1) \otimes \widetilde\Ha_\tau(G_2).\qedhere\]
\edem

\section{Cas d'une base quelconque}
\label{sect8}

Dans les sections précédentes, on a construit des invariants de Hasse partiels 
(et $\mu$-ordinaire) 
associés à chaque groupe de Barsotti-Tate (muni d'une action de $\mathcal O = \mathcal O_F$) sur une base 
lisse. On aimerait étendre la construction à tous les groupes de Barsotti-Tate avec 
$\mathcal O$-structure, quelque soit la base. Pour cela, il suffit de le faire pour le groupe de 
Barsotti-Tate universel sur le champ des $\mathcal O$-modules de Barsotti-Tate tronqués, c'est à dire associé au foncteur $\mathcal{BT}$ qui à tout schéma $S$ sur 
$\kappa_E$ 
associe l'ensemble des classes d'isomorphisme de $\mathcal O$-modules de Barsotti-Tate tronqués (d'échelon $r$). 
Malheureusement ce foncteur n'est pas représentable par un 
schéma, mais seulement par un champ. Par chance, ce champ est lisse. On peut donc construire 
grâce à la section précédente les invariants de Hasse partiels sur une présentation de celui-ci, et 
on aimerait ensuite descendre ces sections. C'est l'objet de cette partie et de l'appendice \ref{appA}.

Rappelons que l'on s'intéresse aux groupes de Barsotti-Tate tronqués, d'échelon disons $r$, munis d'une action de $\mathcal O = \mathcal O_F$, où $F/\QQ_p$ est une extension finie 
non ramifiée de degré $f$. Rappelons aussi que l'on a noté,
\[ \mathcal I = \Hom(\kappa_F,\overline{\FP}) = \Hom_{\QQ_p}(F,\overline{\QQ_p}).\]
Notons $\mathcal F = \NN^{\mathcal I}$ l'ensemble des applications de $\mathcal I$ dans $\NN$. 
Cet ensemble est naturellement muni d'une action de $\Gal(\overline{\QQ_p}/\QQ_p)$ grâce à son
action sur les plongements. Soit $(p_\tau)_{\tau \in \mathcal I} \dans \mathcal F$ une telle 
application. 
\defi
On appelle \textit{corps réflexe associé à} $(p_\tau)_{\tau \in \mathcal I}$ l'extension de $\QQ_p$ 
associée au sous-groupe stabilisateur de $(p_\tau)_{\tau \in \mathcal I}$. On le notera
généralement $E$.
\edefi

\subsection{Descente au champ $\mathcal{BT}$}

On munit la catégorie $Sch/\Spec(\FP)$ de la topologie lisse.
Notons $\mathcal{BT}_{r,\mathcal O}$ le champ sur $\Spec(\FP)$ des groupes de Barsotti-Tate 
tronqués d'échelon $r$, munis d'une action de $\mathcal O$. On peut alors décomposer ce champ en sous-champs ouverts, suivant la $\mathcal O$-hauteur,
\[ \mathcal{BT}_{r,\mathcal O} = \coprod_{h} \mathcal{BT}_{r,\mathcal O,h}.\]

\defi Pour $h$ fixé, étant donné une fonction $(p_\tau)_\tau$, on note $q_\tau = h-p_\tau$. 
On définit alors le sous-champ de ${\mathcal{BT}_{r,\mathcal O,h}}_{|Sch/\Spf(\mathcal O/p)}$ suivant, 
\[{\mathcal{BT}_{r,\mathcal O,h,(p_\tau,q_\tau)}}_{\tau \in \mathcal I},\]
 qui a un schéma $S/\Spec(\mathcal O/p\mathcal O)$ associe les classes 
d'isomorphismes $(G,\iota)$ de $S$-groupes de Barsotti-Tate tronqués d'échelon $r$, 
hauteur $h$, munis d'une action de $\mathcal O$, tels que,
\[ \dim_{\mathcal O_S} \omega_{G,\tau}
 = p_\tau, \quad
\forall \tau \in \mathcal I.\] 

Ce sous-champ descend naturellement à $\Spec(\mathcal O_E/p\mathcal O_E)$ où $E$ est le corps réflexe de $(p_\tau)_{\tau \in \mathcal I}$.
Définissons aussi 
\[X_{r,\mathcal O,h,(p_\tau,q_\tau)},\]
le champ sur $Sch/\Spec(\mathcal O_E/p\mathcal O_E)$ qui à un schéma $S$ associe les classes 
d'isomorphismes $(G,\iota,\alpha)$ où $(G,\iota)$ est un objet de 
$\mathcal{BT}_{r,h, \mathcal O,(p_\tau,q_\tau)}(S)$ et $\alpha$ est une rigidification, i.e. un isomorphisme,
\[ \alpha : \mathcal O_S^{p^{rhf}} \overset{\simeq}{\fleche} \pi_*(\mathcal O_G),\]
où $\pi : G \fleche S$ est le morphisme structural. 
\edefi

On a une action de $GL_{p^{rhf}}$ sur $X_{r,\mathcal O,h,(p_\tau,q_\tau)}$, et un morphisme d'oubli, $GL_{p^{rhf}}$-équivariant,
\[ X_{r,\mathcal O,h,(p_\tau,q_\tau)} \fleche \mathcal{BT}_{r,\mathcal O,h,(p_\tau,q_\tau)}.\]

D'après la proposition 1.8 de \cite{Wed2}, $X_{r,\mathcal O,h,(p_\tau,q_\tau)}$ est représentable par un sous-schéma 
ouvert d'un schéma quasi-affine de type fini (et est donc quasi-affine de type fini) au-dessus de 
$\Spec(\mathcal O_E)$.
De plus, le morphisme précédent induit un isomorphisme,
\[ [GL_{p^{rhf}}\backslash X_{r,\mathcal O,h,(p_\tau,q_\tau)}] \simeq \mathcal{BT}_{r,\mathcal O,h,(p_\tau,q_\tau)}.\]
Cela fait donc de $\mathcal{BT}_{r,\mathcal O,h,(p_\tau,q_\tau)}$ un champ algébrique, qui en plus est lisse d'après 
le corollaire 3.3 de \cite{Wed2}. 

\rem
En fait il existe un champ sur $\Spec(\ZZ_p)$, $\mathcal{BT}_{r,\mathcal O,h}$ qui est aussi algébrique lisse (\cite{Wed2}) dont le champ précédent est (un ouvert de) la réduction
 modulo $p$. À moins d'accepter de parler de champ formel sur $\Spf(\ZZ_p)$, on ne peut cependant pas découper le champ 'en caractéristique 0' selon la signature.
\erem

Pour la suite, on abrégera $X_{r,\mathcal O,h,(p_\tau,q_\tau)}$ en $X$ et $\mathcal{BT}_{r,\mathcal O,h,(p_\tau,q_\tau)}$ en $\mathcal{BT}_r$, en espérant que cela n'entraînera pas de confusion : tous les groupes de Barsotti-Tate considérés seront tronqués d'échelon $r$, de hauteur $h$, et auront pour signature $(p_\tau,q_\tau)$.

Soit $G \fleche \mathcal{BT}_r$ le groupe de Barsotti-Tate universel, muni d'une action de $\mathcal O$, 
et soit $G'$ le produit fibré $G\times_{\mathcal{BT}_r} X$ qui est un schéma en groupe sur $X$, et même un Barsotti-Tate tronqué d'échelon $r$.

\thr
\label{thrdes}
Soit $q \dans \NN$ tel qu'il existe $\tau \dans \mathcal I$ vérifiant $q = q_\tau$. Supposons $r > k_\tau$. 
Alors il existe sur le champ $\mathcal{BT}_r$ un fibré en droite, noté
\[ \det(\omega_{G,q})^{\otimes(p^f-1)},\]
et une section globale $\widetilde\Ha_q$ de ce faisceau sur $\mathcal{BT}_r$, tel que pour tout $\Spec(\kappa_E)$-schéma $U$ lisse, et tout $U$-groupe de Barsotti-Tate tronqué d'échelon 
$r$ $G_u$ muni d'une action de $\mathcal O$, de $\mathcal O$-hauteur $h$, et de signature $(p_\tau)_{\tau \in \mathcal I}$, associé à un $u \dans {\mathcal{BT}_r}_U$, on ait la compatibilité 
avec la définition des sections précédentes,
\[ \widetilde\Ha_q(G_u) = u^*\widetilde\Ha_q(G).\]
En particulier cela définit un invariant de Hasse partiel $\widetilde\Ha_q$, associé à $q$, à tout groupe de Barsotti-Tate muni d'une $\mathcal O$-action, quelque soit la base, et donc un $\mu$-invariant de Hasse $\widetilde{^\mu\Ha}$, par produit.
\ethr

\dem
L'existence du faisceau quasi-cohérent $\det(\omega_{G,q})^{\otimes(p^f-1)}$ provient 
simplement de la fonctorialité de ce faisceau par changement de base, et du fait que $G$ 
est l'objet universel.
Maintenant on peut voir une section de ce faisceau comme un morphisme,
\[ \widetilde\Ha_q : \mathcal O_{\mathcal{BT}_r} \fleche \det(\omega_{G,q})^{\otimes(p^f-1)}.\]
Par le lemme \ref{lemqcoh}, comme $\mathcal{BT}_r$ est un champ algébrique lisse, il suffit donc de voir que pour tout $U$ lisse sur $\Spec(\kappa_E)$, et tout $u \dans {\mathcal{BT}_r}_U$, il existe un morphisme,
\[ (\widetilde\Ha_q)_u : \mathcal O_U \fleche u^*\det(\omega_{G,q})^{\otimes(p^f-1)} = 
\det(\omega_{u^*G,q})^{\otimes(p^f-1)},\]
vérifiant les compatibilités de la définition \ref{defqcoh}. On a construit dans ce cas une section 
$\widetilde\Ha_q(u^*G) \dans \Gamma(U,\det(\omega_{u^*G,q})^{\otimes(p^f-1)})$, c'est-à-dire un morphisme,
\[\widetilde\Ha_q(u^*G) : \mathcal O_U \fleche u^*\det(\omega_{G,q})^{\otimes(p^f-1)}.\]
La compatibilité pour $U' \overset{f}{\fleche} U$, provient de la compatibilité au changement 
de base de $\widetilde\Ha_q$ pour une base lisse, c'est à dire du théorème \ref{chgtbase}, qui assure la commutativité du diagramme,
 \begin{center}
\begin{tikzpicture}[description/.style={fill=white,inner sep=2pt}] 
\matrix (m) [matrix of math nodes, row sep=3em, column sep=2.5em, text height=1.5ex, text depth=0.25ex] at (0,0)
{ 
\mathcal O_{U'}  &  & \det(\omega_{(u\circ f)^*G,q})^{\otimes(p^f-1)} \\
f^*\mathcal O_U  &  & f^* \det(\omega_{u^*G,q})^{\otimes(p^f-1)} \\ 
 };

\path[->,font=\scriptsize] 
(m-1-1) edge node[auto] {$\widetilde\Ha_q((u\circ f)^*G)$} (m-1-3)
(m-1-1) edge node[auto] {$\simeq$} (m-2-1)
(m-2-1) edge node[auto] {$f^*\widetilde\Ha_q(u^*G)$} (m-2-3)
(m-1-3) edge node[auto] {$\simeq$} (m-2-3)
;
\end{tikzpicture}
\end{center}
\qedhere
\edem 

\subsection{$\mathcal O$-modules de Lubin-Tate et groupes $p$-divisibles $\mu$-ordinaires}

Grâce à la compatibilité au changement de base, et à la compatibilité entre inversibilité de $\widetilde{^\mu\Ha}$ et $\mu$-ordinarité vérifiée dans la section \ref{sect3} (proposition
\ref{propolygone}), on peut en particulier utiliser les résultats de Moonen, \cite{Moo}. On va introduire ici quelques notations.

Soit $K/\QQ_p$ une extension finie non ramifiée, et $\mathcal O := \mathcal O_K$ son anneau d'entiers, $\pi$ 
une uniformisante, et $\kappa_K$ son corps résiduel.
Notons $\mathcal I = \Hom(\kappa_K,\overline{\FP})$.

\defi
Soit $\tau \dans \mathcal I$. Soit $R$ une $\mathcal O$-algèbre $p$-adiquement complète et séparée. Il existe un groupe $p$-divisible muni d'une action de 
$\mathcal O$, noté $\mathcal{LT}_\tau$, de hauteur $[K:\QQ_p]$, tel que son module des différentielles
$\omega_{\mathcal{LT}_\tau}$ soit de dimension 1, et sur lequel l'action de $\mathcal O$ soit donnée par $\rho \circ \tau : \mathcal O \fleche R$ où $\rho$ est 
le morphisme structural de $R$. 
On appellera aussi un tel groupe (pour lequel l'action sur l'algebre de Lie se fait via un seul plongement) un $\mathcal O$-module $\pi$-divisible.
\edefi

\rem
Lorsque $R$ est un $\mathcal O$-module sans torsion (et $\tau : \mathcal O \fleche R$ structural, mais on peut en general regarder $\mathcal O \otimes_{\tau, \mathcal O} R$), 
c'est le groupe $p$-divisible associé au display donné par la proposition 29 de \cite{Zink}. 
Il dépend à priori du choix d'une uniformisante de $K$, mais quitte à étendre les scalaires à $\mathcal O_{\widehat{K^{nr}}}$, cette ambiguité disparait.
\erem

\defi
\label{LTA}
On reprend les notations précédentes, supposons que $K/\QQ_p$ soit non ramifiée. Soit $A \subset \mathcal I$ un sous-ensemble. On note alors,
\[ \mathcal{LT}_A,\]
le $\mathcal O$-module $p$-divisible dont le display est donné par le produit tensoriels des displays des différents $\mathcal{LT}_\tau$ pour $\tau \in A$, c'est à dire,
\[ \mathcal P(\mathcal{LT}_A) := \bigotimes_{\tau \in A}^{W(R)\otimes_{\ZZ_p} \mathcal O} \mathcal P(\mathcal{LT}_\tau).\]
On vérifie que cela définit un display sur $R$, et donc d'après le théorème 9 de \cite{Zink}, $\mathcal{LT}_A$ existe, c'est un $\mathcal O$-module 
$p$-divisible sur $R$ de $\mathcal O$-hauteur 1, de signature $p_\tau = 1$ si et seulement si $\tau \in \mathcal A$.
\edefi

\dem
Dans le cas non ramifié, le display est simplement donné par $P = \bigoplus_\tau W(R)e_\tau, Q = \bigoplus_{\tau \not\in A} W(R)e_\tau \oplus \bigoplus_{\tau \in A} I_Re_\tau$, 
$F$ est donné par \[Fe_\tau = 
\left\{
\begin{array}{cc}
pe_{\sigma\tau} & \text{ si } \tau \not\in A \\
 e_{\sigma\tau} & \text{ sinon}
\end{array}
\right.
\]
L'application $V^{-1} : Q \fleche P$ est donnée par $\frac{1}{p}F$, c'est à dire, si $\tau \in A$, $V^{-1}(^Vxe_\tau) = xe_{\sigma\tau}$, et si $\tau \not\in A$,
$V^{-1}e_\tau = e_{\sigma\tau}$. On vérifie que cela définit bien un display sur $R$, muni d'une action de $\mathcal O$, puisque $P = \mathcal O \otimes_{\ZZ_p} W(R)$.
\edem

\exe
Si $A = \emptyset$, alors par convention $\mathcal{LT}_A = \QQ_p/\ZZ_p\otimes_{\ZZ_p}\mathcal O$. 
Si $A = \mathcal I$, $\mathcal{LT}_A = \mu_{p^\infty} \otimes_{\ZZ_p} \mathcal O$.
\eexe

\rem
Il est probable que l'on n'ait pas besoin de l'hypothèse sur la ramification de $K$, les groupes de Lubin-Tate existant même dans ce cas, mais leur display n'est défini que si 
$R$ est sans torsion (cf. \cite{Zink} Proposition 29). Cela est lié au fait que la signature n'est cependant pas bien défini dans le cas ramifié si $pR = 0$.
\erem

Supposons que $\mathcal O = \mathcal O_F$ pour $F/\QQ_p$ une extension finie non ramifiée.
On peut alors réécrire le théorème 1.3.7 de \cite{Moo}. 

\thr[Moonen]
\label{thrmoo}
Soit $k$ un corps algébriquement clos de caractéristique $p>0$. Soit $G/\Spec(k)$ un groupe $p$-divisible muni d'une action de $\mathcal O$, de $\mathcal O$-hauteur $h$ et de signature $(p_\tau,q_\tau)_{\tau \in \mathcal I}$. On utilise la notation suivante,
\[ \{q_\tau : \tau \dans \mathcal I\} = \{ q^{(1)} < q^{(2)}< \dots < q^{(r)}\}, \quad \text{et} \quad q^{(0)} = 0, q^{(r+1)} = h. \]
Alors les conditions suivantes sont équivalentes,
\begin{enumerate}
\item $G$ est $\mu$-ordinaire.
\item On a un isomorphisme de $\mathcal O$-groupes $p$-divisibles,
\[ G \simeq X^{ord}_{(q_\tau)} := \prod_{l =1}^{r+1} \mathcal{LT}_{A_{l-1}}^{q^{(l)}-q^{(l-1)}},\]
où, pour tout $l$,
\[ A_l = \{ \tau \in \mathcal I : q_\tau \leq q^{(l)}\} \quad \text{et} \quad A_0 = \emptyset\]
\item On a un isomorphisme de $\mathcal O$-schémas en groupes,
\[ G[p] \simeq X^{ord}_{(q_\tau)}[p].\]
\end{enumerate}
\ethr

On montrera un théorème similaire (au sens où l'on donne une forme explicite aux groupes $p$-divisibles $\mu$-ordinaires) sur $\mathcal O_C$ dans la section \ref{sect12}.

\rem
Avec les notations du théorème, on a la suite d'inclusions,
\[ A_0 \subset A_1 \subset \dots \subset A_{r+1}.\]
Réciproquement, pour tous sous-ensembles $(A_l)$ de $\mathcal I$ totalement ordonnés pour l'inclusion, et des entiers $(n_l)$, le produit,
\[\prod_{l =1}^{r+1} \mathcal{LT}_{A_{l-1}}^{n_l},\]
défini un groupe $\mu$-ordinaire pour une certaine donnée $\mu$ (que l'on peut expliciter en fonction de la signature du groupe). C'est faux si les $(A_i)$ ne sont pas ordonnés (voir l'exemple de la section \ref{sect7}). 
\erem

\subsection{Exemple de calcul}

Supposons que $S = \Spec(\mathcal O_C/p\mathcal O_C)$, où $C = \widehat{\overline{\QQ_p}}$.

Soit $G/S$ un groupe $p$-divisible tronqué d'échelon $r$, muni d'une action de $\mathcal O$. D'après le début de cette section, on sait associer à $G$ des invariants de Hasse 
partiels, ainsi qu'un $\mu$-invariant, et ce, même si la base n'est pas lisse. A priori pourtant, la construction (relativement) explicite de la section \ref{sect4} ne s'applique pas à 
strictement parler, puisque $\mathcal O_C/p$ n'est pas lisse.
Notons néanmoins toujours $\Sigma = \Spec(W(\kappa_F))$ et regardons $\mathcal E = \mathbb D(G)$ le cristal sur $\Cris(S/\Sigma)$ associé par \cite{BBM} à $G$.
Le site $\Cris(S/\Sigma)$ possède un objet final, l'épaississement à puissances divisées (mais non $p$-adique) introduit par Fontaine, voir par exemple \cite{Che},
\[ A_{cris} \overset{\theta}{\fleche} \mathcal O_C/p.\]
Notons $\phi$ le Frobenius de $A_{cris}$. C'est un relèvement fort du Frobenius de $\mathcal O_C/p$. 
On peut évaluer le cristal $\mathcal E$ sur cet épaississement, et on obtient,
\[H^0(\Cris(S/\Sigma),\mathcal E) = \mathcal E_{A_{cris} \twoheadrightarrow \mathcal O_C/p} =: E,\]
qui est un $A_{cris}/p^r$ module libre de rang $\Ht(G)$ muni d'applications $V$ et $F$, que l'on décompose,
\[E = \bigoplus_{\tau \in \mathcal I} E_\tau.\]
On peut aussi regarder (voir définition \ref{def32}),
\[\Fil E = \Ker\left(E \fleche E\otimes_{\theta} \mathcal O_C/p \fleche \omega_{G}^\vee\right), \quad \text{et} \quad\Fil E_\tau = \Fil E \cap E_\tau.\]
On a alors que 
\[ V(E_{\tau'}) \subset \Fil (E_{\sigma^{-1}\tau'}^{(\phi)}) + pE_{\sigma^{-1}\tau'}^{(\phi)},\]
et si $q_{\tau'} < q_\tau$, alors,
\[\im(\bigotimes^{q_\tau} \Fil (E_{\tau'}^{(\phi)}) \fleche (\bigwedge^{q_\tau} E_{\tau'})^{(\phi)}) \subset 
\left(\Fil^{q_\tau-q_{\tau'}}A_{cris}\bigwedge^{q_\tau} E_{\tau'}\right)^{(\phi)}\subset 
\bigwedge^{q_\tau} E_{\tau'} \otimes_{A_{cris},\phi} p^{q_\tau-q_{\tau'}}A_{cris}.\]
La dernière inclusion étant due au fait que si $M$ est un $A_{cris}-$module, si $z \dans \Fil^iA_{cris}$, $x \dans M$, alors dans $M^{(\phi)} = M \otimes_{A_{cris},\phi} A_{cris}$,
\[ (zx) \otimes 1 = x \otimes \phi(z), \quad \text{et} \quad \phi(z) \dans p^i A_{cris}.\]
On peut donc, pour chaque $j$ tel que $q_{\sigma^j\tau} < q_\tau$ diviser $V^j$ par $p^{q_\tau - q_{\sigma^j\tau}}$ (voir aussi la proposition \ref{pro34}), et donc il existe,
\begin{equation}
\label{zetaOC}
 \zeta_\tau : \bigwedge^q_\tau E_\tau \fleche  \bigwedge^q_\tau E_\tau\otimes_{A_{cris},\phi^f} A_{cris} \pmod{p^{r-k_\tau}A_{cris}}\end{equation}
tel que $p^{k_\tau}\zeta_\tau = V^f$. De plus comme $A_{cris}$ est sans $p$-torsion, et $E$ étant libre sur $A_{cris}/p^r$, un tel $\zeta_\tau$ est unique, or la construction 
général sur le champ $\mathcal{BT}_r$ du début de cette section nous assure l'existence
d'une telle division de $V^f$, on a donc que $\zeta_\tau \equiv \widetilde\Ha_\tau \pmod{p^{r-k_\tau}}$.

On peut restreindre cette application $\zeta_\tau$ à $\Fil E_\tau$ puis la réduire modulo $\Ker \theta$, et ainsi recalculer effectivement sur le cristal $E$ les invariants de Hasse partiels.

\rem
Dans la proposition \ref{prodefouniv}, on donne une autre situation, celle de l'espace des déformations d'un $\mathcal O$-module $p$-divisible, d'une base non lisse 
(mais formellement lisse) sur laquelle on peut construire directement les invariants de Hasse $\Ha_q$. Dans ce cas, à l'aide des displays de Zink on a une 
formule explicite (en termes des coordonnées sur l'espace des déformations) pour ces invariants.
\erem

\subsection{Diviseurs de Cartier}

D'après l'appendice \ref{sect9}, étant donné que le champ $\mathcal{BT}_r$ est lisse, pour
montrer que les invariants de Hasse partiels définissent des diviseurs de Cartier, il suffit de 
montrer que le complémentaire de leur lieu d'annulation est dense (pour la topologie de Zariski).
Pour cela, on va simplement utiliser le résultat suivant de \cite{Wed}. Rappelons simplement 
nos notations, on fixe $r,h \dans \NN^*$, et $(p_\tau)_{\tau \in \mathcal I}$, et on note 
$\mathcal{BT}_r$ le champ algébrique lisse des groupes $p$-divisible tronqué d'échelon $r$, 
hauteur $h$, munis d'une action de $\mathcal O$, et de signature $(p_\tau)_\tau$. On note 
aussi $\mathcal{BT}_\infty$ le champ des groupes $p$-divisibles de hauteur $h$, 
avec une action de $\mathcal O$ de signature $(p_\tau)_\tau$.
On a un morphisme "points de $p^r$-torsion",
\[ \mathcal{BT}_\infty \overset{[p^r]}{\fleche} \mathcal{BT}_r.\]

Le résultat suivant est le résultat principal de \cite{Wed}, mais il n'est pas cité comme tel dans \cite{Wed}, on le retranscrit ici sous une forme qui nous arrange (et on donne une esquisse de sa démonstration telle qu'elle est faite dans \cite{Wed}) :

\thr[Wedhorn]
\label{thrwed}
Soit $F/\QQ_p$ une extension finie non ramifiée, $k$ un corps parfait, et $G/\Spec(k)$ un groupe 
$p$-divisible, muni d'une action de $\mathcal O$ (de hauteur et 
signature fixée). Alors il existe une suite $G = G_0, G_1,\dots,G_n$ 
de groupes $p$-divisibles sur des corps parfaits $k_i$, munis d'une action de $\mathcal O$, 
tels que pour tout $i$, $G$ soit une spécialisation de $G_i$, et 
tel que $G_n/\Spec(k_n)$ soit $\mu$-ordinaire.
\ethr

\dem[Esquisse de la preuve de \cite{Wed} section 4.2]
On peut déformer le display associé à $G_i$ en un display sur $k_i[[T]]$, 
tel que son changement de base à $k_i((T))^{perf} =: k_{i+1}$ soit associé à un groupe 
$p$-divisible avec action de $\mathcal O$ vérifiant que son polygone de Newton est 
strictement plus petit que celui de $G_i$ (si $G_i$ non $\mu$-ordinaire). 
Donc par récurrence, on trouve un groupe $p$-divisible sur $k_n$ qui est $\mu$-ordinaire.
Pour conclure, si $U$ est un ouvert de $|\mathcal{BT}_r|$ qui contient $G_i$, 
et si $\widetilde{G_i}$ est la déformation de $G_i$ sur $k_i[[T]]$ dans $\mathcal{BT_{\infty}}$, 
alors il suffit de voir que $\widetilde{G_i} \dans U(k_i[[T]])$, et on aura donc que 
$G_{i+1} = (\widetilde{G_i})_\eta \otimes k_{i+1} \dans U$.
\edem

\thr
Le lieu $\mu$-ordinaire dans $\mathcal{BT}_r^{\mu-ord} \subset \mathcal{BT}_r$ est 
(topologiquement) dense. En particulier, $\widetilde{^\mu\Ha}$ et $\widetilde{\Ha_b}$, pour tout 
$b \dans \{q_\tau : \tau \dans \mathcal I\}$, sont des diviseurs de Cartier sur $\mathcal{BT}_r$.
\ethr

\dem
Soit $k$ un corps parfait, $x \dans \mathcal{BT}_r(k)$ un $\mathcal O$-module $p$-divisible 
tronqué. D'après \cite{Wed2}, Proposition 3.2, il existe $\tilde x \dans \mathcal{BT}_\infty(k)$ au 
dessus de $x$ (i.e. tel que $\tilde x[p^r] = x$). D'après le théorème \ref{thrwed} précédent, il existe 
$\tilde y \dans  \mathcal{BT}_\infty(k)$ une générisation de $\tilde x$ qui est $\mu$-ordinaire. 
Maintenant $y = \tilde y[p^r] \dans  \mathcal{BT}_r$ est une générisation de $x$, et comme être 
$\mu$-ordinaire ne dépends que de la $p$-torsion (\ref{thrmoo}), $y \dans  \mathcal{BT}_r^{\mu-ord}(k)$.
Comme le champ $\mathcal{BT}_r$ est algébrique, on en déduit que le lieu $\mu$-ordinaire de $|\mathcal{BT}_r|$ est (topologiquement) dense.

Donc par les résultats de l'appendice \ref{AppB}, comme $\mathcal{BT}_r$ est lisse, donc réduit, $^\mu\Ha$ est un diviseur de 
Cartier sur $\mathcal{BT}_r$. Comme c'est un produit, on en déduit donc le même résultat pour $\Ha_b$.
\edem

\section{Multiplicité du diviseur $^\mu\Ha$}
\label{sect10}

Rappelons que l'on est dans la situation de $F/\QQ_p$ une extension non ramifiée, et d'un $\mathcal O = \mathcal O_F$-module $p$-divisible $G$ (éventuellement tronqué), 
de signature $(p_\tau,q_\tau)_\tau$, sur une base $S$ 
(par exemple le champ $\mathcal{BT}_{r,\mathcal O,h,(p_\tau)}$ précédemment défini). 
Afin de démontrer des propriétés sur les invariants de Hasse, notamment qu'ils sont réduits, il serait utile de ramener de tels 
énoncés à des propriétés géométriques sur certaines variétés de Shimura, dont on connait bien la géométrie (cf. \cite{WedVieh} ou \cite{Ham}).

\subsection{Globalisation : Variétés de Shimura}
\label{sect9}
On va définir une donnée PEL simple (cf. \cite{Kot2}) $\mathcal D = (B,\star,U_\QQ,<,>)$, de telle sorte que le groupe $p$-divisible associé à la variété abélienne universelle
sur $Sh_\mathcal D$, la variété de Shimura (entière) associée à $\mathcal D$, soit relié à $\mathcal{BT}_{r,\mathcal O,h,(p_\tau,q_\tau)}$.

À priori les $\mathcal O$-modules $p$-divisibles qui nous intéressent ne sont pas polarisés (mais s'ils le sont on pourra leurs appliquer nos résultats en "oubliant" la polarisation)
et sont donc reliés à des données de Shimura de type unitaire, qui se scindent au dessus du nombre premier considéré (le cas (AL) de \cite{WedVieh}). Commençons donc par définir 
un corps de nombre CM.
Soit $L^+/\QQ$ un corps de nombres totalement réel tel que $p$ soit inerte dans $L^+$ et tel que $O_{L^+,p} = \mathcal O$. Un tel corps de nombres existe, par \cite{Farthese}, 
10.1.2 par exemple. Soit $K/\QQ$ un corps quadratique imaginaire dans lequel $p$ est décomposé, et posons $L = L^+K$, c'est un corps CM sur $\QQ$ 
(de corps totalement réel $L^+$). On peut alors associer à notre donnée locale $(\mathcal O,h,(p_\tau,q_\tau))$ une donnée globale entière (cf. \cite{Farthese} 10.1.1.3, et proposition 10.1.3),
\[ \mathcal D = (B,\star,V,<,>,\mathcal O_{B,p},\Lambda,h_0),\]
où $B/\QQ$ est un algèbre à division de centre $L$, $\star$ est une involution de seconde espèce de $B$ qui induit l'automorphisme non trivial de $\Gal(L/L^+)$, $V = B$, $<,>$  donné par $<x,y> = \tr_{B/\QQ}(xy^\star)$, $\mathcal O_{B,p} = \Lambda$ est un ordre maximal en $p$ de $B$.
On note $G/\QQ$ le groupe des similitudes unitaires associé à $(B,\star,V,<,>)$ et $h_0 : \CC \fleche \End_B(V)$.
Ces choix assurent que si $G_1$ est le groupe unitaire associé à $G$, alors $G_1(\QQ_p) = \GL_h(F)$ et $G_1(\RR) = \prod_{\tau : L^+ \fleche \RR} U(p_\tau,q_\tau)$. On note $E'$ le corps réflexe global associé.

\defi
\label{defiShimura}
On dira que $\mathcal D$ est une donnée de Shimura qui globalise la donnée locale $(\mathcal O,h,(p_\tau,q_\tau))$. Cette donnée $\mathcal D$ n'est pas unique.
\edefi

\rem
On a $h^2 = [B:L]$ et un isomorphisme $\mathcal O_B \otimes \ZZ_p = M_h(\mathcal O_L \otimes_{\ZZ} \ZZ_p)$.
\erem

À cette donnée $\mathcal D$, et un plongement $\nu : \overline{\QQ} \fleche \overline\QQ_p$, Kottwitz associe une tour de Variétés de Shimura $(\mathcal S_{K^p})_{K^p}$ sur 
$\mathcal O_{E'_\nu}$ ($E'_\nu$ est la complétion de $E'$ via $\nu$) pour des sous-groupes compacts ouverts $K^p \subset G(\mathbb A_f^p)$ suffisamment petit. On fixe un tel $K^p$, et on note $\mathcal S = \mathcal S_{K^p}$, puisque le choix d'un niveau hors $p$ n'est pas central dans nos constructions.

$\mathcal S$ paramètre des quadruplets $(A,\lambda,i,\eta^p)$ où $A$ est un schéma abélien de genre $fh^2$ (où $f = [\mathcal O:\ZZ_p]$), $\lambda$ une polarisation, 
\[ i : O_B \fleche \End(A),\]
compatible à l'involution de Rosatti associée à $\lambda$ et $\eta^p$ est une structure de niveau hors $p$.

Considérons $\mathcal A \fleche \mathcal S$ le schéma abélien universel, et $\mathcal A[p^\infty]$ son groupe $p$-divisible. On note $L_p = L \otimes_{\QQ} \QQ_p$. C'est un $O_B \otimes \ZZ_p = M_h(O_{L_p})$-module, que l'on peut donc réécrire, grâce à l'équivalence de Morita,
\[ \mathcal A[p^\infty] = \mathcal O_{L_p}^h \otimes H,\]
où $H$ est un $\mathcal O_{L_p}$-module $p$-divisible de hauteur $2fh$. Or $p$ est décomposé dans $L$, donc $\mathcal O_{L_p} = \mathcal O \times\mathcal O$ et on peut donc écrire grâce à l'involution de Rosatti,
\[ H = G \times G^D,\]
où $G$ est un $\mathcal O$-module (non polarisé), de $\mathcal O$-hauteur $h$ et de signature $(p_\tau,q_\tau)_\tau$. La polarisation de $H$ échange $G$ et $G^D$.

\pro
\label{prof}
À $G[p^r]$ est alors associé une flèche,
\[ f : \mathcal S \fleche \mathcal{BT}_{r,\mathcal O,h,(p_\tau,q_\tau)},\]
et on définit l'invariant de Hasse $\mu$-ordinaire de $\mathcal S$ comme $\widetilde{^\mu\Ha} := \widetilde{^\mu\Ha(G)} = f^*(\widetilde{^\mu\Ha(G^{univ})})$.
\epro

Le lieu $\mu$-ordinaire de $\mathcal S$ correspond au lieu d'inversibilité de $^\mu\Ha$, il est dense (\cite{Wed}) et il est affine dans la compactification minimale.

\rem
Etant donné une donnée de Shimura PEL non ramifiée, plus générale que la précédente, on pourrait encore découper $\mathcal A[p^\infty]$ de manière semblable et associer un (des) invariants de Hasse $\mu$-ordinaires (partiels) à chaque composant.
\erem

\subsection{Espaces de déformations}
\label{ssect83}
Soit $E$ le corps réflex (local) défini dans la section \ref{sect8} et $\kappa_E$ son corps résiduel, de caractéristique $p$.
Soit $S$ un $\mathcal O_E$-schéma.  
Un groupe $p$-divisible sur $S$ avec $\mathcal O$-structure est la donnée $\underline{H} = (H,\iota)$ d'un $S$-groupe $p$-divisible et d'un morphisme,
\[ \iota : \mathcal O \fleche \End(H).\]
Soit $k$  un corps parfait de caractéristique $p$. Soit $\mathcal C_k$ la catégorie des $\mathcal W(k)$-algèbres artiniennes locales $R$ de corps résiduel $k$ (avec les morphismes évidents). 
 
\thr
\label{thr121}
Soit $\underline{\mathbb H} = (\mathbb H, \overline{\iota})$ un $\mathcal O$-module $p$-divisible sur $k$. Notons 
\[\Def_{\underline{\mathbb H}} : \mathcal C_k \fleche \mathfrak{Ens},\] 
le foncteur qui à $R$ associe les classes d'isomorphismes de groupes $p$-divisibles avec 
$\mathcal O$-structure $(H,\iota)$ sur $\Spec(R)$ munis d'une trivialisation,
\[\rho_H : H\times\Spec(k) \overset{\simeq}{\fleche} \mathbb H,\]
qui commute avec $\iota$ et $\overline{\iota}$.
Alors $\Def_{\underline{\mathbb H}}$
est pro-représentable par une $W(k)$-algèbre noetherienne complète $R_{\underline{\mathbb H}}$ de corps résiduel $k$. De plus ce foncteur est formellement lisse, donc $\Def_{\underline{\mathbb H}}$ est (pro-)représentable par 
\[W(k)[[t_1,...,t_r]], \quad \text{où } r = \sum_{\tau \in \mathcal I} p_\tau q_\tau.\]
\ethr

\dem
On sait déjà par \cite{RZ}, par exemple, qu'en oubliant la $\mathcal O$-structure, 
$\Def_{\mathbb H}$ est représentable. Notons $\widetilde{H}$ le groupe $p$-divisible universel 
sur $R_{\mathbb H}$ qui représente ce dernier foncteur. Maintenant $\mathcal O$ agit par 
isogénies sur $\mathbb H$, et par rigidité des quasi-isogénies, par quasi-isogénies sur 
$\widetilde{H}$. Le lieu où $\mathcal O$ agit par isogénies, et donc $\Def_{\underline{\mathbb H}}$, apparaît comme un sous-schéma schéma fermé de $\Spec(R_{\mathbb H})$  \cite{RZ} Proposition 2.9.
La suite découle de \cite{Wed2}, corollaire 2.10. Le calcul de la dimension est donné par la théorie de Grothendieck-Messing.
\edem

\rem
Dans les cas (AU) et (C) de \cite{WedVieh} on a aussi des espaces de déformations semblables, notés $\Def_{\mathbb H}^{AU}$ et $\Def_{\mathbb H}^{C}$. Leurs dimensions respectives sont
\[ \frac{1}{2} \sum_{\tau} p_\tau q_\tau \quad \text{et} \quad f\frac{h}{2}(\frac{h}{2}+1)/2.\]
\erem

Introduisons $\mathcal S = \mathcal S_{K^p}/\Spec(O_{E',\nu})$ la variété de Shimura (de niveau hyperspécial en $p$) associée à la donnée PEL entière $\mathcal D$ 
de type (AL) (cf section \ref{sect9}; définition \ref{defiShimura}) qui globalise la donnée $(\mathcal O,h,(p_\tau)_{\tau \in \mathcal I})$ au sens de la définition \ref{defiShimura}.

\pro
Soit $\overline{\mathcal S}$ est la réduction modulo $p$ de la variété de Shimura $\mathcal S$ et $x \dans \overline{\mathcal S}(k)$, où $k$ est un corps parfait de caractéristique $p$, et $\underline{G_x}$ le $\mathcal O$-module $p$-divisible associé à $x$ par la recette de la sous-section \ref{sect9}; alors,
\[ \Spf(\widehat{\mathcal O}_{\mathcal{S},x}) \simeq \Def_{\underline{G_x}}.\]
\epro

\dem
L'assertion découle de la théorie de Serre-Tate, cf. \cite{KatzST}, Theorem 1.2.1.
\edem

\rem
On peut bien sur écrire de telles identifications des anneaux locaux des variétés de Shimura générales en termes de produits des espaces du type 
$\Def_{\mathbb H}$ et $\Def_{\mathbb H}^{AU}$ dans le cas $(A)$ et $\Def_{\mathbb H}^C$ dans le cas (C). Comme on a supposé qu'il n'y a qu'un idéal au dessus de $p$ dans la 
définition \ref{defiShimura}, on est ici dans un cas particulier où il n'y a qu'un espace de déformations.
\erem

Soit $k$ un corps algébriquement clos de caractéristique $p$.
Soit $\underline{\mathbb H} = (\mathbb H, \overline{\iota})$ un groupe $p$-divisible sur $k$ avec $\mathcal O$-structure, et soit,
\[\Def_{\underline{\mathbb H}} = \Spf(W(k)[[t_1,\dots,t_r]]),\]
l'espace des déformations (par isomorphisme) de $\underline{\mathbb H} = (\mathbb H, \overline{\iota})$.
Soit $H^{univ} \fleche \Spf(W(k)[[t_1,\dots,t_r]])$, la déformation universelle. 

À $H^{univ} \fleche \Spf(k[[t_1,\dots,t_r]])$, sa réduction modulo $p$, est associé un display $\mathcal P = (P,Q,F,V^{-1})$, (cf. \cite{Zink2}, et \cite{Zink} si $\mathbb H$ est formel), où $P$ est un module libre 
sur $\hat W(k[[t_1,\dots,t_r]])$ (même sur $W(k[[t_1,\dots,t_r]])$ si $\mathbb H$ est formel), 
$Q \subset P$ est un sous-module,
$F : P \fleche P$, et $V^{-1} : Q \fleche P$ sont $^F$-linéaires où $^F$ est le Frobenius de $\hat W(k[[t_1,\dots,t_r]])$, vérifiant certaines conditions.

Or à $H^{univ}$ est aussi associé un cristal par Berthelot-Breen-Messing \cite{BBM}, $\mathbb D(H^{univ})$, que l'on peut voir, quitte à l'évaluer sur l'épaississement à puissances 
divisées (cf. \cite{Zink} section 2.3, \cite{Gr} Chapitre IV 3.1, \cite{Lau} Lemma 1.16),
$\hat W(k[[t_1,\dots,t_r]]) \twoheadrightarrow k[[t_1,\dots,t_r]]$, comme un module localement libre sur $\hat W(k[[t_1,\dots,t_r]])$, et on a la relation suivante, \cite{Zink}, Theorem 6, \cite{Lau} Theorem B,
\[ \mathbb D(H^{univ})_{\hat W(k[[t_i]]) \rightarrow k[[t_i]]} = P.\]
Dans \cite{Zink} formule (87), Zink calcule explicitement le display de la déformation universelle d'un groupe $p$-divisible, et on aimerait utiliser celui-ci pour calculer explicitement les invariants
de Hasse partiels.

Tout d'abord introduisons (rapidement) les displays munis d'une $\mathcal O$-action. On suppose que tous les schémas sont au-dessus de $\mathcal O$ dans la suite.

Comme $H^{univ}$ est muni d'une $\mathcal O$-action, on peut décomposer,
\[ P = \bigoplus_\tau P_\tau,\]
idem pour $Q$, et les morphismes $F$ et $V^{-1}$ sont $\sigma$-gradués pour cette décomposition.
Le calcul du display de la déformation universelle d'un groupe $p$-divisible formel (cf \cite{Zink} formule 87) s'adapte ici à la fois aux groupes $p$-divisibles non necessairement formels (voir aussi \cite{Lau} section 3.2) et à l'action de $\mathcal O$, en ne considérant que des objets compatibles 
à la décomposition précédente de $\mathcal P$,
l'outil principal étant la théorie des déformations de Grothendieck-Messing.
Ceci s'adapte dans le cas d'une $\mathcal O$-action, à condition que la déformation de l'algèbre de Lie soit $\mathcal O$-graduée.
On en déduit alors la proposition suivante.

\pro
Soit $k$ algébriquement clos de caractéristique $p$.
Soit $\underline{\mathbb H} = (\mathbb H, \overline{\iota})$ un groupe $p$-divisible sur $k$ avec $\mathcal O$-action et notons
$\mathcal P_0 = (P_0,Q_0,F_0,V_0^{-1}) = (\mathbb D(\mathbb H), V(\mathbb D(\mathbb H)),F,\frac{1}{p}F)$ son display (de Dieudonné) sur $k$, où $(\mathbb D(\mathbb H),F,V)$ 
est le module de Dieudonné de $\mathbb H$.
Choisissons une base de $P_0$, $(e^\tau_i)_{\tau \in \mathcal I, i \in \{1,\dots h\}}$, adaptée à la décomposition, 
\[ P_0 = \bigoplus_{\tau \in \mathcal I} P_0^\tau,\]
de telle manière qu'après réduction modulo $p$, pour tout $\tau$, $(e_1^\tau, \dots ,e_{p_\tau}^\tau)$ induise une base de $P_0^\tau/VP_0^{\sigma\tau} = \omega_{H,\tau}^\vee$.
Notons $B$ la matrice (à coefficient dans $W(k)$) dans cette base de \[F \oplus V^{-1}: P_0 = T \oplus L \fleche P_0,\]
où \[L = \bigoplus L_\tau \quad \text{et} \quad T = \bigoplus T_\tau,\] et $T_\tau = W(k)<e^\tau_1,\dots,e^\tau_{p_\tau}>$, $L_\tau = W(k)<e^\tau_{p_\tau + 1},\dots,e^\tau_h>$.
Le display de la déformation universelle $H^{univ} \fleche \Spf(k[[t_1,\dots,t_r]])$, est alors libre et donné par la matrice de Hasse-Witt (i.e. de $F \oplus V^{-1}$) suivante,
\[ \Diag\left(
\left(
\begin{array}{cc}
\Id_{p_\tau}  &  ([t_{k,l}^\tau])_{(k,l) \in \{1,\dots,p_\tau\}\times \{1,\dots,q_\tau\}} \\
  0 &   \Id_{q_\tau}
\end{array}
\right), \tau \in \mathcal I\right) B,
\]
où $[.] : k[[t_i, i]] \fleche W(k[[t_i, i]])$ est l'application de Teichmuller.
\epro 

\rem
On note pour un ensemble fini $\{1,\dots,n\}$, et $(A_i)_{1\leq i\leq n}$ une collection de matrices, la matrice diagonale par blocs,
\[ \Diag(A_i : i \dans \{1,\dots,n\}) = 
\left(
\begin{array}{cccc}
 A_1 &   &   \\
  &  A_2 &   \\
  &   &   \ddots & \\
    &   &   & A_n \\
\end{array}
\right).
\]
À priori l'ensemble $\mathcal I$ des plongements n'est pas ordonné, mais quitte à choisir un plongement $\tau_0 \dans \mathcal I$, on a une bijection,
\[ 
\begin{array}{ccc}
 \{1,\dots,f\}  & \fleche  &   \mathcal I\\
 i &\longmapsto   &\sigma^i\tau_0   
\end{array},
\]
où $\sigma$ est le Frobenius de $\mathcal O$. Néanmoins on n'insiste pas sur cette remarque, et on espère que la notation est suffisamment claire pour ne pas induire le lecteur en erreur.
\erem

\exe
Soit le groupe $p$-divisible de hauteur 4 sur $k$ avec $\ZZ_{p^2}$-action et signature $((1,1),(0,2))$ donné par la matrice de Verschiebung,
\[
V = \left(
\begin{array}{cc}
0  &  
\left(
\begin{array}{cc}
 1 &      \\
  & 1     
\end{array}
\right)
  \\
\left(
\begin{array}{cc}
    & p  \\
  1  &   
\end{array}
\right)
  &     0
\end{array}
\right)
\]
Il n'est pas $\mu$-ordinaire. 
Sa matrice de Hasse Witt ($F \oplus \frac{1}{p}F$) est donnée dans la même base par,
\[
\left(
\begin{array}{cc}
0  &  
\left(
\begin{array}{cc}
 0 &   1   \\
 1 & 0     
\end{array}
\right)
  \\
\left(
\begin{array}{cc}
  1 &  \\
    &  1 
\end{array}
\right)
  &     0
\end{array}
\right).
\]
Le display de la déformation universelle de $G$ à coefficients dans $W(k[[X]])$ est donc donné "dans la même base" par la matrice de Hasse-Witt,
\[
\left(
\begin{array}{cc}

\left(
\begin{array}{cc}
 1 &     \\
 X & 1     
\end{array}
\right) &0    
  \\
0 & \left(
\begin{array}{cc}
  1 &  \\
    &  1 
\end{array}
\right)
\end{array}
\right) \cdot
\left(
\begin{array}{cc}
0  &  
\left(
\begin{array}{cc}
 0 &   1   \\
 1 & 0     
\end{array}
\right)
  \\
\left(
\begin{array}{cc}
  1 &  \\
    &  1 
\end{array}
\right)
  &     
\end{array}
\right) =
\left(
\begin{array}{cc}
0  &  
\left(
\begin{array}{cc}
 0 &   1   \\
 1 & X     
\end{array}
\right)
  \\
\left(
\begin{array}{cc}
  1 &  \\
    &  1 
\end{array}
\right)
  &     
\end{array}
\right).
\]

\eexe

\lem
Soit $R$ un anneau de caractéristique $p$. Le Frobenius $^F$ de $W(R)$ relève le Frobenius de $R$ mais aussi celui de $W(R)/pW(R)$. Autrement dit, c'est un \textit{relèvement fort}
de Frobenius. C'est encore vrai pour $\hat W(R) \subset W(R)$.
\elem

\dem
Soit $x=(x_0,x_1,\dots) \dans W(R)$, alors écrivons que $x = (x_0,0,\dots) + ^Vy$.
On a donc que $^Fx \dans ^F(x_0,\dots) + pW(R)$ puisque $^{FV}=p$. Mais on a aussi que $^VW(R) = I_R$ et
\[ I_R^p \subset pI_R,\]
puisque $I_R$ est muni de puissances divisées (\cite{Gr} Chap IV, 3.1),
donc $x^p \dans (x_0,\dots)^p + pW(R)$ et comme $(x_0,0,\dots)^p = ^F(x_0,0,\dots)$, cela conclut pour $W(R)$. Pour $\hat W(R)$, on utilise que $pW(R) \cap \hat W(R) = p\hat 
W(R)$ (\cite{Lau} Lemma 1.15 si $p \geq 3$ ou si $p = 2$ et $2R = 0$, la preuve fonctionne encore : si $y =(y_0,y_1,\dots ) \dans W(\mathfrak m_R)
$ tel que $2y \dans \hat W(R)$, alors $2y = (0,y_0^2,y_1^2,\dots)$ or si $y_i \dans \mathfrak m_R$ tel que $y_i^2 \dans \mathfrak m_R^{2n}$, alors $y_i \dans \mathfrak m_R^n$).
\edem

\pro
\label{prodefouniv}
Gardons les notations de la proposition précédente, et notons $\mathbb E =\mathbb D(H^{univ})_{\hat W(k[[t_1,\dots,t_r]]) \rightarrow k[[t_1,\dots,t_r]]}$.
Notons $\phi = ^F$ le morphisme de $\hat W(k[[t_i ; i \dans \{1,\dots, r\}]])$.
Alors il existe pour tout $\tau$ un unique morphisme,
\[ \widetilde{\Ha}_\tau : \bigwedge^{q_\tau} \mathbb E_\tau \fleche \Fil(\bigwedge^{q_\tau} \mathbb E_\tau^{(\phi)}),\]
tel que $p^{k_\tau}\widetilde{\Ha}_\tau = V^f$.
Ce morphisme correspond à l'invariant de Hasse partiel construit comme dans la section \ref{sect8}.
De plus, d'après la proposition précédente et \cite{Zink} Theorem 6, on peut calculer explicitement la matrice de
\[ V : \mathbb E \fleche \mathbb E^{(\phi)},\]
sur $\mathbb E$, elle est donnée par :
\[A\Diag\left(
\left(
\begin{array}{cc}
\Id_{p_\tau}  &  -([t_{k,l}^\tau])_{(k,l) \in \{1,\dots,p_\tau\}\times \{1,\dots,q_\tau\}} \\
  0 &   \Id_{q_\tau}
\end{array}
\right) ; \tau \in \mathcal I\right),\]
où $A \dans M_{fh}(W(k))$ est un relevement de la matrice de $V$ de $\mathbb H$, dans une base adaptée à l'action de $\mathcal O$ et à la filtration de Hodge (comme dans la proposition précédente).
\epro

\dem
$\phi$ est un relèvement fort de Frobenius sur $\hat W(k)[[t_i ; i ]]$, donc par le corollaire \ref{cor33}, $V^f$ est bien divisible par 
$p^{k_\tau}$, et comme $\hat W(k[[t_i ; i]])$ est sans $p$-torsion (on plonge $k[[\underline t]]$ dans son perfectisé, qui induit une injection $W(k[[\underline t]] \subset W(k[[\underline t]]^{perf})$ qui est sans torsion, et $\hat W(R) \subset W(R)$), une telle division est unique, et donc correspond à l'invariant de Hasse construit dans 
la section \ref{sect8}.
$\phi$ n'est pas inversible sur $\hat W(k[[t_i ; i]])$, mais on peut tout de même construire une application $V : \mathbb E \fleche \mathbb E^{(\phi)}$ (notée $V^\#$ dans \cite{Zink}), alors la proposition précédente et le théorème B de \cite{Lau} (ou le théorème 6 de \cite{Zink}) nous assure la forme voulue pour la matrice de $V$.
\edem

\rem
À l'aide de la proposition précédente, on va pouvoir explicitement calculer $\widetilde\Ha_\tau$ sur les espaces de déformations, et donc localement sur les variétés 
de Shimura.

Dans certains calculs (notamment pour la section sur la dualité, mais pas pour tous les calculs de multiplicités), on aurait pu se passer de la théorie de \cite{Zink2},\cite{Lau} et utiliser 
seulement les displays des groupes formels en découpant le groupe $p$-divisible $\mathbb H$ en parties multiplicative, biconnexe, étale : 
$\mathbb H = \mathbb H^m \oplus \mathbb H^{00} \oplus \mathbb H^{et}$ et en 
utilisant la proposition \ref{proproduit} et la remarque \ref{rem62}.

Ce découpage fonctionne sur les anneaux de déformation lorsque $\mathbb H$ est $\mu$-ordinaire (ou qu'il a déjà une partie étale et multiplicative de taille maximale par rapport à 
sa signature) puisque le polygone de Newton diminue par déformation. Malheureusement pour les calculs de multiplicités, si la signature permet une partie étale (i.e. 
$\min_\tau q_\tau \neq 0$), sur la strate de Newton associée au plus petit des $q_\tau$, le groupe $\mathbb H$ que l'on considère est presque formel (sa partie étale n'est pas de 
taille maximale) mais le groupe en fibre générique sur l'espace de déformation a une partie étale maximale (il est $\mu$-ordinaire), autrement dit, même en découpant la partie étale 
de $\mathbb H$, $H^{univ}$ n'est pas formel. Ceci dit, une fois la dualité connue, on aurait tout de même pu s'en sortir quitte à passer au dual de $H^{univ}$ (sauf lorsqu'il y a un lieu 
ordinaire, i.e. $p_\tau = p_{\tau'} (\not\in\{0,h\})$ pour tout $\tau,\tau'$).
\erem

\subsection{Strates de Newton et multiplicité des invariants.}
\hyp
\label{hyp1}
Supposons maintenant que tous les $q_\tau, \tau \dans \mathcal I$, tels que $q_\tau \not\in \{0,h\}$, sont distincts, i.e. 
\[q_\tau \not\in \{0,h\} \text{ et } \tau' \neq \tau \Rightarrow q_\tau \neq q_{\tau'}.\]
\ehyp

Alors on a le théorème suivant :

\thr
\label{thrreduit}
Sous l'hypothèse \ref{hyp1}, le diviseur de Cartier $V(^\mu\Ha) \subset \mathcal{BT}_{r,h,(p_\tau)}$, associé à l'invariant de Hasse 
$\mu$-ordinaire, est réduit.
\ethr

\pro
Soit $S/\Spec(\mathcal O_E)$ la donnée de Shimura PEL associée au problème précédent (cf. \ref{defiShimura}), $\overline{\mathcal S}$ sa réduction modulo $p$. Alors le lieu $\mu$-ordinaire est dense, et sur un ouvert dense du complémentaire, l'invariant de Hasse $^\mu\Ha$ définit une forme linéaire non nulle dans l'espace tangent de $\overline{\mathcal S}$. En particulier $V(^\mu\Ha) \subset \overline{\mathcal S}$ est réduit.
\epro

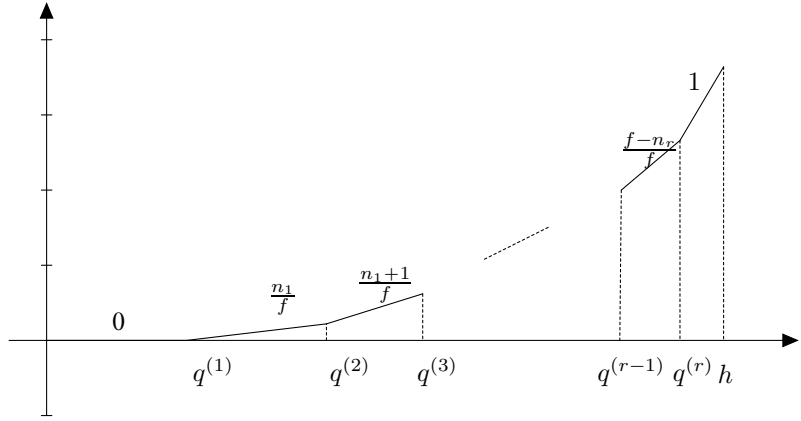
\begin{figure}[h]
\caption{$\mathcal O$-polygone de Hodge associée à la signature $(q_\tau)_{\tau\in \mathcal I}$.}
\label{fig4}
\begin{center}
\begin{tikzpicture}[line cap=round,line join=round,>=triangle 45,x=1.0cm,y=1.0cm]
\draw[->,color=black] (-0.5,0) -- (10,0);
\foreach \x in {,1,2,3,4,5,6,7,8,9}
\draw[shift={(\x,0)},color=black] (0pt,-2pt);
\draw[->,color=black] (0,-1) -- (0,4.5);
\foreach \y in {-1,1,2,3,4}
\draw[shift={(0,\y)},color=black] (2pt,0pt) -- (-2pt,0pt);
\clip(-0.5,-1) rectangle (10,4.5);
\draw (0,0)-- (1.86,0);
\draw (1.86,0)-- (3.72,0.22);
\draw (3.72,0.22)-- (5,0.62);
\draw [dash pattern=on 1pt off 1pt] (5.82,1.08)-- (6.7,1.52);
\draw (7.64,2)-- (8.42,2.66);
\draw (8.42,2.66)-- (9,3.64);
\draw [dash pattern=on 1pt off 1pt] (3.72,0)-- (3.72,0.22);
\draw [dash pattern=on 1pt off 1pt] (5,0.62)-- (5,0);
\draw [dash pattern=on 1pt off 1pt] (7.64,2)-- (7.62,0);
\draw [dash pattern=on 1pt off 1pt] (8.42,2.66)-- (8.42,0);
\draw (1.82,-0.1) node[anchor=north west] {$q^{(1)}$};
\draw (3.64,-0.1) node[anchor=north west] {$q^{(2)}$};
\draw (4.8,-0.1) node[anchor=north west] {$q^{(3)}$};
\draw (7.2,-0.1) node[anchor=north west] {$q^{(r-1)}$};
\draw (8.2,-0.1) node[anchor=north west] {$q^{(r)}$};
\draw [dash pattern=on 1pt off 1pt] (9,3.64)-- (9,0);
\draw (8.8,-0.2) node[anchor=north west] {$h$};
\draw (0.74,0.5) node[anchor=north west] {0};
\draw (2.82,0.88) node[anchor=north west] {$\frac{n_1}{f}$};
\draw (4,1.1) node[anchor=north west] {$\frac{n_1 + 1}{f}$};
\draw (7.5,2.9) node[anchor=north west] {$\frac{f-n_r}{f}$};
\draw (8.40,3.68) node[anchor=north west] {$1$};
\end{tikzpicture}
\end{center}
\end{figure}

\begin{figure}[h]
\caption{Le polygone de Newton $b^l$.}
\label{fig5}
\begin{center}
\definecolor{ffqqqq}{rgb}{1,0,0}
\begin{tikzpicture}[line cap=round,line join=round,>=triangle 45,x=0.5cm,y=0.5cm]
\draw[->,color=black] (0,0) -- (11,0);
\foreach \x in {,2,4,6,8,10}
\draw[shift={(\x,0)},color=black] (0pt,-2pt);
\clip(0,-2) rectangle (11,8);
\draw (1.45,1.68)-- (6.54,2.4);
\draw (6.54,2.4)-- (9.74,6.02);
\draw [color=ffqqqq] (4.82,2.16)-- (8.26,4.34);
\draw [dash pattern=on 1pt off 1pt] (6.54,2.4)-- (6.54,0);
\draw [dash pattern=on 1pt off 1pt] (4.82,2.16)-- (4.81,0);
\draw [dash pattern=on 1pt off 1pt] (8.26,4.34)-- (8.28,0);
\draw (3.7,-0.27) node[anchor=north west] {$q^{(l)}\!-\!1$};
\draw (6.2,-0.3) node[anchor=north west] {$q^{(l)}$};
\draw (7.2,-0.21) node[anchor=north west] {$q^{(l)}\!+\!1$};
\draw [color=ffqqqq] (1.57,5.57)-- (3.1,5.57);
\draw (1.6,4.88)-- (3.07,4.88);
\draw (3.52,5.87) node[anchor=north west] {Newt};
\draw (3.58,5.15) node[anchor=north west] {Hdg};
\end{tikzpicture}
\end{center}
\end{figure}
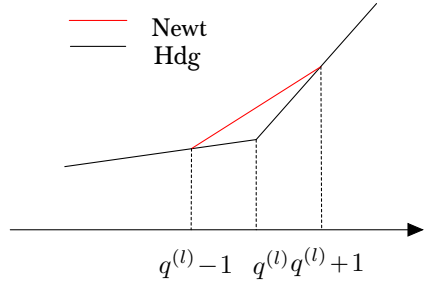

\dem
Le fait que le lieu $\mu$-ordinaire soit dense est déjà connu, cf \cite{Wed}. On va montrer que sur les strates de $V(^\mu\Ha)$ minimales (au sens du polygone de Newton) $^\mu\Ha$ définit une forme linéaire non nulle dans l'espace tangent.
Soit $x \dans \overline{\mathcal S}$. Son $\mathcal O$-polygone de Hodge est celui de la figure \ref{fig4}.
Les $q^{(i)}, i \dans \{1,\dots,f\}$ tels que $q^{(i)} \not\in\{0,h\}$ sont tous distincts par hypothèse. De plus le $\mathcal O$-polygone de Newton de $x$ est au-dessus de son polygone de Hodge, et a même point terminal.
Regardons les points géométriques $x_l$, $l \dans \{1,\dots,f\}$ tels que $q^{(l)} \not\in \{0,h\}$, de $\overline{\mathcal S}$, associé aux groupes $p$-divisible $G/\Spec(k)$, où $k$ algébriquement clos, tels que le polygone de Newton soit égal au polygone de Hodge, sauf autour de $q^{(l)}$ où il est donné par la figure \ref{fig5}.
C'est à dire que l'on s'intéresse aux points $x_l$ dans la strate de Newton de $S$ associée au polygone ci-dessus. On note $\mathcal N_{b^l}$ cette strate. Dans ce cas, par le théorème \ref{thrHN}, on a l'existence de sous groupes $H_1 \subset H_2$ de $G$ tels que l'extension,
\[ 0 \subset H_1 \subset H_2 \subset G,\]
soit scindée, que les polygones de Newton et de Hodge de $H_1$ (respectivement $G/H_2$, 
respectivement $H_2/H_1$) coincident avec la partie des polygones de Newton et de Hodge de 
$G$ définijusqu'à l'abscisse $q^{(l)}-1$ (respectivement, à partir de l'abscisse $q^{(l)}+1$, 
respectivement, entre les abscisses $q^{(l)}-1$ et $q^{(l)}+1$).
Dans ce cas $H_1$ et $G/H_2$ sont $\mu$-ordinaires (pour des signatures $(q_\tau)$ différentes 
de celles d'origine, donnés par la remarque \ref{rem29}). 
Donc par la remarque \ref{rem62} et la proposition \ref{proproduit}, il suffit de regarder 
$H := H_2/H_1$. 

\begin{figure}[h]
\caption{$\mathcal O$-polygones de Hodge et Newton de $H$.}
\label{fig8}
\begin{center}
\definecolor{ffqqqq}{rgb}{1,0,0}
\begin{tikzpicture}[line cap=round,line join=round,>=triangle 45,x=2.0cm,y=2.0cm]
\draw[->,color=black] (-0.5,0) -- (2.3,0);
\foreach \x in {,1,2}
\draw[shift={(\x,0)},color=black] (0pt,-2pt) node[below] {\footnotesize $\x$};
\draw[->,color=black] (0,-0.5) -- (0,2);
\foreach \y in {,0.5,1,1.5}
\draw[shift={(0,\y)},color=black] (2pt,0pt) -- (-2pt,0pt);
\draw[color=black] (0pt,-10pt) node[right] {\footnotesize $0$};
\clip(-0.5,-0.5) rectangle (2.3,2);
\draw (0,0)-- (1,0.37);
\draw (1,0.37)-- (1.95,1.24);
\draw [color=ffqqqq] (0,0)-- (1.95,1.24);
\draw [dash pattern=on 1pt off 1pt] (1,0.37)-- (1,0);
\draw [dash pattern=on 1pt off 1pt] (0,0)-- (-5.9,0);
\draw [dash pattern=on 1pt off 1pt] (1.95,1.24)-- (1.95,0);
\draw [color=ffqqqq] (0.33,1.7)-- (0.67,1.7);
\draw (0.32,1.5)-- (0.69,1.5);
\draw (0.72,1.85) node[anchor=north west] {Newt};
\draw (0.73,1.65) node[anchor=north west] {Hdg};
\draw (0.67,0.32) node[anchor=north west] {$\frac{x}{f}$};
\draw (1.34,0.8) node[anchor=north west] {$\frac{x\!+\!1}{f}$};
\draw (0.7,1.05) node[anchor=north west] {$\frac{2x\!+\!1}{2f}$};
\end{tikzpicture}
\end{center}
\end{figure}
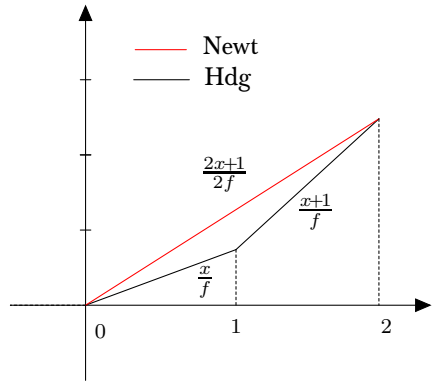
Notons alors,
\[ x := |\{\tau \dans \mathcal I | q_\tau < q^{(l)} \}|, \quad z := |\{ \tau \dans \mathcal I |q_\tau > q^{(l)} \}|.\]
On a par l'hypothèse \ref{hyp1} alors $x + 1 + z = f = [F:\QQ_p]$. Alors $H$ est un groupe $p$-divisible muni d'une action de $\mathcal O$, de $\mathcal O$-hauteur 2 et de signature,
\[ q_\tau (H) = 
\left\{
\begin{array}{ccc}
  0 & \text{ si } q_\tau(G) < q^{(l)}   \\
   1 & \text{ si } q_\tau(G) = q^{(l)}    \\
   2 & \text{ si } q_\tau(G) > q^{(l)}   
\end{array}
\right.
\]
On a donc pour $\tau$ tel que $q_\tau(G) < q^{(l)}$, i.e. $q_\tau(H) = 0$,
$V\mathbb D(H)_{\sigma\tau} = p \mathbb D(H)_\tau,$ 
et pour $\tau$ tel que $q_\tau(G) > q^{(l)}$, i.e. $q_\tau(H) = 2$, 
$V\mathbb D(H)_{\sigma\tau} = \mathbb D(H)_\tau.$
On peut alors fixer une base de $\mathbb D(H)$, i.e. un isomorphisme $\mathbb D(H) = (W(k)^{2})^f$, et un isomorphisme qui conserve l'ordre (normal habituel sur $\{1,\dots,f\}$, donné par le Frobenius sur $\mathcal I$), $\{1,\dots,f\} \simeq \mathcal I$ tel que $\mathbb D(H)_i$ est placé en $i$-position. Avec cet isomorphisme, les matrices de $F$ et $V$ sont données par,
\[
V = \left(
\begin{array}{cccccl}
0 &   &   &   & & p^{n_f}I_2 \\
p^{n_1} I_2& 0  &  & & &\\
 0 &  \ddots & \,\,\ddots & &  & \\
  &   &  A & \,\,\ddots & & \\  
  & 0 &   &  \ddots & \ddots & \\
  &   & &  &  \!\!p^{n_{f-1}}I_2 & 0
\end{array}
\right)
\quad \text{ et,} \]
\[F = \left(
\begin{array}{cccccl}
0 & p^{1-n_1} I_2  &   &   & &  \\
& 0  &  \ddots & & &\\
 &  & \!\!\!\!\ddots & B &  &0 \\
  &   &   & \!\!\ddots &  \ddots &  \\  
  &  &   &  & \ddots & p^{1-n_{f-1}}I_2  \\
  p^{1-n_f}I_2 &   & &  & & 0
\end{array}
\right)
\]
où $n_i \dans \{0,1\}$ et $A,B \dans M_2(W(k))$ sont en position $(k,k+1)$ et $(k+1,k)$ et vérifient $A^\sigma B=B^\sigma A = pI_2$. Comme $q_k = 1$, quitte à changer de base, on a deux choix possibles (d'après la classification des cristaux de hauteur 2 et dimension 1, appliquée au $F^f$-cristal $(W(k)^2,V^f)$, où  $V^f$ est donné par $A$) 
pour $A$ (qui déterminent chacun un unique choix de $B$),
\[ A = 
\left(
\begin{array}{ccc}
  1  &  0 \\
    0 & p  
\end{array}
\right)
\quad \text{ou} \quad A = 
\left(
\begin{array}{ccc}
  0  &  p \\
  1 & 0  
\end{array}
\right).
\]
Le premier cas est le cas où $H$ est $\mu$-ordinaire, ce qui est exclu, on est donc dans le second cas. En particulier, cette strate de Newton est une strate d'Ekedahl-Oort minimale, comme on le vera dans la proposition \ref{pronewteo}. On voit déjà que tout élément de $\mathcal N_{b^l}$ détermine uniquement, grâce à la matrice précédente, sa $p$-torsion.

Retournons alors au calcul de l'invariant de Hasse partiel associé à un point de la strate de Newton-Ekedahl-Oort $\mathcal N_{b^l}$.
D'après la proposition \ref{prodefouniv}, on peut écrire la matrice de $V$ du cristal $\mathcal P$ de la $\mathcal O$-déformation universelle de $H$, dans une base 
($\mathcal P \simeq ((\hat W(k[[X]]))^2)^f$) adaptée aux matrices de $V$ et $F$ précédemment écrites,
\[V^{univ} =  \left(
\begin{array}{ccllcl}
0 &   &   &   & & p^{n_f}I_2 \\
p^{n_1} I_2& 0  &  & & &\\
 0 &  \ddots & \,\,\,\,\ddots & &  & \\
  &   &  
\!\!\!\left(
\begin{array}{ccc}
  0   &p   \\
  1 &   -X
\end{array}
\right)

& \!\!\!\ddots & & \\  
  & 0 &   &  \ddots & \ddots & \\
  &   & &  &  \!\!p^{n_{f-1}}I_2 & 0
\end{array}
\right)
\]
On en déduit alors que,
\[ (V^{univ})^f = p^{\left(\sum_{j \neq k}(n_j)\right)} \Diag\left(\left(
\begin{array}{ccc}
  0   &p   \\
  1 &   -X
\end{array}
\right),\dots,\left(
\begin{array}{ccc}
  0   &p   \\
  1 &   -X
\end{array}
\right)\right).\]
Et donc que $\Ha_k(V^{univ}) = -X$, et donc $^\mu\Ha(G) = -X$ (à un inversible près).
On en déduit donc que l'invariant de Hasse est de multiplicité 1 (i.e. définit une forme linéaire non nulle sur l'espace tangent) sur la réunion,
\[ \bigcup_{l=1}^r \mathcal N_{b^l} \subset V_{S}(^\mu\Ha).\]
Or d'après \cite{Ham}, Theorem 1.1, (b), cette union est dense dans le complémentaire du lieu ordinaire, donc $^\mu\Ha$ est réduit (et c'est aussi le cas de tous les $\Ha_\tau$).
\edem

En fait, on peut dire un petit peu plus sur la géométrie des strates de Newton précédemment considérées :

\pro
\label{pronewteo}
Dans la variété de Shimura simple $\mathcal S$ précédente, la strate de Newton $\mathcal N_{b^l}$ est une strate d'Ekedahl-Oort minimale.
\epro

\dem
On a déjà vu au cours de la démonstration précédente (à l'aide de la filtration Hodge-Newton et du calcul de la signature des gradués) qu'un groupe $p$-divisible dans $\mathcal N_{b^l}$ détermine entièrement (et uniquement) sa $p$-torsion. Pour la réciproque, c'est le lemme suivant.

\lem
\label{lem119}
Soit $k$ un corps algébriquement clos de caractéristique $p$. Soit $N_1,N_2$ deux groupes $p$-divisibles sur $k$, avec action de $\mathcal O$, $\mu$-ordinaires, et $H$ un $\mathcal O$-groupe $p$-divisible sur $k$, de $\mathcal O$-hauteur 2. 
Supposons que toutes les pentes de Newtons vérifient l'inégalité,
\[ \text{pentes de } \Newt_{\mathcal O}(N_1) \leq \text{ pentes de } \Newt_{\mathcal O}(H) \leq \text{ pentes de } \Newt_{\mathcal O}(N_2).\]
Soit $\underline{X}$ un $\mathcal O$-groupe $p$-divisible (sur $k$), dont la $p$-torsion 
(munie de sa $\mathcal O$-structure) vérifie,
\[ \underline{X}[p] \simeq (N_1 \times H \times N_2)[p].\]
Alors \[ \underline{X} \simeq N_1 \times H \times N_2.\]
\elem

\begin{figure}[h]
\caption{$\mathcal O$-polygones de Hodge de $N_1 \times H \times N_2$ et de Newton de $H$}
\label{figNHN}
\centering
\begin{tikzpicture}[line cap=round,line join=round,>=triangle 45,x=1.0cm,y=1.0cm]
\draw[->,color=black] (2.,0.) -- (10.,0.);
\foreach \x in {2.,3.,4.,5.,6.,7.,8.,9.}
\draw[shift={(\x,0)},color=black] (0pt,2pt) -- (0pt,-2pt);
\clip(2.,-1.) rectangle (10.,6.3);
\draw (4.,1.62)-- (6.,2.6);
\draw (6.,2.6)-- (8.,5.28);
\draw [color=ffqqqq] (4.992097411499072,2.1061277316345457)-- (7.009586493060524,3.952845900701102);
\draw (6.,2.6)-- (4.,1.62);
\draw [color=ffqqqq] (7.009586493060524,3.952845900701102)-- (4.992097411499072,2.1061277316345457);
\draw [color=ffqqqq] (3.1540983606557376,5.245207439198856)-- (3.59,5.25);
\draw (3.1737704918032787,4.681258941344779)-- (3.59,4.68);
\draw (3.7442622950819673,5.5167381974248935) node[anchor=north west] {$\Newt_{\mathcal O}(H)$};
\draw (3.7311475409836063,4.931902718168813) node[anchor=north west] {$\Hdg_{\mathcal O}(N_1\times H \times N_2)$};
\draw [dash pattern=on 3pt off 3pt] (6.,2.6)-- (6.,0.);
\draw [dash pattern=on 3pt off 3pt] (7.,0.)-- (7.009586493060522,3.9528459007011008);
\draw [dash pattern=on 3pt off 3pt] (5.,0.)-- (4.992097411499072,2.1061277316345457);
\draw (5.7,-0.07) node[anchor=north west] {$q$};
\draw (4.5,-0.018311874105865) node[anchor=north west] {$q-1$};
\draw (6.5,-0.060085836909870724) node[anchor=north west] {$q+1$};
\end{tikzpicture}
\end{figure}
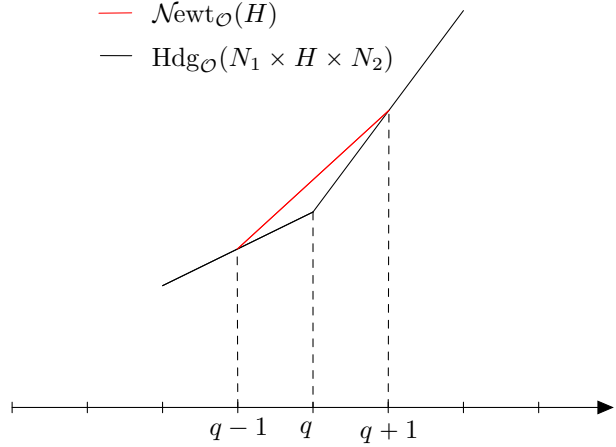

\dem
Par ce qui précede le lemme, il suffit de montrer que $\underline{X}$ est isogène à 
$N_1 \times H \times N_2$, en utilisant que si on a une suite exacte,
\[ 0 \fleche A \fleche B \fleche C \fleche 0\]
de groupes $p$-divisibles sur $k$ avec $\mathcal O$-action, alors $B$ est isogène à 
$A\times C$.
En raisonnant exactement comme dans la démonstration du théorème 1.3.7. de \cite{Moo} (qui découpe petit à petit des parties de $\underline X$ en suivant les pentes de 
$\Hdg_{\mathcal O}(\underline X)$) et en faisant la même chose sur $\underline{X}^D$, on peut se ramener au cas où $H$ n'est pas $\mu$-ordinaire (sinon c'est le théorème 
de Moonen) et le polygone de Newton de $N_1\times H \times N_2$ possède seulement 2 pentes, voir figure \ref{figNHN}. 
Notons que les polygones de Newton de $N_1$ ou de $N_2$ sont égaux à leurs polygones de Hodge (puisqu'ils sont $\mu$-ordinaire) et sont donc isoclins.

Soit $\tau \dans \mathcal I$, $M = \mathbb D(\underline{X})$ le module de Dieudonné de $\underline{X}$ et $V^f : M_\tau \fleche M_\tau$. $V^f$ est divisible par $p^x$, où 
$x = |\{\tau : q_\tau = 0\}|$ et considérons,
\[ \phi_\tau = \frac{1}{p^x}V^f = M_\tau \fleche M_\tau.\]
Notons $W= W(k)$ et $X_0 = N_1 \times H \times N_2$, et $M_0 = M_1 \oplus M_H \oplus M_2$ son module de Dieudonné, alors la matrice de
 $\phi_{\tau,0}$ ($V^f$ est là aussi divisible par $p^x$, puisque $\underline X_0[p] = \underline X[p]$) est,
 \[
\left(
\begin{array}{ccc}
  I_{q-1}&   &   \\
  & 
\left(
\begin{array}{cc}
 0 & p     \\
  1 &  0   
\end{array}
\right)
 &   \\
  &   &   pI_r
\end{array}
\right), \quad \text{où } r = n - q - 1.
\]
Comme $\underline{X}[p] \simeq X_0[p]$, on en déduit que la matrice de $\phi_\tau$ est congrue à celle de $\phi_{\tau,0}$ modulo $p$. Notons,
\[ N = \left(
\begin{array}{cc}
\left(
\begin{array}{cc}

 0 & p     \\
  1 &  0   
\end{array}
\right)
 &   \\
  &    pI_r
\end{array}
\right).
\]
La matrice de $\phi_\tau$ est alors de la forme,
\[
\left(
\begin{array}{cc}
  B_1& B_2   \\
  B_3 &   B_4
\end{array}
\right), \quad \text{où } B_2 \equiv 0, B_3 \equiv 0, B_4 \equiv N \pmod p, \text{ et } B_1 \text{ est inversible.}
\]
Soit $\{f_1,\dots,f_{q-1}\}$ un relèvement d'une base de $N_{1,\tau} = M_{1,\tau}/pM_{1,\tau}$ 
dans $M_{\tau}$, et $\{f_q,\dots f_n\}$ un relèvement dans $M_{\tau}$ d'une base de 
$(M_H \oplus M_2)_\tau \otimes_{W} k$ (compatibles avec les formes précédentes des 
matrices).
On a une bijection,
\[ M_{n-q-1 \times q-1}(W) \overset{\simeq}{\fleche} 
\left\{ 
\begin{array}{ccc}
W-\text{sous-modules } U \subset M_\tau \text{ de rang } q-1, \\
\text{ avec } U/(pM_\tau\cap U) \simeq N_{1,\tau}.
\end{array}
\right\}
\]
qui envoie une matrice $A = (a_{i,j})$ sur l'espace $U_A$ engendré par,
\[ f_1 + p\cdot \sum_{i = q}^n a_{1,j}f_j , \dots , f_{q-1} + p\cdot \sum_{i = q}^n a_{q-1,j}f_j.\]

On peut vérifier par calcul que $\phi_\tau(U_A) = U_{A'}$ où
\[ A' = (\frac{1}{p}B_3 +B_4 \cdot A^{\sigma^f})(B_1 + pB_2\cdot A^{\sigma^f})^{-1}.\]
Donc $U_A$ est stable par $\phi_\tau$ si et seulement si,
\[A(B_1 + pB_2\cdot A^{\sigma^f}) = (\frac{1}{p}B_3 +B_4 \cdot A^{\sigma^f}).\]
On utilise alors la version suivante (légèrement modifiée) du lemme 1.3.6 de Moonen,

\lem

Soit $R$ un anneau local complet de caractéristique résiduelle $p$ (sans $p$-torsion) et $\rho$ un automorphisme de $R$. Soit $C_1,C_2,C_3,C_4$ des matrices de $R$ de taille $r\times s, s\times s, s \times r, r \times r$ respectivement (où $s = q-1, r = n - q +1$), telles que $C_2$ est inversible.
Alors l'équation matricielle,
\begin{equation}
\label{eqmat1} C_1 + XC_2 + (N+pC_3)X^\rho + pXC_4X^\rho = 0,\end{equation}
a une solution $X \dans M_{r\times s}(R)$.
\elem

\dem
Regardons l'équation modulo $p$,
\[ C_1 + XC_2 + NX^\rho \equiv 0 \pmod p.\]
En regardant terme à terme, on a pour tout $(i,j), i \neq 1$,
\[c^1_{i,j} + \sum_k x_{i,k}c^2_{k,l} = 0\]
et pour $i =1$,
\[ c^1_{1,j} + x_{2,j}^\rho + \sum_k x_{1,k}c^2_{k,j} = 0.\]
On peut donc prendre comme solution, $x_{i,j} = (-C_1C_2^{-1})_{i,j}$ pour $i \neq 1$ (comme si $C_3 \equiv 0 \pmod p$) et $x_{1,j} = (-\widetilde{C_1}C_2^{-1})_{1,j}$ où,
\[ \widetilde{C_1} = (\widetilde {c_{i,j}}), \quad \text{où } \widetilde{c_{i,j}} = 
\left\{
\begin{array}{ccc}
c_{i,j}   & \text{si } i \neq 1   \\
c_{1,j} + x_{2,j}^\rho & \text{si } i = 1
\end{array}
\right..
\]
Notons $\Gamma$ cette solution (ou plutôt un relèvement de la solution modulo $p$). Si on pose $X = \Gamma + pX'$, on a que $X$ est solution de (\ref{eqmat1}) si et seulement si $X'$ est solution de \begin{equation} \label{eqmat2} 
C_1' + X'C_2' + (N+pC_3')X^\rho + pX'C_4'X'^\rho = 0,
\end{equation}
où
$C_1' = C_3\Gamma^\rho + \Gamma C_4\Gamma^\rho, C_2' = C_2 + pC_4\Gamma^\rho, C_3' = C_3 + \Gamma C_4,$ et $C_4' = pC_4$.
Donc par approximations successives, comme $R$ est complet, on en déduit une solution à (\ref{eqmat1}).
\edem

\begin{figure}[h]
\caption{Polygones de Hodge et Newton de $H \times N_2$.}
\label{figHN2}
\centering
\begin{tikzpicture}[line cap=round,line join=round,>=triangle 45,x=1.5cm,y=1.0cm]
\draw[->,color=black] (-1.,0.) -- (4.,0.);
\foreach \x in {-1.,1.,2.,3.}
\draw[shift={(\x,0)},color=black] (0pt,2pt) -- (0pt,-2pt);
\draw[->,color=black] (0.,-1.) -- (0.,4.177591320515857);
\foreach \y in {-1.,-0.5,0.5,1.,1.5,2.,2.5,3.,3.5,4.}
\draw[shift={(0,\y)},color=black] (2pt,0pt) -- (-2pt,0pt);
\clip(-1.,-1.) rectangle (4.,4.177591320515857);
\draw (0.,0.)-- (1.0016393442622953,0.40987124463519364);
\draw (1.0016393442622953,0.40987124463519364)-- (2.723770491803279,3.4907010014306157);
\draw [color=ffqqqq] (0.,0.)-- (1.997449458312397,2.1913391867626295);
\draw (1.0016393442622953,0.40987124463519364)-- (0.,0.);
\draw [color=ffqqqq] (1.997449458312397,2.1913391867626295)-- (0.,0.);
\draw [color=ffqqqq] (0.1427863823949168,2.8906447865174494)-- (0.5786880217391788,2.8954373473185937);
\draw (0.17468671700790128,2.4968791259830847)-- (0.5909162252046225,2.4956201846383057);
\draw (0.6516393442622951,3.0813345030103543) node[anchor=north west] {$\Newt_{\mathcal O}(H)$};
\draw (0.6721311475409837,2.66653462611638) node[anchor=north west] {$\Hdg_{\mathcal O}(H \times N_2)$};
\draw [dash pattern=on 1pt off 1pt] (1.0016393442622953,0.40987124463519364)-- (1.0106557377049181,0.);
\draw [dash pattern=on 1pt off 1pt] (2.,0.)-- (1.9974494583123972,2.1913391867626295);
\draw (0.9713114754098361,-0.11855026160030432) node[anchor=north west] {$1$};
\draw (1.8729508196721312,-0.11) node[anchor=north west] {$2$};
\end{tikzpicture}
\end{figure}
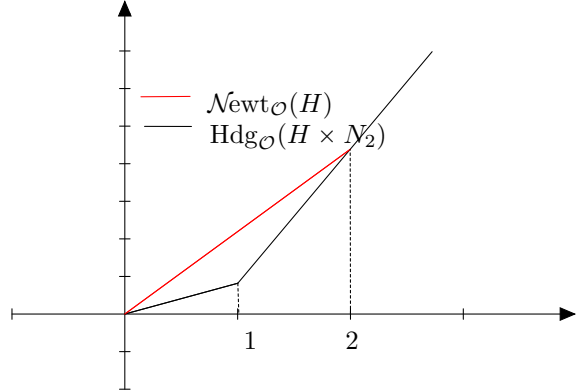

On en déduit donc un sous-$W$-module stable par $\phi_\tau$  $M'_\tau$ de $M_\tau$ et on peut construire comme dans Moonen \cite{Moo} juste avant (1.3.7.4), un sous-$V$-module $M'$ de $M$, avec action de $\mathcal O$, tel que $M' \equiv M_1 \pmod p$. On en déduit au niveau des groupes $p$-divisibles une suite exacte,
\[ 0 \fleche M' \fleche M \fleche Q \fleche 0.\]
Et comme $M' \simeq M_1 \pmod p$, on a que $M'$ est $[p]$-ordinaire, et donc par le théorème 1.3.7 de \cite{Moo}, il est $\mu$-ordinaire, donc $M' \simeq M_1$.
On refait le même raisonnement sur $Q$ en inversant $V$ et $F$ (on pourrait aussi passer à $Q^\vee$, ce qui revient essentiellement au même), qui modulo $p$ est isomorphe à 
$\mathbb D(H \times N_2)$, voir figure \ref{figHN2},
et on en déduit une suite exacte,
\[ 0 \fleche M'' \fleche Q \fleche Q' \fleche 0.\]
où $M'' \simeq M_2$, et $Q'$ est de $\mathcal O$-hauteur 2, donc isomorphe à $\mathbb D(H)$.
On en déduit donc que $\underline{X}$ est isogène à $N_1\times H \times N_2$.
\edem
\edem

\rem
On aurait aussi pu se passer des résultats de \cite{Ham}, et utiliser soit \cite{WedVieh} Théorème 2 (3), ou  \cite{WedVieh} Théorème 2(2), pour montrer que $^\mu\Ha$ est réduit en codimension 1, à l'aide de la proposition précédente, puisque nos strates de Newton $\mathcal N_{b^l}$ sont en fait des strates d'Ekedahl-Oort.
\erem

\dem[Théorème \ref{thrreduit}]
Soit $\overline{\mathcal S}$ une variété de Shimura associée à la donnée $(\mathcal O,h,(p_\tau))$ introduite dans la section \ref{sect9}. 
En considérant la $p^r$-torsion du groupe $p$-divisible "universel" $G[p^r]$ sur 
$\overline{\mathcal S}$, on a une flèche, cf proposition \ref{prof},
\[ f : \overline{\mathcal S} \fleche \mathcal{BT}_r.\]
Considérons \[ \widetilde{S} = \Isom_S(\mathcal O_{G[p^r]},\mathcal O_S^{p^{rh}})\]
qui est un recouvrement étale de $S$, et $X$ la présentation lisse de $\mathcal{BT}_r$ donnée dans la section \ref{sect8}. 
Alors la trivialisation universelle au-dessus de $\widetilde{S}$ de $G[p^r]$ induit une flèche,
\[ \widetilde{f} : \widetilde{S} \fleche X.\]

\lem
Le morphisme $\widetilde{f}$ est lisse et surjectif.
\elem

\dem
En effet, pour la formelle lissité, il suffit de voir que pour tout relèvement infinitésimal d'algèbres locales $R' \fleche R$ à corps résiduels parfaits, et tout diagramme,
\begin{center}
\begin{tikzpicture}[description/.style={fill=white,inner sep=2pt}] 
\matrix (m) [matrix of math nodes, row sep=3em, column sep=2.5em, text height=1.5ex, text depth=0.25ex] at (0,0)
{ 
 \Spec(R) & \widetilde{S} \\
 \Spec(R') & X\\
 };

\path[->,font=\scriptsize] 
(m-1-1) edge node[auto] {$a_0$} (m-1-2)
(m-2-1) edge node[auto] {$x$} (m-2-2)
(m-1-2) edge node[auto] {$\widetilde{f}$} (m-2-2)
(m-1-1) edge node[auto] {} (m-2-1)
;
\end{tikzpicture}
\end{center}
On peut trouver un morphisme,
\[ a : \Spec(R') \fleche \widetilde S\]
qui relève $a_0$ et tel que le diagramme commute,
\begin{center}
\begin{tikzpicture}[description/.style={fill=white,inner sep=2pt}] 
\matrix (m) [matrix of math nodes, row sep=3em, column sep=2.5em, text height=1.5ex, text depth=0.25ex] at (0,0)
{ 
\Spec(R') & &  \\
& \Spec(R) & \widetilde{S} \\
& \Spec(R') & X\\
 };

\path[->,font=\scriptsize] 
(m-2-2) edge node[auto] {$a_0$} (m-2-3)
(m-3-2) edge node[auto] {$x$} (m-3-3)
(m-2-3) edge node[auto] {$\widetilde{f}$} (m-3-3)
(m-2-2) edge node[auto] {} (m-3-2)
(m-2-2) edge node[auto] {} (m-1-1)
;
\draw (m-1-1) edge[double distance=2pt] (m-3-2);
\draw[->] (m-1-1) to[bend left=20] node[right] {$a$} (m-2-3) ;
\end{tikzpicture}
\end{center}
Mais $\widetilde{S} = S \times_{\mathcal{BT}_r} X$. Considérons $\mathcal {BT}_{\infty}$ le champ des groupes $p$-divisibles avec actions de $\mathcal O_F$ et signature fixée. La flèche,
\[ \mathcal {BT}_{\infty} \overset{[p^r]}{\fleche} \mathcal{BT}_r,\]
est formellement lisse d'après \cite{Wed2} 2.8. 
On en déduit donc un morphisme, \[\overline{x}^\infty : \Spec(R') \fleche \mathcal{BT}_\infty.\]
Mais d'après la théorie des déformations de Serre-Tate, la flèche,
\[ S \fleche \mathcal {BT}_\infty\]
est formellement lisse, on en déduit donc un morphisme,
\[ \overline{a} : \Spec(R') \fleche S.\]
On en déduit donc par propriété universelle du produit fibré, la flèche cherchée,
\[a : \Spec(R') \fleche \widetilde{S}. \]
La lissité découle ensuite du fait que $\widetilde{S}$ et $X$ sont tous les deux de présentation finie. Enfin, pour la surjectivité, c'est le théorème 6 de \cite{WedVieh}.
\edem

De plus, en considérant les sous-schémas  fermés $V_{X}(^\mu\Ha)$ et $V_{\widetilde{S}}(^\mu\Ha)$ on a (par définition),
\[  V_{\widetilde{S}}(^\mu\Ha) = \tilde f^*V_{X}(^\mu\Ha).\]
Or on a vu d'après la proposition précédente que  $V_{\widetilde{S}}(^\mu\Ha)$ était un sous-schéma réduit, c'est donc aussi le cas, par lissité et surjectivité 
de $\tilde f$, de $V_{X}(^\mu\Ha)$. Par l'appendice \ref{appB}, le diviseur de Cartier $^\mu\Ha$ sur $\mathcal{BT}_r$ est donc réduit.
\edem

\rem
\label{remred}
Si l'hypothèse \ref{hyp1} n'est pas vérifiée, alors il existe $\tau \neq \tau'$ tels que $q_\tau = q_{\tau'}$, notons $\tau' = \sigma^j\tau$. On peut alors construire 
une première application,
\[ \zeta := \frac{V^j}{p^{l_\tau}} : \det(\omega_{G^D,\tau}) \fleche \det(\omega_{G^D,\sigma^j\tau})^{(p^j)},\]
où $l_\tau = \sum_{i=1}^j \max(q_\tau - q_{\sigma^i\tau},0)$, et une seconde application,
\[ \xi :=\frac{V^{f-j}}{p^{m_\tau}} : \det(\omega_{G^D,\sigma^j\tau}) \fleche \det(\omega_{G^D,\tau})^{(p^{f-j})},\]
où $m_\tau = \sum_{i=j+1}^f \max(q_\tau - q_{\sigma^i\tau},0)$. Alors $l_\tau + m_\tau = k_\tau$ et $\widetilde{\Ha}_\tau = \zeta\otimes\xi^{\otimes(p^j)}$ et donc $\widetilde{\Ha}_\tau$
n'est pas réduit...
\erem

\section{Compatibilité à la dualité}
\label{sect11}

Soit $S$ un schéma de caractéristique $p$, et $G/S$ un groupe $p$-divisible tronqué 
d'échelon $r$, muni d'une action de $\mathcal O$, de hauteur $h$ et signature $(p_\tau,q_\tau)_\tau$ 
(après extension des scalaires). Soit $q \dans \{ q_\tau : \tau \dans \mathcal I\}$ une abscisse de 
rupture de $G$, et supposons que $r > k_\tau$, pour n'importe quel $\tau$ tel que $q = q_\tau$.

\pro
\label{pro131}
Il y a un isomorphisme canonique de $\mathcal O_S$-modules, 
\[\det(\omega_{G^D,q})^{\otimes(p^f-1)} \simeq \det(\omega_{G,q})^{\otimes (p^f-1)}.\]
\epro

La démonstration suit les mêmes lignes que celle de la proposition 2, 2.2.3 de \cite{Far}.

\dem
Notons 
\[\mathcal E = \mathcal{E}xt^1_{cris}(G^D, \mathcal O_{S/\Sigma})_{(S \rightarrow S)},\]
le cristal de Dieudonné covariant, évalué sur $(S\overset{id}{\fleche} S)$, et $\mathcal E_{\tau}$ ses facteurs directs, donnés par l'action de $\mathcal O$. On a la suite exacte, avec les notations de l'introduction,
\[ 0\fleche \omega_{ G^D,\tau} \fleche \mathcal E_{\tau} \fleche (\omega_{G,\tau})^\vee \fleche 0.\]
Celle-ci s'insère dans un morphisme de suite exactes,
 \begin{center}
\begin{tikzpicture}[description/.style={fill=white,inner sep=2pt}] 
\matrix (m) [matrix of math nodes, row sep=3em, column sep=2.5em, text height=1.5ex, text depth=0.25ex] at (0,0)
{ 
0&\omega_{G^D,\tau} & \mathcal E_{\tau} & (\omega_{G,\tau})^\vee& 0\\
0&\omega_{G^D,\sigma^{-1}\tau}^{(p)} & \mathcal E_{\sigma^{-1}\tau}^{(p)}  & (\omega_{G,\sigma^{-1}\tau})^{(p),\vee}& 0\\
 };

\path[->,font=\scriptsize] 
(m-1-1) edge node[auto] {} (m-1-2)
(m-1-2) edge node[auto] {} (m-1-3)
(m-1-3) edge node[auto] {} (m-1-4)
(m-1-4) edge node[auto] {} (m-1-5)
(m-2-1) edge node[auto] {} (m-2-2)
(m-2-2) edge node[auto] {} (m-2-3)
(m-2-3) edge node[auto] {} (m-2-4)
(m-2-4) edge node[auto] {} (m-2-5)

;
\draw[->] (m-1-2) to[bend left=15,left] node[right]  {$V$} (m-2-2);
\draw[->] (m-2-2) to[bend left=15,left] node {$F$} (m-1-2);
\draw[->] (m-1-3) to[bend left=15,left] node[right]  {$V$} (m-2-3) ;
\draw[->] (m-2-3) to[bend left=15,left] node  {$F$} (m-1-3) ;
\draw[->] (m-1-4) to[bend left=15,left] node[right]  {$V$} (m-2-4) ;
\draw[->] (m-2-4) to[bend left=15,left] node {$F$} (m-1-4) ;
\end{tikzpicture}
\end{center}
Comme $F$ restreint à $\omega_{G^D,\sigma^{-1}\tau}^{(p)}$ est nul, on en déduit la factorisation,
\[ F^: (\omega_{G,\sigma^{-1}\tau})^{(p),\vee} \fleche \mathcal E_{\tau}.\]
De même, on peut factoriser,
\[ V : \mathcal E_{\tau} \fleche \omega_{G^D,\sigma^{-1}\tau}^{(p)}.\]
On en déduit un morphisme de complexes (les complexes étant les lignes) qui est exact,
 \begin{center}
\begin{tikzpicture}[description/.style={fill=white,inner sep=2pt}] 
\matrix (m) [matrix of math nodes, row sep=3em, column sep=2.5em, text height=1.5ex, text depth=0.25ex] at (0,0)
{ 
  &\omega_{G^D,\tau} & \omega_{G^D,\sigma^{-1}\tau}^{(p)}\\
(\omega_{G,\sigma^{-1}\tau})^{(p),\vee} & \mathcal E_{\tau} &   \omega_{G^D,\sigma^{-1}\tau}^{(p)}\\
(\omega_{G,\sigma^{-1}\tau})^{(p),\vee}  & (\omega_{G,\tau})^{\vee}   & \\
 };

\path[->,font=\scriptsize] 
(m-1-2) edge node[auto] {} (m-2-2)
(m-1-2) edge node[auto] {$V$} (m-1-3)
(m-2-1) edge node[auto] {$F$} (m-2-2)
(m-2-2) edge node[auto] {$V$} (m-2-3)
(m-2-2) edge node[auto] {} (m-3-2)
(m-3-1) edge node[auto] {$F$} (m-3-2)
;
\draw (m-1-3) edge[double distance=2pt] (m-2-3);
\draw (m-2-1) edge[double distance=2pt] (m-3-1);
\end{tikzpicture}
\end{center}
On note \[0 \fleche C_1 \fleche C_2\fleche C_2\fleche 0\] ce complexe.
On a alors $\det(C_1) \otimes \det(C_3) \simeq \det(C_2).$ Or $\det(C_1) \simeq 
\omega_{G^D,\sigma^{-1}\tau}^{\otimes p}\otimes\omega_{G^D,\tau}^{\vee}$ et $\det(C_3) \simeq 
\omega_{G,\tau}\otimes\omega_{G^D,\sigma^{-1}\tau}^{\vee,\otimes p}$. 
De plus, comme $C_2$ est exact, car $G$ est un groupe de Barsotti-Tate tronqué, 
$\det(C_2) \simeq \mathcal O_S$, et donc on a un isomorphisme,
\[ \phi_\tau : \omega_{G^D,\sigma^{-1}\tau}^{\otimes p}\otimes\omega_{G^D,\tau}^{\vee} \overset{\simeq}{\fleche} 
\omega_{G,\sigma^{-1}\tau}^{\otimes p}\otimes\omega_{G,\tau}^{\vee}.\]
En combinant les isomorphismes $\phi_\tau$ pour tout $\tau$, on obtient un isomorphisme,
\[ \phi_{\sigma^{1-f}\tau}^{\otimes(p^{f-1})}\circ\phi_{\sigma^{2-f}\tau}^{\otimes(p^{f-2})}\circ\dots\circ\phi_\tau : 
\omega_{G^D,\tau}^{\otimes (p^f-1)} \overset{\simeq}{\fleche} \omega_{G,\tau}^{\otimes(p^f-1)}.\]
\edem

\rem
Cette démonstration n'utilise pas l'hypothèse sur $r$, elle marcherait pour $r=1$.
\erem
\pro
\label{prodet}
Notons $d = \sum_\tau p_\tau$ la dimension de $G$. Supposons $r > d$. 
Considérons $\mathcal E_\tau$, le facteur direct du cristal de $G$. 
Alors le Frobenius-cristal,
\[ \det(\mathcal E_\tau,V^f),\]
se factorise en,
\[ (\mathcal O_{S/\Sigma}/p^r,p^d)\otimes(\mathcal L_\tau,\phi),\]
où $\mathcal L_\tau$ est un $(\mathcal O_{S/\Sigma})/p^{r}$-module localement libre de rang 1 et 
$(\mathcal L_\tau,\phi)$ est un cristal unité, c'est-à-dire que $\phi$ est un isomorphisme, 
unique modulo $p^{r-d}$. On peut voir $\phi$ comme une section de $\det(\mathcal E_\tau)^{\otimes(p^f-1)}$.
En particulier, si $G$ est un $\mathcal O$-module $p$-divisible sur un schéma $S$ sur lequel $p$ nilpotent, au-dessus de $\mathcal O$, alors on peut factoriser,
\[\det(\mathcal E_\tau,V^f) = (\mathcal O_{X/\Sigma},p^d)\otimes_{\mathcal O} L_\tau,\]
où $L_\tau$ est un $\mathcal O$-système local de rang 1, localement trivial pour la topologie étale.
\epro

\dem
On va montrer le résultat sur $\mathcal{BT}_{r}$, le champ de la section \ref{sect8}, sur lequel on note (encore) $G$ le groupe universel.
Tout d'abord pour montrer la première assertion, supposons $r > d$, et il suffit de le faire après tirée en arrière par la présentation $X \fleche \mathcal{BT}_r$, 
puis de montrer que la factorisation redescend.
On a donc $G/X$ un $\mathcal O$-module $p$-divisible tronqué d'échelon $r$ sur $X$ lisse sur $\mathcal O/p = \kappa_F$.
Dans ce cas si on considère le cristal $\det \mathcal E_\tau$, alors par le même raisonnement que la proposition \ref{pro34}, on a l'inclusion de faisceau dans $(X/\Sigma)^{cris}$, $\Sigma = \Spec(\ZZ_p)$,
\[ \im\left(\det(\Fil \mathcal E_{\tau'}) \fleche \det(\mathcal E_{\tau'})\right) \subset J_{X/\Sigma}^{h-q_{\tau'}}\det(\mathcal E_{\tau'}).\]
On en déduit par le même raisonnement que le corollaire \ref{cor33}, que \[ \det V^f : \det \mathcal E_\tau \fleche \det \mathcal E_\tau^{(p)},\]
est divisible (sur les épaississements $p$-adiques $(U,T,\delta)$) par $p^d$, puisque,
\[ d = \sum_\tau p_\tau = \sum_\tau h - q_\tau.\]
Mais comme $\det \mathcal E_\tau$ et $\det \mathcal E_\tau^{(p)}$ sont des cristaux, et que $X \otimes \kappa_F$ est lisse, on peut localement le relever en un schéma lisse sur 
$\mathcal O$, c'est-à-dire qu'il existe un recouvrement affine $(U \subset X)$ et des épaississements $p$-adiques $(U,T,\delta)$ avec $T/\Sigma$ affine lisse et 
$U = T \otimes \kappa_F$.
Dans ce cas le théorème 6.6. de \cite{BerO} nous dit que notre morphisme $\det V^f : \det \mathcal E_\tau \fleche \det \mathcal E^{(p)}$ est équivalent à un morphisme entre deux modules libres sur $\mathcal O_T/p^r$ munis d'une connexion, et comme l'épaississement est $p$-adique, $\det V^f$ est divisible par $p^d$.

Il existe donc un unique morphisme de cristaux $\phi_\tau : \det \mathcal E_\tau \pmod{p^{r-d}\mathcal O_{X/\Sigma}} \fleche \det \mathcal E_{\tau}^{(p)} \pmod{p^{r-d}\mathcal O_{X/\Sigma}}$, tel que $p^d\phi_\tau = V^f$ 
(la compatibilité à la connection est comme dans le théorème \ref{thruni}).

Montrons que $\phi_\tau$ est un isomorphisme. Il suffit de le voir fibres à fibres, et même sur les points géométriques.
Considérons alors $k$ un corps algébriquement clos, alors l'application $\phi_\tau$ sur le cristal se calcule sur le module de Dieudonné sur $W(k)$ correspondant, et dans ce cas c'est
simplement la théorie des diviseurs élémentaires.
L'unicité de $\phi_\tau$ (modulo $p^{r-d}$) assure que la construction redescend le long de $X \fleche \mathcal{BT}_r$ (voir lemme \ref{lemqcoh} et théorème \ref{thrdes}).
On a donc une (unique modulo $p^{r-d}$) décomposition sur $\mathcal{BT}_r$, 
\[\det(\mathcal E_\tau,V^f) = (\mathcal O_{\mathcal{BT}_r/\Sigma},p^d)\otimes_{\mathcal O} (\mathcal L_\tau,\phi_\tau),\]
où $\mathcal L_\tau$ est un $\mathcal O_{\mathcal{BT}_r/\Sigma}$ inversible, muni d'un isomorphisme,
\[\phi_\tau : \mathcal L_\tau \fleche \mathcal L_\tau^{(p^f)},\]
C'est-à-dire un cristal unité. Il est donc associé à un $\mathcal O$-système local de rang 1, $L_\tau$.
Le cas d'une base générale $S$ s'en suit par tiré en arrière à partir de $\mathcal{BT}_r$.
\edem

\rem
Le $\FF_{p^f}$-système local de rang 1 $L_\tau \pmod p$ a donc une monodromie à valeurs dans $\FF_{p^f}^\times$, et donc $L_\tau^{\otimes p^f-1}$ est un système local de rang 
1 constant, or on a la filtration de Hodge,
\[ 0 \fleche \omega_{G^D,\tau} \fleche (\mathcal E_\tau)_{\mathcal{BT}_{d+1} \overset{\id}{\fleche} \mathcal{BT}_{d+1}} \fleche \omega_{G,\tau}^\vee \fleche 0,\]
de telle sorte que $\det(\mathcal E_\tau)_{\mathcal{BT}_{d+1}\overset{\id}{\fleche} \mathcal{BT}_{d+1}}^{\otimes(p^f-1)}$ soit trivial. 
On en déduit que si $S/\FF_{p^f}$ est un schéma, et $G$ un $\mathcal O$-module de 
Barsotti-Tate tronqué d'échelon $d+1$ sur $S$, on a donc un isomorphisme, \cite{KnudMum},
\[ \det(\omega_{G^D,\tau})^{\otimes(p^f-1)} \otimes \det(\omega_{G,\tau})^{\otimes(1-p^f)} \overset{\simeq}{\fleche} O_S,\]
et donc on retrouve la proposition \ref{pro131}, sous une hypothèse en plus sur $r$ (la démonstration de la proposition \ref{pro131} fonctionnant pour $r = 1$)
\[\det(\omega_{G^D,\tau})^{\otimes(p^f-1)} \simeq \det(\omega_{G,\tau})^{\otimes(p^f-1)}.\]
\erem

\lem
Soit $x_r \dans \mathcal{BT}_r(k)$, correspondant à un $\mathcal O$-module $p$-divisible $\underline{G_r}$ tronqué d'échelon $r$, alors il existe une préimage $x \dans \mathcal{BT}_\infty(k)$ de 
$x_r$ par la flèche,
\[ \mathcal{BT}_\infty \overset{[p^r]}{\fleche} \mathcal{BT}_r.\]
C'est-à-dire qu'il existe $\underline{G}/k$ un $\mathcal O$-module $p$-divisible (non tronqué) tel que $\underline{G}[p^r] = \underline{G_r}$.
\elem

\dem
Comme $\mathcal O$ est non ramifié, c'est \cite{Wed2} proposition 3.2. 
\edem

\lem
\label{lem106}
Soit $x_r,x,G$ comme précédemment. Soit $y_r$ une préimage de $x_r$ par la présentation $X \fleche \mathcal{BT}_r$. Regardons alors le produit fibré,

 \begin{center}
\begin{tikzpicture}[description/.style={fill=white,inner sep=2pt}] 
\matrix (m) [matrix of math nodes, row sep=3em, column sep=2.5em, text height=1.5ex, text depth=0.25ex] at (0,0)
{ 
D &  \widehat{X}_{y_r} \\
\Def_{\underline G} \otimes k & \mathcal{BT}_r  \\
};

\path[->,font=\scriptsize] 
(m-1-1) edge node[auto] {$\pi$} (m-1-2)
(m-1-1) edge node[auto] {$p$} (m-2-1)
(m-2-1) edge node[auto] {$$} (m-2-2)
(m-1-2) edge node[auto] {} (m-2-2)
;
\end{tikzpicture}
\end{center}
Alors les flèche $p,\pi$ sont lisses surjectives. En particulier, les sections $\widetilde{^\mu\Ha(G^{univ})}$ et $\widetilde{^\mu\Ha(G^{univ,D})}$ coincident sur $\widehat{X}_y$ si et seulement si elles coincident sur $D$.
\elem

\rem
On va donc pouvoir ramener le calcul des invariants de Hasse au calcul sur
$\Def_{\underline{G}}$, où l'on a des coordonnées explicites. 
En particulier, pour montrer que $\widetilde{\Ha}(G[p^r]) =\widetilde{\Ha}(G[p^r]^D)$, il suffira de le faire sur $\Def_{\underline{G}}$.
\erem

\dem
$D$ s'identifie aussi à la complétion formelle du produit fibré de $X \times_{\mathcal{BT}_r} \Def_{\underline G}$, et $X \fleche \mathcal{BT}_r$ est surjective, donc $p$ aussi.
D'après Illusie, Wedhorn \cite{Wed2} la flèche $\mathcal{BT}_\infty \fleche \mathcal{BT}_r$ est formellement lisse, et donc $\pi$ aussi. De plus $\pi$ est surjective puisque $\widehat{X_{y_r}}$ n'a qu'un point.
%
\edem

\thr
\label{thrdual}
Supposons fixée une signature $(p_\tau,q_\tau)_\tau$, choisissons un $\tau \in \mathcal I$, et notons 
\[k_{\tau} =  \sum_{\tau' : q_{\tau'} < q_\tau} q_\tau - q_{\tau'} \quad \text{et} \quad k_{\tau}^D =  \sum_{\tau' : p_{\tau'} < p_\tau} p_\tau - p_{\tau'}.\] 
Soit $r > \max(d,k_\tau,k_\tau^D)$. Sur $\mathcal{BT}_r \otimes \kappa_F$, sous l'isomorphisme,
\[ \det\omega_{G^D,\tau}^{\otimes(p^f-1)} \simeq \det\omega_{G,\tau}^{\otimes(p^f-1)},\]
on a l'égalité entre les sections $\widetilde\Ha_\tau(G)$ et $\phi \otimes \widetilde\Ha_\tau(G^D)$, où $\phi$ est une section inversible déterminée par la proposition précédente.
Ceci induit bien sûr une égalité, à une section inversible prés, dans $H^0(\mathcal{BT}_r,\det(\omega_{G^D,q_\tau})^{(p^f-1)})$ des sections $\widetilde{\Ha_{q_\tau}(G)}$ et 
$\widetilde{\Ha_{p_\tau}}(G^D)$, ainsi que des sections $\widetilde{^\mu\Ha}(G)$ et $\widetilde{^\mu\Ha}(G^D)$, (pour $r$ est suffisamment grand).

En particulier on en déduit que si $G$ est un $\mathcal O$-module $p$-divisible tronqué d'échelon $r$ sur $\mathcal O_C$ et $r > \max(d,k_\tau(G),k_\tau(G^D))= \max(d,k_\tau, k_\tau - d + fq_\tau)$, alors, \[ \Ha_\tau(G) = \Ha_\tau(G^D).\]
\ethr

\dem
Identifions les deux fibrés $\det\omega_{G^D,\tau}^{\otimes(p^f-1)}$ et $\det\omega_{G,\tau}^{\otimes(p^f-1)}$. Comme $X \fleche \mathcal{BT}_r$ est lisse surjective, il suffit de voir l'égalité sur $X$. De plus, comme $X$ est lisse, donc réduit, il suffit de voir que 
les deux sections coïncident sur un ouvert dense, par exemple le lieu $\mu$-ordinaire, cf \cite{GW} Corollary 9.9 par exemple.
De plus, il suffit de voir que les sections coïncident en chaque anneau local,
\[ \mathcal O_{X,y_r},\]
pour tout point géométrique $y_r = \Spec(k)$ au dessus de $x_r$ disons, $\mu$-ordinaire, de $X$.
Or la flèche\[\mathcal O_{X,y_r} \fleche \widehat{\mathcal O_{X,y_r}},\] 
est injective (puisque $X$ est localement noetherien (Artin-Rees)), il suffit donc de voir l'égalité sur,
\[ \Spf(\widehat{\mathcal O}_{X,y_r}).\]
Mais d'après le lemme précédent, si on note $\mathbb G$ un relèvement à $\mathcal{BT}_\infty$ de $x_r$, et $D = (\Def_{\mathbb G} \otimes_{\mathcal O_E} k) \times_{\mathcal BT_r} \widehat{X}_{y_r}$, les sections coïncident sur $ \widehat{X}_{y_r}$ si elles coïncident sur $D$. Or on a une flèche,
\[ \Phi :  D \fleche \Def_{\mathbb G} \otimes_{\mathcal O_E} k,\]
qui identifie les deux $\mathcal O$-module tronqués universels sur $D$ (ceux provenant de $X$ et $\Def_{\mathbb G}$).
Or on a un isomorphisme (théorème \ref{thr121}),
\[\Def_{\mathbb G} \otimes_{\mathcal O_E} k \simeq \Spf(k[[t_{k_i,l_i}^i,i = \{1,\dots,f\}, k_i \dans \{1,\dots,p_{\sigma^i\tau}\}, l_i \dans \{1,\dots,q_{\sigma^i\tau}]]),\] 
pour lequel on a un épaississement à puissances divisées $\hat W(k[[t_{k_i,l_i}^i]])$ de $W(k)$, et sur lequel on connaît le cristal de la déformation universelle (cf proposition 
\ref{prodefouniv}).
Comme $G^{univ} \times \Spec(k) = \mathbb G$ est $\mu$-ordinaire, son cristal (i.e. son module de Dieudonné) se décompose et sa matrice de Verschiebung s'écrit 
(d'après le théorème de Moonen \ref{thrmoo}) dans une base adaptée à la décomposition $P_0 = \bigoplus_{i=0}^f P_0^{\sigma^i\tau}$ et à la filtration de Hodge,
\[
A^\flat = \left(
\begin{array}{cccc}
0  &   & & A^\flat_0   \\
A^\flat_1  & 0  & &  \\
  & \ddots  &   \ddots &  \\
 0  & & A^\flat_{f-1} &   0  \\\end{array}
\right), \quad A^\flat_i = 
\left(
\begin{array}{cc}
\Id_{q_{\sigma^i\tau}}  &      \\
  &   p\Id_{p_{\sigma^i\tau}}  
\end{array}
\right).
\]

Notons $P = \bigoplus P_\tau$ le cristal sur $\hat W(k[[t_i]])$ du groupe universel $G$, alors sa matrice de Verschiebung donnée par \cite{Zink} (87), et la proposition \ref{prodefouniv}, est,
\[ A = A^\flat\Diag(
\left(
\begin{array}{cc}
  \Id_{p_{\sigma^i\tau}}  & -[t_{k,l}^i]   \\
0  &   \Id_{q_{\sigma^i\tau}}   
\end{array}
\right), i = 0,\dots,f-1),\]
où $A^\flat$ est le relèvement trivial de la matrice du Verschibung de $G\times \Spec(k)$ à $\hat W(k[[t_i]])$.
De plus, si on note 
\begin{equation}
\label{mat}A_{i}^{(\phi^{f-1})}A_{i+1}^{(\phi^{f-2})}\cdots A_{i-2}^{(\phi)}A_{i-1} =
\left(
\begin{array}{cc}
D^i_1  &D^i_2   \\
D^i_3  &D^i_4  
\end{array}
\right),\end{equation}
avec $D^i_1 \dans M_{q_{\sigma^i\tau}}(\hat W(k[[t_j]]))$, alors $\Ha_{\sigma^i\tau}(G)$ est donné par $p^{-k_{\sigma^i\tau}}\det D^i_1$ (si $r > k_{\sigma^i\tau}$), 
d'après la proposition (\ref{prodefouniv}).

On peut calculer effectivement $A^{(\phi^{f-1})}A^{(\phi^{f-2})}\cdots A^{(\phi)}A$, et on trouve une matrice diagonale dont le $i$-eme bloc est de la forme,
\[ \left(
\begin{array}{cccccccc}
1 & \\
& \ddots &\\
&&1 &&&& (*)\\
&&&p\\
&&&& \ddots & \\
&0&&&&p\\
&&&&&&\ddots\\
&&&&&&&p^t
\end{array} 
\right),\quad \text{où } t = f - |\{ \tau : p_\tau = 0\}|.\]
On peut donc préciser l'écriture (\ref{mat}), qui est donc de la forme,
\[\left(
\begin{array}{cc}
D^i_1  &D^i_2   \\
0 &D^i_4  
\end{array}
\right)
 ,\]
 où $\det D_1^i$ est divisible par $p^{k_{\sigma^i\tau}}$ si $r > k_{\sigma^i\tau}$. Appliquons cela à $\tau$, c'est à dire $i = 0$, on en déduit donc que 
 $\det V^f_{|\mathcal{P}_{\tau}} = \det D_1^0\det D_4^0$.
Or on sait d'après la proposition précédente que $\det V^f_{|\mathcal P_{\tau}} = p^ds$ où $s$ est un inversible dans $\hat W(k[[\underline t]])$ (c'est la spécialisation de la 
section inversible sur $\mathcal{BT}_r$ et d'un choix d'un isomorphisme $\det(\omega_{G^D,\tau})^{(p^f-1)} \simeq \det(\omega_{G,\tau})^{(p^f-1)}$).
Mais on a que 
\[\det D_1^0 = p^{k_{\tau}}\widetilde{\Ha}_{\tau}^0(G)\quad \text{et}  \quad \det D_4^0 = p^{p_{\tau}f - k_{\tau}(G^D)}\widetilde\Ha^0_{\tau}(G^D)^{-1},\]
où $\widetilde\Ha^0_{\tau}(G) \equiv \widetilde\Ha_{\tau}(G) \pmod{(p)}$ et idem pour $G^D$.
On en déduit donc qu’au-dessus de $\Def_{\mathbb G}$,
\[ p^{k_{\tau}}\widetilde\Ha^0_{\tau}(G)=p^dsp^{k_{\tau}(G^D)-fp_{\tau}}\widetilde\Ha^0_{i\tau}(G^D),\]
Mais un calcul direct donne que pour tout $\tau'$, $k_{\tau'}(G) - k_{\tau'}(G^D) = d - fp_{\tau'}$. Autrement dit, en divisant par $p^{k_{\tau}}$ et en réduisant modulo $p$, 
on en déduit que sur $\Def_{\mathbb G}$ et donc sur $D$ en tirant en arrière par $\Phi$, puis par le lemme \ref{lem106}, les anneaux locaux de $X$,
\[ \widetilde\Ha_{\tau}(G) = s \otimes \widetilde\Ha_{\tau}(G^D),\]
où $s$ est la section inversible dans $H^0(\mathcal{BT}_r,\mathcal O_{\mathcal{BT}_r})$ qui correspond au $\phi$ de la proposition \ref{prodet}.
En particulier, comme cette égalité est vraie pour les anneaux locaux en chaque point d'un ouvert Zariski dense, on en déduit que (comme diviseurs de Cartier) 
$\widetilde\Ha_\tau(G) = \widetilde\Ha_\tau(G^D)$ 
sur tout $\mathcal{BT}_r$, pour $r$ assez grand, et donc par la proposition \ref{proechelon} dès lors que $r > \max(d,k_\tau(G),k_\tau(G^D))$, alors, quelque soit la base, 
au dessus de $\kappa_F$, si $G$ est un $\mathcal{BT}_r$,
\[ \widetilde{\Ha_\tau}(G) = \widetilde{\Ha_\tau}(G^D).\qedhere\]
\edem

En particulier, si $S = \Spec(O_K)$, $K$ une extension de $\QQ_p$, $v$ sa valuation normée par $v(p)=1$ (par exemple $C = \widetilde{\overline{\QQ_p}}$), si on trivialise sur $O_K$ $\det(\omega_{G^D}^{\otimes(p^f-1)})$ et $\det(\omega_{G^D}^{\otimes(p^f-1)})$, et que dans ces bases on note $^\mu\Ha(G),^\mu\Ha(G^D) \dans [0,1]$ (et si $K$ contient $W(\kappa_F)$, resp. $\Ha_\tau$) les valuation des élément $\widetilde{^\mu\Ha}(G)$ et  $\widetilde{^\mu\Ha}(G)$ (resp.  $\widetilde{\Ha}_\tau$), qui ne dépendent pas des trivialisation, on a l'égalité,
\[ ^\mu\Ha(G) =  {^\mu\Ha}(G^D) \quad \text{(resp. } \quad \Ha_\tau(G) = \Ha_\tau(G^D)).\]

\section{$\mathcal O$-modules $p$-divisibles $\mu$-ordinaires sur $\mathcal O_C$}

\label{sect12}

Soit $\mathcal O = \mathcal O_F$ l'anneau des entiers d'une extension $F/\QQ_p$ non ramifiée, et $\tau \in \Hom(F,C)$ un plongement. Alors il existe un groupe $p$-divisible (strict) muni d'une $\mathcal O$-action au-dessus de $\mathcal O$ (et donc sur toutes les extensions de $\mathcal O$), noté $\mathcal{LT_\tau}$, de $\mathcal O$-hauteur 1 et tel que
$\omega_{\mathcal{LT_\tau}}$ soit de dimension 1, et sur lequel l'action de $\mathcal O$ soit donnée par $\tau$.
Alors pour tout $A \subset \Hom(F,C)$, on a construit, cf. définition \ref{LTA},
\[ "\mathcal{LT}_A = \bigotimes_{\tau \in A}^{\mathcal O_K} \mathcal LT_\tau",\]
où le produit tensoriel est à prendre au sens du produit tensoriel des displays, et provient d'un $\mathcal O$-module $p$-divisible de $\mathcal O$-hauteur 1.
Sur $k=\overline{\FP}$, $\mathcal{LT_A}\otimes_{\mathcal O} k$ a pour module de Dieudonné, $(W(k)^f,V)$, où
\[ V =
\left(
\begin{array}{cccc}
 0 &   &   & \delta_{\sigma\tau}\\
 \delta_\tau& 0  &   \\
0  &  \delta_{\sigma^{-1}\tau} &   0\\
 0 & 0  & \ddots &  0\\
  0 &   &0 & \delta_{\sigma^{f-1}\tau} \\
\end{array}
\right), \quad \delta_{\tau'} = 
\left\{
\begin{array}{cc}
 p & \text{si } \tau \in A  \\
  1 & \text{sinon}
\end{array}
\right.\]
En particulier il est isocline de pente $\frac{|A|}{f}$.

On a alors le lemme suivant,

\pro
\label{pro26}
Si $G/\mathcal O_C$ est un groupe $p$-divisible avec une action de $\mathcal O$, où $\mathcal O = \mathcal O_F$ et $F/\QQ_p$ est non ramifiée, et soit $(p_\tau,q_\tau)$ sa signature,
et $h$ sa $\mathcal O$-hauteur. 
Alors $G$ est $\mu$-ordinaire si et seulement si
\[ G \simeq X^{ord}_{(q_\tau)} := \prod_{l =1}^{r+1} \mathcal{LT}_{A_{l-1}}^{q^{(l)}-q^{(l-1)}},\]
où l'on note $\{q_\tau : \tau \in \Hom(K,C)\} = \{ q^{(1)},\dots, q^{(r)}\}$, $q^{(0)} = 0$ et $q^{(r+1)} = h = p_\tau + q_\tau$, et pour tout $l$,
\[ A_l = \{ \tau \in \mathcal I : q_\tau \leq q^{(l)}\} \quad \text{et} \quad A_0 = \emptyset\]
\epro

\dem
On sait d'après un théorème de Moonen (\ref{thrmoo}) que c'est le cas sur $\mathcal O_C/\mathfrak m_C = k$ qui est algébriquement clos de caractéristique $p$.
De plus, toujours par Moonen \cite{Moo}, ou par Shen \cite{Shen}, puisque les polygones de Hodge et Newton coïncident, pour tout $q^{(l)}$ (i.e. tous les points de rupture) on a un 
sous-$\mathcal O$-module de $G$, et donc on a une filtration,
\[ 0 \subset G_1 \subset \dots \subset G_r = G,\]
telle que les gradués $E_l = G_l/G_{l-1}$, $l \dans \{2,\dots,r\}$, sont des $\mathcal O$-modules. Par Moonen, 
\[E_l \otimes_{\mathcal O_C} k \simeq \mathcal{LT}_{A_{l-1}}^{n_{l-1}}, \quad n_l = q^{(l-1)} - q^{(l-2)}.\]
Donc $E_l$ est une déformation sur $\mathcal O_C$ d'un $\mathcal{LT}_A^n$. Mais une telle déformation est triviale ; en effet, d'après \cite{SW} elle correspond à un triplet $(T,W,\alpha)$ muni 
d'une action de $\mathcal O$, c'est à dire comme on connaît la signature et la hauteur, $n$, à $(\mathcal O^n,\bigoplus_{\tau \not\in A} C_\tau^n,\bigoplus_{\tau} \alpha_\tau)$.
Or la flèche $\alpha_\tau : \mathcal O^n \otimes_{\ZZ_p} C \fleche W_\tau$ est surjective, donc quitte à changer de base $W_\tau$ 
(c'est-à-dire à faire un isomorphisme dans la catégorie des triplets), on peut supposer que le $\mathcal O$-triplet est de la forme,
\[(\mathcal O^n, \bigoplus_{\tau \not\in A} C_\tau^n, \bigoplus_{\tau \not\in A} \tau^{\oplus n}).\]
En particulier il n'existe à isomorphisme près sur $\mathcal O_C$ qu'un seul tel triplet, et donc $E_l \simeq \mathcal{LT}_A^n$.
Donc $G$ est filtré par ce que l'on souhaite, il nous reste seulement à montrer que toute telle filtration est scindée.
Mais encore une fois utilisons la description de Scholze-Weinstein de la catégorie des groupes $p$-divisibles sur $\mathcal O_C$, cf. \cite{SW} Theorem B, et on se retrouve avec 
une filtration par des triplets $(T,W,\alpha)$ munis d'une action de $\mathcal O$, telle que le gradué est donné par des (puissances de) triplets de Lubin-Tate généralisés, et le lemme suivant nous 
dit que la filtration est scindée.
\edem

\lem
\label{lem27}

Soit $C = (T,W,\alpha)$ un triplet de Scholze-Weinstein muni d'une action compatible de $\mathcal O$.
Supposons qu'il existe une filtration de $C$ par des $\mathcal O$-triplets, 
\[C_0 = 0 \subset C_1 = (T_1,W_1,\alpha_1) \subset C_2 = (T_2,W_2,\alpha_2) \subset \dots \subset C_r = (T_r,W_r,\alpha_r) = (T,W,\alpha),\]
telle que $C_{i+1}/C_i \simeq (T_{A_i},W_{A_i},\alpha_i)^{n_i}$ est le triplet de Scholze-Weinstein de $\mathcal{LT}_{A_i}^{n_i}$, et $A_1 \supset A_2\supset \dots \supset A_r$. Alors, $C$ est l'extension triviale des $C_{i+1}/C_i$.
De plus l'application $\alpha : T \fleche W$ est $E$-rationnelle, où $E$ est la composée de tous les $E_{A_i}$.
\elem

\dem
Sans perte de généralité, puisque les $T_i$ sont libres, on peut écrire $T = \mathcal O^h, T_i = \mathcal O^{m_i}$.
Le triplet associé à $\mathcal{LT}_A$ est donné par
\[ (\mathcal O, W = \bigoplus_{\tau \not\in A} C_\tau, \alpha_\tau : \mathcal O \fleche C_\tau = \tau).\]
Donc $(T_2, W_2)$ est l'extension de $\mathcal{LT}_{A_1}^{n_1}$ par $\mathcal{LT}_{A_{2}}^{n_{2}}$, et donc est donné par,
\[ (\mathcal O_K^{n_1 + n_{2}}, \bigoplus_{\tau} C^{m_\tau}, (\alpha_\tau)_\tau),\]
où \[m_\tau = 
\left\{
\begin{array}{cc}
  n_1 + n_{2} & \text{ si}\tau \not\in A_{2}   \\
  n_1 & \text{ if}\tau \not\in A_1 \priv A_{2}    \\
 0  &   \text{ sinon}   
\end{array}
\right.
\] et \[\alpha_\tau : \mathcal O^{n_1 + n_{2}} \fleche C^{m_\tau}
\left\{
\begin{array}{cc}
\left(
\begin{array}{cc}
 I_{n_1} & M_\tau    \\
    0 &   I_{n_{2}}
\end{array}
\right)
& \text{ si } \tau \not\in A_{1}   \\
\left(
\begin{array}{cc}
 0_{n_1} & 0    \\
    0 &   I_{n_{2}}
\end{array}
\right)
  & \text{ si }\tau \in A_1 \priv A_{2}    \\
 0 &   \text{ sinon}   
\end{array}
\right.\]
Mais dans le cas non scindé (i.e. $\tau \not\in A_{1}$), en modifiant la base de $W_\tau$, grâce à l'isomorphisme de triplet
\[(Id,
\left(
\begin{array}{cc}
 Id & -M_\tau    \\
  & Id     
\end{array}
\right)_\tau)
\]
on peut rendre la filtration scindée (i.e. $M_\tau = 0$).
On peut raisonner identiquement (ou par récurrence) pour montrer \[(T_i,W_i) = \prod_{k=1}^i C_{k}/C_{k-1}.\]
\edem

\rem
C'est faux sans l'hypothèse $A_1 \supset A_2 \supset \dots \supset A_r$ : il existe (sur $\mathcal O_C$) des extensions dans la catégorie des triplets de Scholze-Weinstein de 
$\mu_{p^\infty}$ par $\QQ_p/\ZZ_p$ qui ne sont pas triviales (ni ordinaires donc) !
Le problème dans l'algorithme précédent est que , en reprenant les notations, si on a $\tau \dans A_2 \priv A_1$, et donc,
\[ \alpha_\tau : \mathcal O^{n_1+n_2} \fleche C^n_2,\]
qui est donnée par la matrice $(I_{n_2},M_\tau)$, il faut pouvoir faire le morphisme de triplets,
\[ 
\left(
\begin{array}{ccc}
 I_d   &   -M_\tau  \\
&  I_d    
\end{array}
\right)
\]
pour rendre $\alpha_\tau$ de la forme $(Id,0)$, ce qui n'est possible que si $M_\tau$ est à coefficients dans $\mathcal O[1/p]$ (et qui ne dépend pas de $\tau$).
\erem

\rem
Le résultat de cette section peut s'interpreter géométriquement de la manière suivante. Soit $\mathcal S$ une variété de Shimura PEL sur $C$ (disons celle de la section 
\ref{sect9}, mais plus généralement PEL non ramifiée en $p$). Regardons $S_\infty \fleche S$ le revêtement adique pro-étale perfectoide associé à $\mathcal S$ 
par Scholze dans \cite{Schtor}. Alors Caraiani et Scholze \cite{SchCar} ont construit une application,
\[ \pi_{HT} : S_\infty \fleche \mathcal{F}\ell,\]
où $\mathcal F\ell$ est une certaine Grassmanienne, qui raffine l'application de Hodge-Tate de \cite{Schtor}. Leur construction est en fait beaucoup plus générale, et remarquons que dans le cas considéré, avec l'aide des structures entières, celle-ci peut se construire de manière beaucoup plus simple.
On peut alors interpreter le théorème précédent en disant que le lieu $\mu$-ordinaire de $S_\infty$ est envoyé sur les points $E$-rationnels de $\mathcal{F}\ell$ (où $E$ est le corps réflexe $p$-adique associé à $(p_\tau,q_\tau)$), ce qui est en fait un cas particulier de la compatibilité entre strates de Newton de $S$ et de $\mathcal{F}\ell$ de \cite{SchCar} (dans le cas extrêmement 'basique'–si l'on ose dire, bien qu'il ne s'agisse pas de la strate basique, bien au contraire ! – de la strate $\mu$-ordinaire).
\erem

\appendix

\section{Rappels sur les Champs et modules quasi-cohérents}
\label{appA}

Les résultats de ces appendices n'ont rien d'originaux, mais il semble difficile d'en trouver des références précises dans la littérature, c'est pourquoi on a préférer les réécrire. On renvoie à \cite{LMB} et \cite{Stackproj} pour les notions utilisées ici.
Soit $S$ un schéma de base. On munit $Sch/S$ de la topologie lisse 
Soit $\mathfrak X \fleche Sch/S$ un champ algébrique, et soit \[p : X \fleche \mathfrak X,\]
une présentation par un schéma $X$.

\defi
\label{defqcoh}
Un $\mathcal O_{\mathfrak X}$-module quasi-cohérent $\mathcal E$ sur $\mathfrak X$ est la donnée, pour tout schéma $U$ dans $Sch/S$, et pour tout $u \dans \mathfrak X_U = \mathfrak X(U)$, i.e. pour tout morphisme,
\[ u : U \fleche \mathfrak X,\]
d'un $\mathcal O_U$-module quasi-cohérent $\mathcal E_u$ sur $U$, tel que, pour toute flèche 
$U' \overset{f}{\fleche} U$, on ait un isomorphisme,
\[ \alpha_f : f^* \mathcal E_u \overset{\simeq}{\fleche} \mathcal E_{f^*u}.\]
De plus, ces isomorphismes doivent vérifier la condition de cocycle, c'est-à-dire, si on a des 
flèches,
\[ U'' \overset{g}{\fleche} U' \overset{f}{\fleche} U,\]
alors, le diagramme suivant commute,
\begin{center}
\begin{tikzpicture}[description/.style={fill=white,inner sep=2pt}] 
\matrix (m) [matrix of math nodes, row sep=3em, column sep=2.5em, text height=1.5ex, text depth=0.25ex] at (0,0)
{ 
\mathcal E_{(f\circ g)^*u}  &  & g^*\mathcal E_{f^*u} & & g^*f^*\mathcal E_u \\
 &  & (g\circ f)^*\mathcal E_u & &  \\
 };

\path[->,font=\scriptsize] 
(m-1-1) edge node[auto] {$\alpha_g$} (m-1-3)
(m-2-3) edge node[auto,right] {can} (m-1-5)
(m-1-3) edge node[auto] {$g^*\alpha_f$} (m-1-5)
(m-1-1) edge node[auto,left] {$\alpha_{f\circ g}$} (m-2-3)
;
\end{tikzpicture}
\end{center}

que l'on écrit plus rapidement,
\[ \alpha_{f\circ g} = g^*\alpha_f \circ \alpha_g.\]

Un morphisme entre deux modules quasi-cohérents sur $\mathfrak X$, noté
\[\phi : \mathcal E \fleche \mathcal F,\]
est la donnée, pour tout $U \dans Sch/S$, pour tout $u \dans \mathfrak X_U$, d'un morphisme 
de $O_U$-modules quasi-cohérents,
\[ \phi_u : \mathcal E_u \fleche \mathcal F_u,\]
tel que pour toute flèche, $U' \overset{f}{\fleche} U$, le diagramme suivant commute,
\begin{center}
\begin{tikzpicture}[description/.style={fill=white,inner sep=2pt}] 
\matrix (m) [matrix of math nodes, row sep=3em, column sep=2.5em, text height=1.5ex, text depth=0.25ex] at (0,0)
{ 
\mathcal E_{f*u}  &  & \mathcal F_{f^*u} \\
f^*\mathcal E_{u}  &  & f^*\mathcal F_{u}  \\
 };

\path[->,font=\scriptsize] 
(m-1-1) edge node[auto] {$\phi_{f^*u}$} (m-1-3)
(m-2-1) edge node[auto] {$\alpha_{f}$} (m-1-1)
(m-2-3) edge node[auto] {$\beta_{f}$} (m-1-3)
(m-2-1) edge node[auto] {$f^*\phi_u$} (m-2-3)
;
\end{tikzpicture}
\end{center}
où l'on a noté $\alpha_\cdot$ les isomorphismes structuraux de $\mathcal E$, et $\beta_\cdot$ 
ceux de $\mathcal F$. On note $QCoh(\mathfrak X)$ la catégorie ainsi définides modules quasi-cohérents sur $\mathfrak X$.
\edefi

\pro
\label{proqcoh}
Soit $\mathfrak X \fleche Sch/S$ un champ algébrique, et soit,
\[ p : X \fleche \mathfrak X,\]
une présentation (lisse) par un $S$-schéma $X$.
Alors la catégorie $QCoh(\mathfrak X)$ est équivalente à la catégorie des modules 
quasi-cohérents sur $X$ munis d'une donnée de descente, notée $Des(X)$, 
c'est à dire, des modules quasi-cohérents $F$ sur $X$, munis d'un isomorphisme,
\[ \phi : p_1^*F \fleche p_2^*F,\]
tel qu'on ait la condition de cocycle,
\[ p_{23}^*\phi \circ p_{12}^*\phi = p_{13}^*\phi.\]
On a noté les projections sur les facteurs correspondants,

\[\xymatrix{ X\times_{\mathfrak X} X \times_{\mathfrak X} X\ar@<-16pt>[d]^{p_{12}}\ar@<0pt>[d]^{p_{13}}\ar@<16pt>[d]^{p_{23}} \\
  X\times_{\mathfrak X} X \ar@<-6pt>^{p_1}[d] \ar@<6pt>[d]^{p_2} \\
  X }\]
Les morphismes de $Des(X)$ sont les morphismes $h : F \fleche G$ tels que,
\[ p_2^*h \circ \phi_F = \phi_G \circ p_1^*h.\]
C'est à dire tels que le diagramme suivant commute.
\begin{center}
\begin{tikzpicture}[description/.style={fill=white,inner sep=2pt}] 
\matrix (m) [matrix of math nodes, row sep=3em, column sep=2.5em, text height=1.5ex, text depth=0.25ex] at (0,0)
{ 
p_2^*F &  & p_2^*G \\
p_1^*F &  & p_1^*G\\
 };

\path[->,font=\scriptsize] 
(m-1-1) edge node[auto] {$p_2^*h$} (m-1-3)
(m-2-1) edge node[auto] {$\phi_{F}$} (m-1-1)
(m-2-3) edge node[auto] {$\phi_{G}$} (m-1-3)
(m-2-1) edge node[auto] {$p_1^*h$} (m-2-3)
;
\end{tikzpicture}
\end{center}
\epro

\dem
C'est essentiellement une reformulation de la descente fpqc.
\edem

\rem
\label{remqcoh}
On a énoncé la proposition pour les modules quasi-cohérents, mais c'est encore vrai en remplaçant quasi-cohérent par localement libre de rang fini, localement libre de rang $n$ fixé (donc aussi inversible), d'après la proposition 2.5.2 de \cite{EGA4}, volume 2.
\erem

\exe
\label{exeinv}
Soit $D \overset{i}{\hookrightarrow} \mathfrak X$ un sous-champ fermé d'un champ algébrique. 
Alors il existe un $\mathcal O_{\mathfrak X}$-module quasi cohérent $\mathcal I_D$ associé à $D$ :
Soit $U \dans Sch/S$, et $u \dans \mathfrak X_U$ plat. Alors, on a le diagramme
\begin{center}
\begin{tikzpicture}[description/.style={fill=white,inner sep=2pt}] 
\matrix (m) [matrix of math nodes, row sep=3em, column sep=2.5em, text height=1.5ex, text depth=0.25ex] at (0,0)
{ 
U\times_{u,\mathfrak X,i}D  &  &D \\
U &  & \mathfrak X  \\
 };

\path[->,font=\scriptsize] 
(m-1-1) edge node[auto] {} (m-1-3)
(m-1-1) edge node[auto] {$j$} (m-2-1)
(m-1-3) edge node[auto] {$i$} (m-2-3)
(m-2-1) edge node[auto] {$u$} (m-2-3)
;
\end{tikzpicture}
\end{center}
et $j$ est une immersion fermée. On a donc un faisceau d'idéaux $\mathcal I_{D,u}$ sur $U$ 
associé à $U\times D$. Soit maintenant $U' \overset{f}{\fleche} U$ plat, alors grâce à 
l'identification canonique,
\[ f^*(U\times D) \cong U' \times_{u\circ f, \mathfrak X, i} D,\]
on en déduit un morphisme canonique $\alpha_f : \mathcal I_{D,u\circ f} \cong 
f^*\mathcal I_{D,u}$. Comme $\alpha_f$ est canonique, il vérifie la condition de cocycle, 
et on en déduit donc un faisceau quasi-cohérent de $\mathcal O_X$-module $\mathcal I_D$ sur 
$\mathfrak X$, associé à $D$, en considérant $U = X$ une présentation de $\mathcal X$ et $U' = X \times_{\mathcal X} X$ etc.. en utilisant la proposition précédente.
Néanmoins il n'est pas clair (et probablement faux) que pour tout $U \dans Sch/S$ et $u \dans \mathcal X_U$ on ait que
$\mathcal I_{D,u}$ soit l'ideal associé au fermé $D \times_u U$ de $U$. 
On voudrait dire que $\mathcal I_D$ est un faisceau d'idéaux, seulement on n'est pas certain que pour tout $U/S$ et tout $u \dans \mathfrak X_U$, $\mathcal I_{D,u} \subset \mathcal O_{\mathfrak X,u} = \mathcal O_u$, mais c'est le cas si $u$ est plat. Néanmoins on a une flèche $\mathcal I_D \fleche \mathcal O_{\mathfrak X}$, qui est par exemple donnée par l'inclusion de faisceaux sur une présentation de $\mathfrak X$ (voir proposition précédente), qui est un monomorphisme sur tout $u : U \fleche \mathfrak X$ plat. On dira qu'un $\mathcal O_{\mathfrak X}$-module $\mathcal I$ qui vérifie cette propriété est un faisceau d'idéaux sur $\mathfrak X$.
\eexe

\lem
\label{lemqcoh}
Si $\mathfrak X$ est un champ algébrique lisse, alors la pour la donnée d'un faisceau quasi-
cohérent sur $\mathfrak X$ ou d'un morphisme dans $QCoh(\mathfrak X)$, on peut se restreindre 
au données de la définition, mais seulement sur $U,U',U''$ des $S$-schémas lisses (avec les 
mêmes compatibilités).
\elem

\dem
Soit donc $S$ un schéma. Supposons donc que l'on a un champ algébrique 
$\mathfrak X \fleche S$ lisse. C'est-à-dire qu'il existe une présentation (lisse),
\[ p : X \fleche \mathfrak X,\]
avec $X$ lisse sur $S$.
Supposons maintenant donné pour tout schéma lisse $U$ sur $S$, et tout 
$u \dans \mathfrak X_U$, la donné d'un faisceau quasi-cohérent $\mathcal F_u$ sur $U$ 
(avec les compatibilités de la définition \ref{defqcoh} lorsque $U',U''$ lisses sur $S$). 
Alors comme $X$ est lisse, on a la donnée d'un faisceau quasi-cohérent $F := \mathcal F_p$ 
sur $X$, et muni d'une donnée de descente de la même manière que dans la 
proposition \ref{proqcoh} 
(car $X \times_{\mathfrak X}X$ et $X\times_{\mathfrak X}X\times_{\mathfrak X}X$ sont 
des schémas lisses sur $S$). 
On en déduit donc que $F$ est dans $Des(X)$, et donc qu'il existe un 
unique faisceau quasi-cohérent sur $\mathfrak X$, $\widetilde{\mathcal F}$, 
qui lui correspond. 
Par construction, 
$\mathcal F_u \simeq (\widetilde{\mathcal F})_u$ (canoniquement) 
pour tout $u : U \fleche \mathfrak X$ avec $U$ lisse sur $S$.
Soit maintenant pour tout $U/S$ lisse, et tout $u \dans \mathfrak X_U$, des faisceaux 
quasi-cohérents $\mathcal F_u$ et $\mathcal G_u$ sur $U$, et d'un morphisme,
\[\phi_u : \mathcal F_u \fleche \mathcal G_u,\]
avec les compatibités de la définition \ref{defqcoh}, pour tout $U',U''$ lisses sur $S$.
Alors comme $X$ est lisse, il existe un morphisme,
\[ \phi_p : F \fleche G, \quad \text{où } F = \mathcal F_p, G = \mathcal G_p.\]
On a vu que $F, G$ sont dans $Des(X)$ et de plus, comme dans la proposition \ref{proqcoh}, 
$\phi_p$ est un morphisme dans $Des(X)$. On en déduit donc qu'il existe un unique morphisme,
\[\widetilde{\phi} :\widetilde{\mathcal F} \fleche \widetilde{\mathcal G}\]
dans $QCoh(\mathfrak X)$ tel que le diagramme suivant commute,
\begin{center}
\begin{tikzpicture}[description/.style={fill=white,inner sep=2pt}] 
\matrix (m) [matrix of math nodes, row sep=3em, column sep=2.5em, text height=1.5ex, text depth=0.25ex] at (0,0)
{ 
\widetilde{\mathcal F}_u  &  & \widetilde{\mathcal G}_u \\
\mathcal F_u  &  & \mathcal G_u  \\
 };

\path[->,font=\scriptsize] 
(m-1-1) edge node[auto] {$\widetilde{\phi}_u$} (m-1-3)
(m-1-1) edge node[auto] {$\simeq$} (m-2-1)
(m-1-3) edge node[auto,left] {$\simeq$} (m-2-3)
(m-2-1) edge node[auto] {$\phi_u$} (m-2-3)
;
\end{tikzpicture}
\end{center}
pour tout $u \dans \mathfrak X_U$ avec $U$ lisse sur $S$.
\edem

On utilise ce dernier lemme pour construire les invariants de Hasse partiel sur le champ des 
groupes de Barsotti-Tate tronqués d'échelon $r$, munis d'une action de $\mathcal O_F$, car ces derniers
 ne sont définis (à priori) que sur une base lisse. 
 
\section{Diviseurs de Cartier et Champs}
\label{appB}
Soit $\mathfrak X$ un champ sur $S$, suivant \cite{LMB}, on définit son espace topologique sous-jacent,
\[ |\mathfrak X| = \left(\coprod_{K corps} ob\, \mathfrak X_{\Spec K}\right)_{\big{/\sim}},\]
où $u \sim v$, $u \dans \mathfrak X_{\Spec K}, v \dans \mathfrak X_{\Spec(K')}$ si il existe une 
extension commune $K''$ de $K,K'$, $i : K \subset K'', j : K' \subset K''$ et un isomorphisme de 
$i^*u$ sur $j^*v$ dans $\mathfrak X_{\Spec K''}$. On munit $|\mathfrak X|$ de la topologie de Zariski, donnée par les $|\mathfrak U|$, où $\mathfrak U$ est un sous-champ ouvert.

\defi
\label{deficart}
Soit $S$ un schéma, et $\mathfrak X$ un champ algébrique sur $Sch/S$. 
Soit $D \subset \mathfrak X$ 
un sous-champ fermé. Alors les conditions suivantes sont équivalentes (voir exemple 
\ref{exeinv}):
\begin{enumerate}
\item \label{1}$\mathcal I_D$ est un $\mathcal O_{\mathfrak X}$-module inversible.
\item \label{2} Il existe $X \overset{p}{\fleche} \mathfrak X$, avec $X$ un schéma, surjectif et lisse tel que 
$D\times_{\mathfrak X} X$ est un diviseur de Cartier sur $X$.
\item \label{3} Pour tout $X \overset{p}{\fleche} \mathfrak X$, avec $X$ un schéma, surjectif et lisse, $D\times_{\mathfrak X} X$ est un diviseur de Cartier sur $X$.
\item \label{4} Il existe un recouvrement affine lisse quasi-compact i.e. pour tout $x \dans \mathfrak X$,  des voisinages lisses $U \fleche \mathfrak X$ de $x$, tel que si $U = \Spec A$, alors $D\times_{\mathfrak X} U= \Spec(A/(f))$ et $f \dans A$ n'est 
pas un diviseur de zéro. 
\end{enumerate}
Dans ce cas, on dit que $D$ est un \textit{Diviseur de Carter effectif} de $\mathfrak X$.
\edefi

\dem
\ref{1} $\Rightarrow$ \ref{3} : On a donc $\mathcal I_D$ un faisceau inversible, donc $p^*\mathcal I_D$ est un faisceau inversible sur $X$. Donc $D\times_{\mathfrak X} X$ est un diviseur de Cartier sur $X$, car $p$ est lisse donc plat.
\\
\ref{2} $\Rightarrow$ \ref{1} : On a donc sur $X$ l'existence d'un $(\mathcal I_D)_p \subset \mathcal O_X$ 
qui est un faisceau d'idéaux inversible associé à $D\times_{\mathfrak X} X$, munie d'une donnée 
de descente par la proposition $\ref{proqcoh}$. De plus, par la remarque \ref{remqcoh}, $\mathcal I_D$ est inversible.
\\
\ref{2} $\Rightarrow$ \ref{4} : Soit $X \fleche \mathfrak X$ lisse surjectif une présentation par un 
schéma. Prenons un ouvert de $X$,  $U, \fleche X$, avec $U = \Spec A$ affine.
Alors,
\[ D \times_{\mathfrak X} U = (D\times_{\mathfrak X} X) \times_X U = \Spec(A/(f)),\]
avec $(f)$ non diviseur de zéro, car $D\times_{\mathfrak X} X$ est un diviseur de cartier sur $X$. De plus, comme $X \fleche \mathfrak X$ est lisse, donc quasi-compact, le recouvrement précédent est quasi-compact.\\
\ref{4} $\Rightarrow$ \ref{2} : On écrit un morphisme surjectif, et localement lisse, donc plat quasi-compact,
\[ U = \coprod_{x \in |\mathfrak X|} U_x \fleche \mathfrak X.\]
Ce morphisme n'est pas lisse a priori, mais quitte à prendre une présentation $p: X \fleche \mathfrak X$, $\coprod_{x \in |\mathfrak X|} (U_x \times_{\mathfrak X} X)$ est un recouvrement localement lisse de $X$, donc on peut extraire un recouvrement $U'$ quasi-compact, puisque $X$ est quasi-séparé. Maintenant le morphisme
\[ U' \fleche \mathfrak X,\]
est quasi-compact, plat et surjectif, et $D \times_{\mathfrak X} U'$ est localement un diviseur de Cartier (sur chaque $U_x$) par hypothèse, $D\times_\mathfrak X X$ est donc un diviseur de Cartier sur $X$.
\edem

\rem
Soit $\mathfrak X$ un champ algébrique, $X \overset{p}{\fleche} \mathfrak X$ une présentation. Alors on a une surjection d'espaces topologiques,
\[ |X| \overset{|p|}{\fleche} |\mathfrak X|.\]
D'après \cite{LMB}, Proposition 5.6, $|p|$ est un morphisme ouvert. Plus précisément pour tout 1-morphisme de champs, plat localement de présentation finie, $f$, $|f|$ est ouvert.
\erem

\pro
Soit $\mathfrak X$ un champ algébrique. On a une équivalence de catégories,
\begin{IEEEeqnarray*}{ccc}
\left\{ \begin{array}{c}
\text{Diviseurs de cartier effectifs} \\
\text{de } \mathfrak X
 \end{array}
\right\} & \overset{\sim}{\fleche} & 
\left\{
\begin{array}{c}
s : \mathcal O_{\mathfrak X} \fleche \mathcal L \text{ où} \\
\mathcal L \text{ est un faisceau inversible et } \\
$s$ \text{ génériquement trivialisante.}
 \end{array}
\right\} 
\\
D & \longmapsto & (\mathcal O_{\mathfrak X} \fleche \mathcal I_D^\vee)\\
V(s) \subset \mathfrak X & \longmapsfrom & (s : \mathcal O_X \fleche \mathcal L)
\end{IEEEeqnarray*}
où $V(s)$ est un sous-schéma fermé de $\mathfrak X$ donné par,
\[ V(s)_U = \{ u \dans \mathfrak X_U : u^*s = 0\}. \]
\epro

\dem
C'est vrai sur un schéma (\cite{Stackproj} DIVISORS lemma 2.20). On se ramène à ce cas en considérant une présentation : il suffit de combiner la définition précédente avec la proposition \ref{proqcoh} 
(voir aussi remarque \ref{remqcoh}).
\edem

\pro
Soit $\mathfrak X$ un champ algébrique réduit (i.e. tel qu'il existe une présentation 
$X \overset{p}{\fleche}\mathfrak X$ avec $X$ un schéma réduit).
Soit $V(s)$ le sous-champ fermé de $\mathfrak X$ associé à un couple $(\mathcal L,s)$ où $s : \mathcal O_{\mathfrak X} \fleche \mathcal L$ est une section. 
Notons $D(s)$ le sous-champ ouvert de $\mathfrak X$ définipar,
\[D(s)_U = \{ u \dans \mathfrak X_U : u^*s :\mathcal O_U \overset{\sim}{\fleche} \mathcal L_u\}. \]
Alors $V(s)$ est un diviseur de Cartier si et seulement si $|D(s)|$ est dense dans $|\mathfrak X|$.
\epro

\dem
Soit $p : X \fleche \mathfrak X$ une présentation (réduite) de $\mathfrak X$. 
Supposons que $V(s)$ soit un diviseur de Cartier, alors $p^*V(s)$ est un diviseur de Cartier sur 
$X$, donc $X\priv p^*V(s)$ est un ouvert dense (car $X$ réduit) et comme $|p|$ est ouverte et 
surjective, $|p(X\priv p^*V(s))| = |D(s)|$ (car $p$ surjectif) est un ouvert dense.\\
Réciproquement, si $|D(s)|$ est dense, $|p|^{-1}(|D(s)|) = |X \priv p^*V(s)|$ est un ouvert dense de $X$, donc $p^*V(s)$ est un diviseur de Cartier sur $X$, et la définition \ref{deficart} conclut.
\edem

\rem
Être réduit pour un champ ne dépend pas de la présentation : Si $(X,\pi_1),(X',\pi_2)$ sont deux présentations de $\mathfrak X$ et $X$ est réduit, alors $X \times_\mathfrak 
X X' = \pi_2^*X$ est un schéma (car $\pi_1$ relativement représentable) et $X \times_{\mathfrak X} X'$ est réduit car $X$ l'est et $\pi_2$ est lisse (donc $\pi_1^*\pi_2$ 
aussi). Donc comme $\pi_2^*\pi_1 : X \times_{\mathfrak X}X' \fleche X'$ est lisse, $X'$ est réduit. 
\erem

 \nocite{*}
\bibliographystyle{alpha-fr} 
\bibliography{hasseinvariantsv0.6} 

\backmatter

\end{document}